%% file: 2-epsilon_exp.tex
\documentclass{amsart}

\usepackage{amssymb}
\usepackage{tikz}
\usepackage{verbatim}
\usepackage{hyperref}
\usepackage{orcidlink}
\usepackage[english]{babel}
\usepackage{wasysym}
\usepackage[shortlabels]{enumitem}

\usepackage[colorinlistoftodos]{todonotes}
\usepackage{color, soul}
\usepackage{mathtools}
\usepackage{mathrsfs}

\usepackage{thmtools}
\usepackage{thm-restate}

\newtheorem{theorem}{Theorem}[section]
\newtheorem{lemma}[theorem]{Lemma}
\newtheorem{corollary}[theorem]{Corollary}

\theoremstyle{definition}
\newtheorem{definition}[theorem]{Definition}
\newtheorem{notation}[theorem]{Notation}
\newtheorem{example}[theorem]{Example}

\newtheorem{proposition}[theorem]{Proposition}

\newtheorem{facts}[theorem]{Facts}
\newtheorem{observation}[theorem]{Observation}

\newtheorem{conjecture}[theorem]{Conjecture}

\theoremstyle{remark}
\newtheorem{remark}[theorem]{Remark}
\theoremstyle{remark}

\numberwithin{equation}{section}

\newcommand{\aut}{\operatorname{Aut}}

\newcommand{\sym}{\operatorname{Sym}}

\newcommand{\agl}{\operatorname{AGL}}
\newcommand{\pgl}{\operatorname{PGL}}
\newcommand{\pgammal}{\operatorname{P\Gamma L}}
\newcommand{\id}{\operatorname{id}}

\newcommand{\csp}{\operatorname{CSP}}
\newcommand{\dom}{\operatorname{Dom}}

\newcommand{\rk}{\operatorname{rk}}

\newcommand{\inter}{\operatorname{I}}
\def\acts{\curvearrowright}

\newcommand{\acl}{\operatorname{acl}}

\newcommand{\fa}{\mathfrak{A}}
\newcommand{\fb}{\mathfrak{B}}
\newcommand{\fc}{\mathfrak{C}}

\newcommand{\ff}{\mathfrak{F}}

\newcommand{\fm}{\mathfrak{M}}
\newcommand{\fn}{\mathfrak{N}}
\newcommand{\fv}{\mathfrak{V}}
\newcommand{\sign}{\Sigma}

\newcommand{\bigeq}{\nabla}
\newcommand{\smalleq}{\Delta}
\newcommand{\partition}{P}

\newcommand{\gone}{\mathcal{H}}
\newcommand{\gones}{\mathscr{H}}

\newcommand{\biggroup}{\mathcal{G}}
\newcommand{\kernel}[2]{\mathcal{K}^{#1}_{#2}}

\newcommand{\fbh}{\mathscr{B}}

\newcommand{\Ker}{\operatorname{Ker}}

\newcommand{\same}{\approx}
\newcommand{\notsame}{\not\approx}
\newcommand{\age}{\operatorname{Age}}
\newcommand{\forb}{\operatorname{Forb}}

\newcommand{\subex}{\mathcal{S}}
\newcommand{\subexs}{\mathscr{S}}
\newcommand{\fh}{\mathcal{F}}
\newcommand{\fhs}{\mathscr{F}}
\newcommand{\group}{\mathcal{G}}
\newcommand{\fin}{\operatorname{fin}}

\newcommand{\betw}{\operatorname{Betw}}
\newcommand{\cyc}{\operatorname{Cyc}}
\newcommand{\sep}{\operatorname{Sep}}

\newcommand{\orbits}{\mathcal{O}}
	\newcommand{\orb}{o}
\newcommand{\oi}{\ell}
\newcommand{\os}{u}

\newcommand{\flip}{\leftrightarrow}
\newcommand{\turn}{\circlearrowleft}

\newcommand{\data}{D}
\newcommand{\lift}{\mathcal{L}}

\newcommand{\unsim}{\mathord{\sim}}

\begin{document}

\title{Structures with not too fast unlabelled growth}
\author{Bertalan Bodor \orcidlink{0009-0003-6679-6355}}
\address{HUN-REN Alfréd Rényi Institute of Mathematics}
\email{bodor@renyi.hu}
\thanks{The author has been funded by the European Research Council (Project POCOCOP, ERC Synergy Grant 101071674). Views and opinions expressed are however those of the authors only and do not necessarily reflect those of the European Union or the European Research Council Executive Agency. Neither the European Union nor the granting authority can be held responsible for them.}
\maketitle

\begin{abstract}
	Let $\subexs$ be the class of all structures whose growth rate on orbits of subsets of size $n$ is not faster than $\frac{2^n}{p(n)}$ for any polynomial $p$. In this article we give a complete classification of all structures in $\subexs$ in terms of their automorphism groups. As a consequence of our classification we show that $\subexs$ has only countably many structures up to bidefinability, all these structures are first-order interpretable in $(\mathbb{Q};<)$ and they are interdefinable with a finitely bounded homogeneous structure. Furthermore, we also show that all structures in $\subexs$ have finitely many first-order reduct up to interdefinability, thereby confirming Thomas' conjecture for the class $\subexs$.
\end{abstract}

\input intro.tex
\input prelims.tex
\input fincover_q.tex
\input fincover_h.tex

\input finish.tex
\input thomas.tex

\bibliographystyle{alpha}
\bibliography{local.bib}

\end{document}

%% file: intro.tex

\section{Introduction}

	There are three natural sequences of counting orbits attached to a permutation group $G\leq \sym(X)$.

We write
	
\begin{itemize}
\item $\orb_n(G)$ for the \emph{number of $n$-orbits of $G$}, i.e., the number of orbits of the coordinate-wise action of $G$ on $X^n$
\item $\oi_n(G)$ for the \emph{number of injective $n$-orbits of $G$}, i.e., the number of those $n$-orbits of $G$ which contain tuples with pairwise different entries,
\item $\os_n(G)$ for the \emph{number of orbits of $n$-subsets of $G$}, i.e., the number of orbits of the natural action of $G$ on $n$-element subsets of $X$.
\end{itemize}

	We call the sequences $(\oi_n(G))_{n\in \omega}$ or $(\os_n(G))_{n\in \omega}$, the \emph{labelled} and \emph{unlabelled growth} of $G$, respectively.

	In this paper we always assume that the group $G$ is \emph{oligomorphic}, i.e., any (or equivalently all) of the sequences defined above have finite values. In this case we can assume without loss of generality that the domain set $G$ is closed in the topology of pointwise convergence, since switching to the closure of the group does not change its orbit growth functions. By a Theorem of Engeler, Ryll-Nardzeski, and Svenonius we know that this is equivalent to $G$ being an automorphism group of an $\omega$-categorical structures~\cite{HodgesLong}. We also know that in this case every orbit of $G$ is a realization of a complete type in $\fa$~\cite{HodgesLong}. By the Downward L\"{o}wenheim--Skolem theorem we can (and we will) also assume that the underlying set $X$ is countable (see also the discussion in~\cite{Oligo}, Section 1). We use the notation $\orb_n,\oi_n,\os_n$ and the corresponding terminology for the group and the corresponding structure $\fa$ interchangeably. The study of the behavior of the sequences above for oligomorphic groups was initiated by Cameron and Macpherson~\cite{CameronCounting,Oligo}, and it has been a subject of active research since.
	
	There are many structural results available about $G$, or the corresponding structure $\fa$ as above, when we put more restrictive upper bounds on the sequences $(\oi_n(G))_{n\in \omega}$ or $(\os_n(G))_{n\in \omega}$. Clearly, a structure is finite if and only if any of these sequences are eventually 0. For this reason one can argue that the slower the labelled or the unlabelled growth of a structure is, the closer it is to a finite structure. This intuition is also confirmed by some the concrete known results we mention below.

	As for $(\oi_n(G))$ we know that it has an exponential upper bound if and only if $G$ is a closed supergroup of a finite direct product of symmetric groups (with the action on the disjoint union)~\cite{bodirsky2021permutation}. In fact, it is shown in the same paper that $\oi_n(G)<cn^{dn}$ for some $d<1$ if and only if $\fa$ is \emph{cellular}, i.e., it is a first-order reduct of a finite cover of a unary structure with finitely many predicates (see more details in~\cite{braunfeld2022monadic,bodor2024classification} or in Subsection~\ref{sect:her}).
	
	Compared to the labelled growth the situation tends to be much more complicated if we put similar upper bounds for the unlabelled growth of a group. Already the case when $\os_n(G)=1$ for all $n$, these groups are called \emph{highly homogeneous}, is far from obvious, and these groups were classified by Cameron: they are exactly the automorphism groups of the 5 reducts of $(\mathbb{Q},<)$ (see Subsection~\ref{sect:high}). Those closed groups for which $\os_n(G)$ has a polynomial upper bound are classified in~\cite{falque2020classification} (or in~\cite{falque2019classification}). It turns out, as it is shown in the same paper, that in this case the sequence $(\os_n(G))$ is in fact always \emph{quasi-polynomial}, i.e., it can be written as $$\os_n(G)=a_k(n)n^k+\dots+a_0(n)$$ where all the coefficients $a_i(n)$ are periodic functions.
	
	In the case when the structure $\fa$ is stable an interesting dichotomy holds.
	
\begin{theorem}[\cite{braunfeld2022monadic}, Theorem 3.10]\label{sam_intro}
	Let $\fa$ be a stable $\omega$-categorical structure. Then exactly one of the following holds.
\begin{enumerate}
\item\label{it:monad} $\os_n(\fa)<c^n$ for all $c>1$ if $n$ is large enough.
\item $\os_n(\fa)>c^n$ for all $c$ if $n$ is large enough.
\end{enumerate}
\end{theorem}

	This means that for automorphism groups of stable structures there is no difference between `slower than exponential' and `at most exponential' unlabelled growth. We also know a structure description of those structures which satisfy item~\ref{it:monad} in the above theorem. A structure is called \emph{hereditarily cellular} iff it can be built up from finite structures by taking finite disjoint unions, wreath products with the pure set, and first-order reducts (see more details in Subsection~\ref{sect:her}). Hereditarily cellular structures are introduced, under a different terminology in~\cite{lachlan1992}, and their automorphism groups are discussed in detail in~\cite{bodor2024classification}. By the results of~\cite{braunfeld2022monadic} we also know that a stable $\omega$-categorical structure has slower than exponential unlabelled growth (item~\ref{it:monad} in the theorem above) if and only if it is hereditarily cellular. From this we can conclude the following.
	
\begin{theorem}\label{intro:ms}
	Let $\fa$ be a countable structure. Then the following are equivalent.
\begin{enumerate}
\item $\fa$ is hereditarily cellular.
\item $\fa$ is stable and $\os_n(\fa)<c^n$ for all $c>1$ if $n$ is large enough.
\item $\fa$ is stable and $\os_n(\fa)<c^n$ for some $c$ if $n$ is large enough.
\end{enumerate}
\end{theorem}
	
	If we remove the stability assumption then the conclusion of Theorem~\ref{sam_intro} no longer holds. The simplest counterexample is the wreath product $Z_2\wr \aut(\mathbb{Q})$ whose unlabelled growth is the Fibonacci sequence. In general, it seems that the larger $c$ is in an upper bound $c^n$ for the unlabelled growth, the more complicated the situation gets. 

	In this paper we focus our attention to the case when $c<2$. We define $\subexs$ to be the class of those countable structures $\fa$ such that for no polynomial $p$ do we have $\os_n(\fa)\geq \frac{2^n}{p(n)}$. The class $\subexs$ was introduced and studied in a recent paper by Simon~\cite{simon2025omega} (where these were called structures with ``few finite substructures''). The following is shown in the aforementioned paper.

\begin{theorem}\label{simon_intro}
	Let $\fa\in \subexs$ be a structure with a primitive automorphism group. Then $\os_n(\fa)=1$ for all $n$, (and thus $\aut(\fa)$ is one of the 5 groups in Cameron's classification).
\end{theorem}
	
	We note that in the original formulation of Theorem~\ref{simon_intro} in~\cite{simon2025omega} also allows $\fa$ to be strictly stable (stable but not $\omega$-stable) in the conclusion, however this case can be excluded by an argument presented in~\cite{braunfeld2022monadic}. Theorem~\ref{simon_intro} confirms a long-standing conjecture by Macpherson. We also note that the number 2, hidden in the definition of $\subexs$, in the assumption of Theorem~\ref{simon_intro} is optimal. This is witnessed by the structure $S(2)$, called the \emph{local order}, whose automorphism is primitive and its unlabelled growth is approximately $2^n/n$.
	
	Finally, we mention a recent structural result for general exponential upper bound. A first-order $\phi(x_1,\dots,x_m,y_1,\dots,y_n)$ formula is said to have the \emph{independence property} if we can find an infinite set $A$ of $m$-tuples and a family $\bar{b}_I: I\subset A$ of $n$-tuples such that $\phi(\bar{a},\bar{b}_I)$ holds if and only if $\bar{a}\in I$. A structure $\fa$ is called \emph{NIP} if there is no formula in any elementary extension of $\fa$ that has the independence property. A structure is called \emph{monadically NIP} if all of its unary expansions are NIP.
	
\begin{theorem}[\cite{braunfeld2021characterizations}, Theorem 5.3]\label{mnip_intro}
	Let $\fa$ be a structure such that $\os_n(\fa)<c^n$ for some $c$. Then $\fa$ is monadically NIP.
\end{theorem}

	In fact, it is conjectured that the converse of Theorem~\ref{mnip_intro} also holds for $\omega$-categorical structures (\cite{braunfeld2022monadic}, Conjecture 2).

\subsection{The main result}

	In the present paper we give a complete classification of the class $\subexs$ in terms of their automorphism groups. This is a common generalization of the classification results given in~\cite{falque2020classification} and~\cite{bodor2024classification}, and in particular contains the classification of all closed permutation groups with subexponential unlabelled growth. Our classification results are based on the description given~\cite{simon2025omega} and the description of hereditarily cellular structures/groups as given in~\cite{bodor2024classification}. We also obtain a description of structures in $\subexs$ similar to that of hereditarily cellular structures we mentioned earlier: a structure is contained in $\subexs$ if and only if it can be built up from certain finite covers of $\mathbb{Q}$ and finite structures by taking finite disjoint unions, wreath products with the pure set, and first-order reducts. The full details are provided in Section~\ref{sect:finish}.
	
	As a consequence of our classification we verify two conjectures made in~\cite{braunfeld2022monadic}: one about a description of structures with subexponential orbit growth (Conjecture 1, see Theorem~\ref{final_struct}), and one about gaps in unlabelled growths in the class $\subexs$ (Conjecture 3). The latter can be formulated as follows.
	
\begin{restatable}{theorem}{thmgaps}\label{gaps_intro}
	Let $\fa\in \subexs$. Then there exists some $d\in \omega, d\geq 1$ such that $\lim_{n\rightarrow \infty}(u_n(\fa))^{1/n}=\gamma_d$ where $\gamma_d$ denotes the largest real root of the polynomial $x^d-x^{d-1}-\dots-1$.
\end{restatable}

	We also show that $\subexs$ contains only countably many structures up to bidefinability.

	Our classification also gives us the following the characterization of structures with at most polynomial unlabelled growth.
	
\begin{restatable}{theorem}{polyconst}\label{poly_const_struct}
	Let $\fa$ be a countable structure. Then the following are equivalent.
\begin{enumerate}
\item $\os_n(\fa)$ is at most polynomial.
\item There exists an expansion of $\fa$ by finitely many constants which is bidefinable with a disjoint union of finitely many copies of $(\mathbb{Q};<)$ and some cellular structure.
\end{enumerate}
\end{restatable}

	Theorem~\ref{poly_const_struct} is especially relevant in the context of infinite-domain CSPs; see our discussion below.

\subsection{Model theoretical consequences}\label{sect:model}

\subsubsection*{Finite homogeneity and finite boundedness}

    We say that a structure $\fa$ is \emph{homogeneous} if every partial isomorphism between substructures of $\fa$ extends to an automorphism of $\fa$. It follows from an easy counting argument that if $\fa$ is homogeneous and has a finite relational signature then $\orb_n(\fa)\leq 2^{cn^k}$ for some $c>0$ where $k$ is the maximum of 2 and the maximal arity of $\fa$, in particular $\fa$ is $\omega$-categorical. We say that a structure $\fa$ is \emph{finitely homogenizable} if $\fa$ is first-order interdefinable with a homogeneous structure in a finite relational signature. Note that the upper bound above does not guarantee finite homogenizability: for instance the infinite dimensional vector space $\fv$ over a finite field $K$ we have $\orb_n(\fv)\leq {|K|^n\choose n}< |K|^{n^2}$ but $\fv$ is not finitely homogenizable~\cite{BodorCameronSzabo,macpherson1991interpreting}. In the present paper, however, we show that every structure in $\subexs$ is finitely homogenizable. In fact, we show that every structures in $\subexs$ is interdefinable with a homogeneous relational structure $\fa$ which is also \emph{finitely bounded}, meaning that the class of finite structures that embed into $\fa$ can be described by finitely many forbidden substructures.

\subsubsection*{Interpretability in $(\mathbb{Q};<)$}

    If $\fa$ and $\fb$ are structures then we say that $\fb$ is (first-order) interpretable in $\fa$ if there exist a $d\in \omega$ and a partial surjective map $I\colon A^d\rightarrow B$ such that for all relations defined by atomic formulas in $\fb$, including $x=y$ and $x=x$, are definable in $\fa$. Our classification shows that every structure in $\subexs$ is interpretable in $(\mathbb{Q};<)$. We also note that for this result the number 2 in the definition of $\subexs$ is tight as it is shown in~\cite{bodirsky2025taking} that $S(2)$ is not interpretable in $(\mathbb{Q};<)$.

\subsection{Thomas' conjecture}

	In~\cite{RandomReducts} Thomas made the conjecture that every countable homogeneous structure with a finite relational signature has finitely many reducts up to first-order interdefinability. An equivalent formulation of the conjecture states that automorphism groups of all aforementioned structures have finitely many closed supergroups contained in the symmetric group acting on the domain set of the structure. The conjecture has been verified for many well-known homogeneous structures. The list includes the rationals with the usual ordering~\cite{Cameron5}, the countably infinite random graph~\cite{RandomReducts}, the homogeneous universal $K_n$-free graphs~\cite{Thomas96}, the expansion of $({\mathbb Q};<)$ by a single constant~\cite{JunkerZiegler}, the universal homogeneous partial order~\cite{Poset-Reducts}, the random ordered graph~\cite{42}, and many more~\cite{agarwal,AgarwalKompatscher,BodJonsPham,BBPP18}. Note that if we drop the assumption that the signature of the homogeneous structure $\fa$ is relational, then Thomas' conjecture is false even if we keep the assumption that $\fa$ is $\omega$-categorical: already the countable atomless Boolean algebra has infinitely many first-order reducts~\cite{BodorCameronSzabo}.
	
	Note that all aforementioned results only solve Thomas' conjecture for specific structures mostly using some general techniques from Ramsey theory, and they all provide the exact number of reducts of such structures. Another way to approach this conjecture is to prove (or potentially disprove) it for more restrictive classes of structures, but potentially without specifying the exact number of reducts. We mention two results in this direction. In~\cite{bodor2024classification} the conjecture has been verified for hereditarily cellular structures, and in~\cite{simon2022nip} it is verified for $\omega$-categorical primitive NIP structures with thorn rank 1. 
	
	In this paper we show that Thomas' conjecture also holds for every structure in the class $\subexs$. Since this class contains all hereditarily cellular structures and all expansion of $(\mathbb{Q};<)$ with finitely many constants, our result generalizes the results of both~\cite{bodor2024classification} and~\cite{JunkerZiegler} with regard to Thomas' conjecture.

\subsection{Connection to CSPs}

	Homogeneous and in general $\omega$-categorical structures also play an important role in the study of \emph{Constraint Satisfaction Problems}, or \emph{CSPs} for short. These are usually formulated as follows. The CSP over a structure $\fb$ with a finite relational signature, denoted by $\csp(\fb)$, is the decision problem where the input is a finite structure $\fa$ with the same signature as $\fb$, and we need to decide whether there exists a homomorphism from $\fa$ to $\fb$. By the famous result of Bulatov of Zhuk, we know that if $\fb$ is finite then the computational complexity of $\csp(\fb)$ satisfies a dichotomy: it is either in $\mathbf{P}$ or it is $\mathbf{NP}$-complete (\cite{BulatovFVConjecture,zhuk2020proof}). Both proofs use various techniques and results from universal algebra many of which can also be generalized to $\omega$-categorical structures, as well. In order to formulate a complexity dichotomy we need further model-theoretical restriction since in general even homogeneous structures with a finite relational signature can have undecidable CSPs (see for instance~\cite{BodirskyNesetrilJLC}). One way to guarantee that $\csp(\fb)$ is at least in $\mathbf{NP}$ is to assume that $\fb$ is homogeneous and finitely bounded, or a reduct of such a structure. This motivates the following complexity dichotomy conjecture, originally presented in~\cite{bodirsky2021projective}.

\begin{conjecture}[Infinite-domain CSP dichotomy conjecture]\label{conj:csp}
	Let $\fb$ be a first-order reduct of a finitely bounded homogeneous relational signature. Then $\csp(\fb)$ is either in $\mathbf{P}$ or it is $\mathbf{NP}$-complete.
\end{conjecture}

	For a general introduction and more details on infinite domain CSPs we refer the reader to~\cite{bodirsky2021complexity,pinsker2022current}.

	A natural subclass of structures for which the above conjecture is sometimes considered is the class of structures which are interpretable in $(\mathbb{Q};<)$ (see for instance~\cite{Bodirsky}, partial results in this direction have been obtained in~\cite{ecsps,tcsps-journal,BodMot-Unary,bodor2022csp,bodirsky2024complexity,bodirsky2025structures}. The conjecture has also been verified for various other classes of structures, see~\cite{BodPin-Schaefer-both,posetCSP18,bodirsky2021canonical,Phylo-Complexity,bodirsky2018universal,bodor2022csp,mottet2024order,mottet2024smooth}.
	 
	As we mentioned in Subsection~\ref{sect:model} we know that all structures in $\subexs$ are interpretable in $(\mathbb{Q};<)$. In~\cite{bodor2022csp} Conjecture~\ref{conj:csp} has been verified for hereditarily cellular structures which are exactly the stable structures in $\subexs$ by Theorem~\ref{intro:ms}.
	
	A more approachable goal would be to verify Conjecture~\ref{conj:csp} for structures with at most polynomial unlabelled growth. By using Theorem~\ref{poly_const_struct} one can argue that this is enough to show it in the case when the structure $\fb$ is bidefinable with a disjoint union of finitely many copies of $(\mathbb{Q};<)$ and some cellular structure (see the discussion at the end of Section~\ref{sect:orbit}). This looks promising since the infinite-domain CSP dichotomy conjecture has already been solved for all reducts of $(\mathbb{Q};<)$~\cite{tcsps-journal} and all hereditarily cellular structures~\cite{bodor2022csp}. We want to emphasize, however, that in order to recover the complexity of a CSP of structure usually a much finer distinction than just bidefinability is necessary, namely we need to identify relations which are \emph{primitively positively} definable. Therefore, even in the simple-looking case outlined above the solution to Conjecture~\ref{conj:csp} poses some serious challenges.

%% file: prelims.tex

\section{Preliminaries}\label{sect:prelims}

\subsection{Permutation group notation}
	
	For a set $X$, we write $\sym(X)$ for the group of all permutations of $X$. We say that a group $G$ is a \emph{permutation group} $G\leq \sym(X)$ for some set $X$. Note that in this case the set $X$ is unique, and we call it the \emph{domain} of $G$, denoted by $\dom(G)$. The cardinality of $\dom(G)$ is called the \emph{degree} of $G$. We say that two permutation groups $G$ and $H$ are \emph{isomorphic as permutation groups} if there is bijection $e\colon \dom(X)\rightarrow \dom(H)$ such that $H=eGe^{-1}=\{e^{-1}ge: g\in G\}$. All groups appearing in this paper will be permutation groups, and when we say that two permutation groups are isomorphic we always mean that they are isomorphic \emph{as permutation groups}.

	Assume that $G\leq \sym(X)$ is a permutation group, and $Y\subset X$. We introduce the following notation.

\begin{itemize}[noitemsep,topsep=0pt]
\item $G_{Y}$ denotes the \emph{pointwise stabilizer} of $Y$, that is, $$G_Y=\{g\in G\colon\forall y\in Y(g(y)=y)\}.$$
\item $G_{\{Y\}}$ denotes the \emph{setwise stabilizer} of $Y$, that is, $$G_{\{Y\}}=\{g\in G\colon\forall y\in Y(g(y)\in Y)\}.$$
\item $G|_Y$ denotes the \emph{restriction of $G$ to $Y$}, that is, $$G|_Y=\{h|_Y\colon h\in G\},$$ provided that $Y$ is invariant under $G$.
\item $G_{(Y)}:=G_{\{Y\}}|_Y=G_{\{X\setminus Y\}}|_Y$.
\item $G_{((Y))}:=G_{\dom(G)\setminus Y}|_Y$.
\end{itemize}
	If $Y$ is finite, say $Y=\{x_1,\dots,x_n\}$, then we simply $G_{x_1,\dots,x_n}$ for $G_Y$.

	Every permutation group $G\leq \sym(X)$ is equipped with a topology called the \emph{topology of pointwise convergence}. This topology is defined to be the subspace topology of the product space $X^X$ where $X$ is equipped with the discrete topology. We say that a permutation group is \emph{closed} if it is closed in the topology of pointwise convergence. Equivalently, a group $G\leq \sym(X)$ is closed if the following holds: for all $g\in\sym(X)$, if for every finite $F\subset X$ there exists $g'\in G$ such that $g'|_F=g|_F$ then $g\in G$.
	
	Let us assume that $E$ is an equivalence relation on some set $X$. Then for an element $x\in X$ we write $[x]_E$ for the unique $E$ class that contains $E$. Note that in this case $[x]_E=[y]_E$ holds if and only if $(x,y)\in E$. We denote by $X/E$ the set of $E$ classes on $X$, i.e., $X/E=\{[x]_E: x\in X\}$. If $G\leq \sym(X)$ is a group then we say that $E$ is a \emph{congruence} of $G$ if $E$ is preserved by every element of $G$, that is, for all $x,y\in X$ and $g\in G$ we have $(x,y)\in E$ if and only if $(g(x),g(y))\in E$. If $E$ is a congruence of $G$ then there is a natural homomorphism $G\rightarrow \sym(X/E)$ defined by $g\mapsto g/E$ where $g/E: [x]_E\mapsto [g(x)]_E$. We write $G/E$ for the image of this homomorphism, i.e., $G/E=\{g/E: g\in G\}$.

\subsection{Orbit growth functions}

	Let $G\leq \sym(X)$ be a permutation group, and let $n\in \omega$. Then we introduce the following notation.
\begin{itemize}
\item $\orbits_n(G)$ denotes the \emph{set of $n$-orbits of $G$}, i.e., the set of orbits of the coordinate-wise action of $G$ on $X^n$,
\item for $a\in X^n$ we denote by $o_G(a)$ the orbit of the tuple $a$, i.e, the unique $n$-orbit of $G$ containing $a$,
\item for $a,b\in X^n$ we write $a\same_G b$ if $o_G(a)=o_G(b)$,
\item $\orb_n(G)$ denotes the \emph{number of $n$-orbits of $G$}, i.e., $\orb_n(G)=|\orbits_n(G)|$,
\item $\oi_n(G)$ denotes the \emph{number of injective $n$-orbits of $G$}, i.e., the number of those $n$-orbits of $G$ which contain tuples with pairwise different entries,
\item $\os_n(G)$ denotes the \emph{number of orbits of $n$-subsets of $G$}, i.e., the number of orbits of the natural action $G\acts {X\choose n}=\{Y\subset X\colon|Y|=n \}$.
\end{itemize}

\begin{remark}\label{ois}
	It follows easily from the definitions that for any permutation group $G$ we have $\os_n(G)\leq \oi_n(G)\leq n!\os_n(G)$ and $\orb_n(G)=\sum_{k=1}^n{n\brace k}\oi_k(G)$ where ${n\brace k}$ denotes the Stirling number of the second kind.
\end{remark}

	A permutation group is called \emph{oligomorphic} if and only if $\orb_n(G)$ is finite for all $n\in \omega$. By Remark~\ref{ois} it follows that the oligomorphicity of a group $G$ is also equivalent to $\os_n(G)$ (or $\oi_n(G)$) being finite for all $n\in \omega$. We call the sequence $(\oi_n(G))_{n\in \omega}$ the \emph{\textbf{l}abelled growth profile} of $G$, and $(\os_n(G))_{n\in \omega}$ the \emph{\textbf{u}nlabelled growth profile} of $G$. The orbit growth functions $\oi_n(\cdot)$ and $\os_n(\cdot)$ are introduced and discussed in general in~\cite{CameronCounting,Oligo}. Note that at some places in the literature the sequences $\orb_n,\oi_n$ and $\os_n$ are denoted by $F_n^*,F_n$ and $f_n$, respectively. In this paper we are mostly concerned with the unlabelled growth profile of permutation groups.
	
	Let $G\leq \sym(X)$ be a permutation group, and let $S\subset X$. Then we write $\acl_G(X)$ for the union of finite $G$-orbits of elements contained in $X$.\footnote{The corresponding notion for structures is called \emph{algebraic closure}, hence the notation.} Note that if $X$ is finite and $G$ is oligomorphic then $\acl_G(X)$ is also finite.

\subsection{Automorphisms of structures}

	In this paper we denote structures by Fraktur letters $\fa,\fb,\fc,\dots$ (sometimes with subscripts or other diacritics), and we write $A,B,C,\dots$ for the corresponding domain set. The signature of a structure $\fa$ is denoted by $\sign(\fa)$. For a structure $\fa$, we write $\age(\fa)$ (called the \emph{age} of $\fa$) for the class of all finite
structures which embed into $\fa$. We write $\aut(\fa)$ for the automorphism group of a structure $\fa$. It is easy to see that the automorphism group of every structure is closed. Conversely, every closed permutation group is the automorphism group of some structure. In this paper every structure is assumed to be countable, unless otherwise mentioned. We also assume for convenience that all structures in our discussions have a purely relational signature.\footnote{Constants can be considered to be unary relations represented by a single element.}

	We are considering the following two notions of equivalence on first-order structures.
	
\begin{definition}
	Let $\fa$ and $\fb$ be structures. Then we say that 
\begin{itemize}
\item $\fb$ is a \emph{reduct} of $\fa$ if $A=B$ if all relations, and constants of $\fb$ are definable in $\fa$ by a first-order formula,
\item $\fa$ and $\fb$ are \emph{interdefinable} if they are reducts of one another
\item $\fa$ and $\fb$ are \emph{bidefinable} if $\fb$ is isomorphic to some structure which is interdefinable with $\fa$.
\end{itemize}
\end{definition}

\begin{remark}
	Note that bidefinability is a symmetric condition. Indeed, if $\fa$ and $\fb'$ are interdefinable and $\iota:\fb'\rightarrow \fb$ and isomorphism then the structure $\fa'\coloneqq \iota(\fa)$ is interdefinable with $\fb$. 
\end{remark}

	We say that a countable structure $\fa$ is $\omega$-categorical if every countable model of the first-order theory of $\fa$ is isomorphic to $\fa$. By the theorem of Engeler, Ryll-Nardzewski, and Svenonius, we know that a countable structure $\fa$ is $\omega$-categorical if and only if $\aut(\fa)$ is oligomorphic (see for instance~\cite{HodgesLong}, Theorem 7.3.1). The connection between $\omega$-categorical structures and their automorphism groups is summarized in the following theorem.

	
\begin{theorem}\label{aut_inv}
	Let $\fa$ and $\fb$ structures, and let us assume that $\fa$ is $\omega$-categorical. Then the following hold.
\begin{enumerate}
\item\label{it:reduct} $\fb$ is a reduct of $\fa$ if and only $\aut(\fa)\leq \aut(\fb)$.
\item\label{it:inter} $\fa$ and $\fb$ are interdefinable if and only if $\aut(\fa)=\aut(\fb)$.
\item\label{it:bi} $\fa$ and $\fb$ are bidefinable if and only if $\aut(\fa)$ and $\aut(\fb)$ are isomorphic (as permutation groups).
\item\label{it:red_inter} The map $\fc\mapsto \aut(\fc)$ defines a bijection between the reducts of $\fc$ up to interdefinability and the closed supergroups of $\aut(\fa)$.
\item\label{it:red_bi} The map $\fc\mapsto \aut(\fc)$ defines a bijection between the reducts of $\fc$ up to bidefinability and the closed supergroups of $\aut(\fa)$ up to isomorphism.
\end{enumerate}
\end{theorem}

\begin{proof}
	We know that a relation $R$ is first-order definable in an $\omega$-categorical structure $\fa$ if and only if $R$ is preserved by $\aut(\fa)$ (see for instance~\cite{HodgesLong}, Corollary 7.3.3). This immediately implies item~\ref{it:reduct}. The implications~\ref{it:reduct}$\rightarrow$~\ref{it:inter} and~\ref{it:inter}$\rightarrow$~\ref{it:bi} are obvious. 
	
	For item~\ref{it:bi} we only need to argue that closed supergroups of $\aut(\fa)$ are exactly the automorphism groups of reducts of $\fa$. Clearly, if $\fb$ is a reduct of $\fa$ then $\aut(\fa)\leq \aut(\fb)$ and $\aut(\fb)$ is closed. Conversely, if $\aut(\fa)\leq G$ and $G$ is closed then, as we discussed earlier, $G=\aut(\fb)$ for some structure $\fb$, and then $\fb$ is a reduct of $\fa$ by item~\ref{it:reduct}.
\end{proof}

	Theorem~\ref{aut_inv},~\ref{it:bi} tells us that in order to classify some class $\mathscr{C}$ of $\omega$-categorical structures which is closed under bidefinability, it is enough to classify the automorphism groups of structures in $\mathscr{C}$ up to isomorphism.

	In this paper we use the following notational conventions: we use \texttt{mathcal} letters ($\mathcal{H},\mathcal{S},\mathcal{F}$) for classes of permutations groups, and \texttt{mathscr} letters ($\mathscr{H},\mathscr{S},\mathscr{F}$) for the corresponding structures, (for instance we write $\mathscr{C}$ for the class of all structures with automorphism groups in $\mathcal{C}$). If all groups in $\mathcal{C}$ are closed and then we have $\mathcal{C}=\{\aut(\fa): \fa\in \mathscr{C}\}$.

\subsection{Direct products and disjoint unions}\label{sect:direct}

	Let $\{X_i\colon i\in I\}$ be a family of pairwise disjoint sets, and for each $i\in I$ let $G_i\leq \sym(X_i)$ be a permutation group. Then we define the \emph{direct product} of the permutation groups $G_i$, denoted by $\prod_{i\in I}G_i$, to be the set of all permutations of $X:=\bigcup\{X_i\colon i\in I\}$ which can be written as $\bigcup\{f(i)\colon i\in I\}$, for some $f\colon I\to\prod_{i\in I}G_i$ that satisfies $f(i)\in G_i$ for all $i\in I$. In this case $\prod_{i\in I}G_i$ as a group is also the usual direct product of the groups $G_i\colon i\in I$. It is clear from the definition that direct products of closed groups are closed.
	
	Direct products of groups correspond to disjoint unions of structures in the following sense.
	
\begin{definition}\label{def:union}
	Let $(\fa_i)_{i\in I}$ be a family of relational structures with pairwise disjoint domains. 

	Then we define the \emph{disjoint union} $\fb$ of the structures $\fa_1,\dots,\fa_n$, denoted by $\biguplus_{i=1}^n\fa_i$ such that
\begin{itemize}
\item the domain set of $\fb$ is $\bigcup_{i=1}^nA_i$,
\item its signature of $\fb$ is the union of the signatures of $\fa_i: i\in I$ together with some (distinct) unary symbols $U_i\not\in \bigcup_{i\in I}\sign(\fa_i)$,
\item $U_i^\fb=A_i$ for all $1\leq i\leq n$ and,
\item $R^\fb=R^{\fa_i}$ for all $R\in \sign(\fa_i)$.
\end{itemize}
\end{definition}

	It is easy to see that if $(\fa_i)_{i\in I}$ are as above then $\aut(\biguplus_{i=1}^n\fa_i)=\prod_{i\in I}\aut(\fa_i)$. It is also clear that finite unions of $\omega$-categorical structures are $\omega$-categorical which in turn implies that oligomorphic groups are closed under taking finite direct products.

\subsection{Wreath products}

	If $G\leq \sym(X)$ is a permutation group, and $S$ is a set then we define the power $G^S$ to be $\prod_{s\in S}G_s$ where $G_s:=\{(g,s)\colon g\in G\}$ and $(g,s)(g',s)=(gg',s)$. This makes sense since the sets $\dom(G_i)=X\times \{i\}$ are pairwise disjoint. Then $\dom(G^S)=X\times S$, and since all $G_i$ are isomorphic to $G$ it follows that $G^S$ is isomorphic to the power $G^{|S|}$ in the usual group theoretical sense. 
	
	Now let $H\leq \sym(Y)$ be another permutation group. Then we define \emph{wreath product} of the groups $G$ and $H$, denoted by $G\wr H$, as follows. The domain set of $G\wr H$ is $X\times Y$ and a permutation $\sigma\in \sym(X\times Y)$ is contained in $G\wr H$ if and only if there exist $\beta\in H$ and $\alpha_y\in G$ such that $\sigma(x,y)=(\alpha_y(x),\beta(y))$. Note that in this case we have a natural surjective homomorphism $\phi$ from $G\wr H$ to $H$ by mapping every $\sigma$ to its action of the second coordinates. The kernel of this homomorphism is exactly the power $G^Y$, as defined above. Moreover, the map $\iota\colon \beta\mapsto ((x,y)\mapsto (x,\beta(y))$ splits the homomorphism $\phi$, that is, $\phi\circ \iota=\id(Y)$. This implies that $G\wr H$ can be written as a semidirect product $G^Y\rtimes e(H)$. 
	
	We can also define wreath products of structures as follows.

\begin{definition}\label{def:copies}
	Let $\fa$ and $\fb$ be structures with disjoint relational signatures. Then we define the structure $\fa\wr \fb$, called the \emph{wreath product of $\fa$ and $\fb$} as follows. 
\begin{itemize}
\item The domain set of $\fa\wr \fb$ is $A\times B$, 
\item the signature of $\fa\wr \fb$ is $\sign(\fa)\cup \sign(\fb)\cup\{E\}$ where $E$ is a distinguished binary symbol not contained in $\sign(\fa)\cup \sign(\fb)$,
\item $E^{\fa\wr \fb}=\{(a_1,b),(a_2,b)\colon a_1,a_2\in A, b\in B\}$,
\item $R^{\fa\wr \fb}=\{(a_1,b),\dots,(a_k,b)\colon (a_1,\dots,a_k)\in R^{\fa}, b\in B\}$ for all $R\in \sign(\fa)$,
\item $S^{\fa\wr \fb}=\{(a_1,b_1),\dots,(a_k,b_1)\colon (b_1,\dots,b_k)\in S^{\fb}\}$ for all $S\in \sign(\fb)$.
\end{itemize}
\end{definition}

	The correspondence between wreath products of groups and wreath products of structures is as follows. Let $\fa$ and $\fb$ be a structure as in Definition~\ref{def:copies}. Then we have $\aut(\fa\wr \fb)=\aut(\fa)\wr \aut(\fb)$. This implies immediately that if $G$ and $H$ are closed then so is $G\wr H$. Moreover, if $\fa$ and $\fb$ are $\omega$-categorical then so is $\fa\wr \fb$ which implies that if $G$ and $H$ are oligomorphic groups then $G\wr H$ is also oligomorphic.

\subsection{Covers of permutation groups}

	In this subsection we introduce the notion of covers of permutation groups which plays a crucial role in this article. This construction will allow us to decompose certain oligomorphic permutation groups into simpler parts from which the original group can be recovered. Our definition of covers is largely based on the ones introduced in~\cite{ivanov1999finite,EvansIvanovMacpherson,pastori,bodirsky2021permutation}. In some cases covers are considered for structures and in some other cases it is only considered for \emph{closed} permutation groups (these are called \emph{permutation structures} in~\cite{EvansIvanovMacpherson} and~\cite{pastori}). One important difference in our definition compared to the earlier ones is that we also define everything for general permutation groups without any closedness assumption. Although all permutation groups showing up during our constructions will be closed, this will sometimes not follow immediately from the definition, so we prefer this more general treatment.
	
\begin{notation}
	Let $\pi\colon X\rightarrow Y$ be a map. Then we write $\sim_{\pi}$ for the equivalence relation $\{(a_1,a_2): \pi(a_1)=\pi(a_2)\}$. If $G\leq \sym(X)$ and $\sim_{\pi}$ is a congruence of $G$ then $\pi$ gives rise to a homomorphism $X\rightarrow \sym(Y)$ defined by $\mu_{\pi}(g)(a)=\pi(g(\pi^{-1}(a)))$. We that mapping $\pi$ is clear from the context we simply write $\mu$ for $\mu_{\pi}$.
\end{notation}

\begin{definition}\label{def:cover}
	Let $G$ and $G^*$ be permutation groups. Then a triple $\Phi=(G,G^*,\pi)$ is called a \emph{cover} of $G^*$ if
\begin{enumerate}
\item $\pi$ is a surjective map from $\dom(G)$ to $\dom(G^*)$,
\item the equivalence relation $\sim_{\pi}$ is preserved by $G$,
\item the image of $G$ under $\mu_{\pi}$ equals $G^*$.
\end{enumerate}
	The sets $\pi^{-1}(w)\colon w\in \dom(G^*)$ are called the \emph{fibers} of $\Phi$.

	Let $a\in \dom(G^*)$, and let $A:=\pi^{-1}(a)$. Then we call the groups
\begin{itemize}
\item $G_{(A)}$,
\item $\mu_{\pi}^{-1}(G_a)$,
\item $G_{((A))}$
\end{itemize}
the \emph{fiber group}, the \emph{binding group} and the \emph{pointwise binding group} of $\Phi$ at $a$, respectively. The \emph{kernel of $\Phi$} is the kernel of the homomorphism $\mu_{\pi}$. 
	
	If all fibers of $\Phi$ are finite then $\Phi$ is called a \emph{finite cover} of $G^*$.

	Note that the binding group and the pointwise binding group at a given point are always normal subgroups of the corresponding fiber group. We say that $\Phi$ has \emph{finite fiber factors} if all factor groups $G_{(A)}/G_{((A))}\colon A\in \{\pi^{-1}(a):a\in \dom(G^*)\}$ are finite. Clearly, all finite covers have finite fiber factors.
	
	Note that if $G$ is closed then the direct product of the pointwise binding groups of $\Phi$ is a normal subgroup of $G$. We say that a $\Phi$ is \emph{linked} if all of its pointwise binding groups are trivial.
	
	We say that two covers $(G,G^*,\pi)$ and $(G',(G')^*,\pi')$ are \emph{isomorphic} if there exist bijections $e\colon \dom(G)\rightarrow \dom(G'),e^*\colon \dom(G^*)\rightarrow \dom((G')^*)$ such that $\mu_{e^*}\circ \pi=\pi'\circ \mu_e$. So the isomorphism of two covers essentially means that there are the same up to renaming the domains of the involved groups.

	We say that a group $G$ is a \emph{cover} of $G^*$ if there exists a $\pi\colon G\rightarrow G^*$ such that $\Phi=(G,G^*,\pi)$ is a cover of $G^*$. A structure $\fa$ is called a cover of $\fb$ is $\aut(\fa)$ is a cover of $\aut(\fb)$. Finite covers of groups and structures are defined analogously.
\end{definition}

\begin{remark}
	Let $G\leq \sym(X)$ be an arbitrary permutation group, and let $\sim$ be a congruence of $G$. Let us define $\pi\colon a\mapsto [a]_\sim$. Then it is easy to see that $\sim=\sim_{\pi}$ and $(G,G/\unsim,\pi)$ is a cover. In fact, every cover a permutation group arises from the construction above. Indeed, let $(G,G^*,\pi)$ be a cover. Let $\sim=\sim_{\pi}$. Then there is natural bijection between $X/\unsim_{\pi}$ and $Y=\dom(G^*)$ defined by $w\mapsto \pi^{-1}(w)$. By renaming the elements of $Y$ along this bijection we obtain $G^*=G/\unsim$ and $\pi(a)=[a]_{\sim_{\pi}}$.
\end{remark}

\begin{notation}
	For a tuple $a=(a_1,\dots,a_k)$ we write $\pi_i(a)=a_i$.
\end{notation}

\begin{observation}
	Let $(G,G^*,\pi)$ be a cover. Then we can relabel the elements of the domain of $G$ such that the following hold.
	
\begin{enumerate}[(a)]
\item every element of $\dom(G)$ is a pair
\item $\pi=\pi_2$
\item if $a\same_{G^*}b$ then $\pi_1(\pi_2^{-1}(a))=\pi_1(\pi_2^{-1}(b))$, and the map $(f,a)\mapsto (f,b)$ extends to a permutation in $G$. 
\end{enumerate}
\end{observation}

\begin{proof}
	Let us first relabel each element $a\in \dom(G)$ by the pair $(a,\pi(a))$. Then items (a) and (b) are already satisfied.

	Now let $O$ be an orbit of $G^*$, and let us fix an element $a\in O$. Now for each $b\in O$ let us consider an element $g_b\in G$ such that $\mu(g_b)(a)=b$, and let us relabel each element of the form $(u,b)$ to $(g_b^{-1}(u),b)$. By doing this element for all orbits of $G^*$ it is easy to check that item (c) will also be satisfied.
\end{proof}
	
	To ease the notation in our later discussions we always consider covers satisfying conditions (a)-(c) above. This way every cover of some group $G^*$ is already specified by the group $G$. Also, in this case we can (and we will) identify the fiber groups, binding groups and pointwise binding groups with their action of the first coordinates. By condition (c) these groups are the same, not just isomorphic, at elements coming from the same orbit. From this point on when we say $G$ is a cover of $G^*$ we always automatically assume that the triple $(G,G^*,\pi_2)$ satisfies conditions (a)-(c).
	
\begin{definition}
	We say that a cover $(G,G^*,\pi)$ is
\begin{itemize}
\item \emph{trivial} its kernel is trivial (only contains the identity permutation),
\item \emph{strongly trivial} if all of its fiber groups are trivial,
\item \emph{(strongly) split} if there is a $G_0\leq G$ such that $(G_0,G^*,\pi)$ is (strongly) trivial.
\end{itemize}
\end{definition}
	
\begin{example}\label{strongly_trivial}
	Let $F$ be an arbitrary set, and let $\id(F)\wr G^*\leq G\leq \sym(F)\wr G^*$, and $\pi=\pi_2$. Then $(G,G^*,\pi)$ is a strongly split cover of $G^*$.
\end{example}

	It follows from the discussion in~\cite{bodor2024classification}, Section 2.9 that if $G^*$ is transitive then every strongly split cover of $G^*$ is isomorphic to one given in Example~\ref{strongly_trivial}.
	
	Let us assume that $(G,G^*,\pi)$ is a cover with kernel $K$, and let us assume that $\phi$ is some injective map from $G^*\rightarrow G$ such that $\mu_{\pi}\circ \phi=\id$. Then $G$ can be written as $G=\ker(\mu)\phi(G^*)$. Note that $\Phi$ is split if and only if $\phi$ can be chosen to be a group homomorphism. If $G$ is as in Example~\ref{strongly_trivial} then $\phi$ can be chosen such that $\phi(\sigma)=\hat{\sigma}$ where $\hat{\sigma}$ is defined to be the permutation which maps the $(f,a)$ to $(\sigma(f),a)$ for $(f,a)\in F\times \dom(G^*)$.
	
	Now we introduce some further notation which will be useful in the discussion of the kernels of strongly split covers. Let $X$ and $F$ be sets, and let $G\coloneqq \sym(F)\wr \sym(X)$. Then for all $a\in X:=\dom(G)$ and $g\in G$ we define $\delta_a(g)\colon f\mapsto \pi_1(g(f,a))$. Then clearly $g(f,a)=(\delta_a(f),\mu_{\pi_2}(a))$. Let $H\triangleleft L\leq \sym(F)$. Then we define $$\kernel{L}{H}=\{g\in G\cap \ker(\mu_{\pi_2})\colon \forall a,b\in X(\delta_a(g)\in L\wedge \delta_a(g)H=\delta_b(g)H)\}.$$ Then $\kernel{L}{H}$ always forms a subgroup of $G$. Moreover, if $H$ and $L$ are closed (in particular if $F$ is finite) then so is $\kernel{L}{H}$.

\begin{notation}
	Let $S\subset \dom(G^*)$ and  $\sigma\in \sym(F)$. Then we denote by $\kappa(\sigma,S)$ the unique permutation $\gamma\in G\cap \ker(\mu_{\pi_2})$ such that
\[
\delta_a(\kappa)=
\begin{cases}
\sigma &\text{if $a\in S$}\\
\id(F) &\text{otherwise.}
\end{cases}
\]

	We also write $\kappa(\sigma,a)$ for $\kappa(\sigma,\{a\})$.
\end{notation}

\subsection{Hereditarily cellular groups and structures}\label{sect:her}

	In this subsection we briefly recall the classification of hereditarily cellular structures mostly taken from~\cite{bodor2024classification}.
	
	We first need to introduce a couple of notations.
	
\begin{definition}\label{def:bigg}
	Let $N_i\leq \sym(Y_i)$ be permutation groups for some pairwise disjoint sets $Y_0,\dots,Y_k$, and let $H$ be a permutation group which normalizes $N\coloneqq \prod_{i=0}^kN_i$. Then we write $\group(H;N_0,\dots,N_k)$ for the permutation group with domain $(Y_0\times \{0\})\cup \bigcup_{i=1}^k(Y_i\times\omega)$ generated by $(N_0\wr (\id(\{0\}))\times \prod_{i=1}^k(N_i\wr \sym(\omega))$ and the permutations $(a,n)\mapsto (h(a),n)\colon h\in H$.
\end{definition}

	Note that the definition of the group $\group(H;N_0,\dots,N_k)$ differs from the one given in~\cite{bodor2024classification}, however the equivalence of the two definitions follows from the discussion in Section 3 in the aforementioned paper.
	
	It is clear from the definition that $$\group(N;N_0,\dots,N_k)=(N_0\wr (\id(\{0\}))\times \prod_{i=1}^k(N_i\wr \sym(\omega)).$$
	
	Now we are ready to give the definition of hereditarily cellular groups/structures.
	
\begin{definition}\label{def:gone}
	We define the classes of closed permutation groups $\gone_i\colon i\in \{-1\}\cup \omega$ recursively as follows.
\begin{itemize}
\item $\gone_{-1}=\emptyset$,
\item For $n\geq 0$ a group $G$ is contained in $\gone_n$ if and only if $G$ is isomorphic to $\group(H;N_0,\dots,N_k)$ for some $N_0,N_1,\dots,N_k,H\in \gone_{n-1}$ where $Y_0=\dom(N_0)$ is finite, and $N_0=\id(Y_0)$.
\end{itemize}

	We write $\gone$ for the union of all classes $\gone_n\colon n\in \omega$. A group $G$ is called \emph{hereditarily cellular} if $G\in \gone$, and a structure is called \emph{hereditarily cellular} if its automorphism group is hereditarily cellular. Following our notational convention we write $\gones$ for the class of hereditarily cellular structures.
	
	The \emph{rank} of a hereditarily cellular group $G$, denoted by $\rk(G)$ is the smallest $n$ such that $G\in \gone_n$. The rank of a hereditarily cellular structure is the rank of its automorphism group.
	
	Hereditarily cellular groups/structures with rank at most 1 are called \emph{cellular}.
\end{definition}

\begin{remark}\label{rank0}
	It is clear from the definition that $\gone_0$ is exactly the class of permutation groups with finite degree.
\end{remark}

	Note that in~\cite{bodor2024classification} the classes $\gone_n$ is defined via \emph{$\omega$-partitions}, a notation introduced in~\cite{lachlan1992}, but again these definitions are known to be equivalent. We recall this equivalence below.
	
\begin{definition}\label{def:omega_part}
	Let $G\leq \sym(X)$ be a permutation group. Then a triple $K, \bigeq, \smalleq)$ is called an \emph{$\omega$-partition} of $G$ iff
\begin{enumerate}
\item $K\subset X$ is finite, and it is fixed setwise by $G$.
\item $\smalleq$ and $\bigeq$ are congruences of $G|_{X\setminus K}$ with $\smalleq\subset \bigeq$.
\item $\bigeq$ has finitely many classes.
\item Every $\bigeq$-class is a union of $\aleph_0$ many $\smalleq$-classes.
\item For every $C\in X/\bigeq$ we have $G_{((C))}/\smalleq=\sym(C/\smalleq)$.
\end{enumerate}
	The \emph{components} of an $\omega$-partition $(K, \bigeq, \smalleq)$ are the permutation groups 
$G_{((Y))}$ for all $Y\in X/\smalleq$.
\end{definition}

\begin{theorem}[~\cite{bodor2024classification}, Theorem 3.41]\label{omega_part}
	Let $n\geq 0$. Then $G\in \gone_n$ if and only if $G$ has an $\omega$-partition with components in $\gone_{n-1}$.
\end{theorem}

\begin{remark}
	We know from the results of~\cite{lachlan1992} (and Theorem~\ref{omega_part}) that if $\fa\in \gones$ then $\rk(\fa)$ also equals to the Morley rank of $\fa$. In particular, the Morley rank of every hereditarily cellular structure is finite.
\end{remark}

	Next we recall a few important facts about hereditarily cellular groups that we will need in this paper.
	
\begin{theorem}[\cite{bodor2024classification}, Corollary 3.43]\label{thm:build_groups}
	$\gone$ is the smallest class of groups which contains $\id(\{\emptyset\})$, closed under isomorphisms, and is closed under taking
\begin{itemize}
\item finite direct products,
\item wreath products with $\sym(\omega)$, and
\item finite-index supergroups.
\end{itemize}
	In particular, every hereditarily cellular group is closed.
\end{theorem}

	Theorem~\ref{thm:build_groups} can also be formulated in terms of structures. We say that a reduct $\fb$ of $\fa$ is of \emph{finite index} if $|\aut(\fb):\aut(\fa)|$ is finite.
	
\begin{theorem}[\cite{bodor2024classification}, Theorem 3.44]\label{thm:build_structures}
	A class $\gones$ is the smallest class of structures which contains the one-element pure set, closed under isomorphisms, and is closed under taking
\begin{itemize}
\item finite disjoint unions,
\item wreath products with the pure set, and
\item finite-index reducts.
\end{itemize}
\end{theorem}

\begin{notation}
	For a permutation group $G$ we denote by $G_{\fin}$ the intersection of all finite-index closed subgroups of $G$.
\end{notation}

\begin{remark}
	Note that if $|G:H|$ is finite, then $H$ contains a normal subgroup $N$ which has also finite index in $G$. Indeed, we can take $N$ to be the kernel of the action of $G$ by right multiplication on the set of right cosets of $H$. In this case $N=\bigcup_{g\in G}g^{-1}Hg$, so if $H$ is closed then so is $N$. This implies that $G_{\fin}$ is also equal to the intersection of finite-index \emph{normal} closed subgroups of $G$.
\end{remark}

\begin{lemma}\label{fin_index_gone}
	Let $G\in \gone_n$. Then $|G:G_{\fin}|$ is finite, and for every group $G_{\fin}\leq G_*\leq G$ we have $G_*\in \gone_n$.
\end{lemma}

\begin{proof}
	By Lemma 3.34 in~\cite{bodor2024classification} we know that if $G=\biggroup(H;N_0,\dots,N_k)$ and $|G:G_*|$ is finite then $G_*\supset \biggroup(N;N_0,\dots,N_k)\in \gone_n$. It also follows from the proof of Lemma 3.34 in~\cite{bodor2024classification} that $|G:\biggroup(N_0,\dots,N_k)|$ is itself finite. Therefore, $G_{\fin}=\biggroup(N;N_0,\dots,N_k)$. In particular, $G_{\fin}\in \gone_n$. Since $|G_*:G_{\fin}|<|G:G_{\fin}|$ is finite it follows that $G_*$ is closed, and thus $G_*\in \gone_n$.
\end{proof}

	By the results of~\cite{lachlan1992} we know that a structure is hereditarily cellular if and only if it is $\omega$-categorical and monadically stable.
	
\begin{definition}
	We say that a structure $\fa$ is \emph{monadically stable} iff every expansion of its theory by unary predicates is stable.
\end{definition}

	The class of hereditarily cellular structures can also be described by their unlabelled growth profile. We say that a function $f\colon\omega\rightarrow \mathbb{R}$ is \emph{slower (faster) than exponential} if for all $c>1$ we have $f(n)<c^n$ ($f(n)>c^n$) if $n$ is large enough. Using these notions we have the following dichotomy result for the unlabelled growth profile of $\omega$-categorical stable structures.

\begin{theorem}[\cite{braunfeld2022monadic}, Theorem 3.10]\label{sam}
	Let $\fa$ be an $\omega$-categorical stable structure. Then exactly one of the following holds.
\begin{enumerate}
\item $\os_n(\fa)$ is slower than exponential, and $\fa$ is monadically stable.
\item $\os_n(\fa)$ is faster than exponential, and $\fa$ is not monadically stable.
\end{enumerate}
\end{theorem}

	Next we recall a description of cellular structures using their unlabelled growth profile.

\begin{definition}
	A permutation group $G$ is called \emph{$P$-oligomorphic} if $\os_n(\fa)$ has a polynomial upper bound.
\end{definition}

	The following two theorems follow from the results of~\cite{falque2019classification} and~\cite{braunfeld2022monadic}.

\begin{theorem}\label{poly}
	Let $\fa$ be a countable structure. Then the following are equivalent.
\begin{enumerate}
\item $\os_n(\fa)$ is P-oligomorphic.
\item $\os_n(\fa)\sim cn^k$ for some $c\geq 0$ and $k\in \omega$.
\item For all $c>0$ we have $\os_n(\fa)< e^{c\sqrt{n}}$ if $n$ is large enough.\footnote{Here we also allow $c=0$ to include finite structures.}
\end{enumerate}
\end{theorem}

\begin{theorem}\label{cell}
	Let $\fa$ be a countable structure. Then $\fa$ is cellular if and only if $\fa$ is stable and $\aut(\fa)$ is $P$-oligomorphic.
\end{theorem}

	For further discussion on the unlabelled growth profile of hereditarily cellular structures we refer the reader to~\cite{braunfeld2022monadic}.
	
\subsection{Homogeneous and finitely bounded structures}

	Next we recall some definitions and basic results concerning homogeneity and finite boundedness of structures.

\begin{definition}
	A structure $\fa$ is called \emph{homogeneous} if every isomorphism between finite substructures of $\fa$ can be extended to an automorphism of $\fa$.
	
	We say that a structure is \emph{$m$-homogenizable} if it is interdefinable with a homogeneous structure with a finite relational signature all of whose relations have arity of at most $m$. We say that a structure is \emph{finitely homogenizable} if it is $m$-homogenizable for some $m\in \omega$.
\end{definition}

\begin{remark}
	In the general definition of homogeneity we have to consider \emph{finitely generated} rather than finite substructures. However, since we only consider purely relational structures, this does not make a difference in our setting.
\end{remark}

	An easy calculation on types shows that every homogeneous structure with a finite relational signature is $\omega$-categorical; in fact if $\fa$ is homogeneous in a finite relational signature and the maximal arity of relations of $\fa$ is $m$ then $\orb_n(\fa)\leq 2^{cn^{\max(2,m)}}$ for some $c>0$. 
		
	It turns out that $m$-homogenizable and in turn finitely homogenizable structures can also be characterized in a natural way in terms of their automorphism groups.

\begin{notation}
	For any $n$-tuple $a$ and a function $p\colon [m]\rightarrow [n]$ we write $a_{\pi}$ for the tuple $(a_{p(1)},\dots,a_{p(m)})$.
\end{notation}

\begin{definition}
	We say that a permutation group $G$ is \emph{$m$-sensitive} if for all $u,v\in A^n$ if $u\notsame_G v$ there exists a function $p\colon [m]\rightarrow [n]$ such that $u_p\notsame_G v_p$. We say that $G$ is \emph{finitely sensitive} if it is $m$-sensitive for some $m\in \omega$.
\end{definition}

\begin{lemma}\label{when_finhom}
	Let $G$ be closed oligomorphic group, and let $m\geq 2$.
Then the following are equivalent.
\begin{enumerate}

\item\label{it:homgroup} $G$ is $m$-sensitive.
\item $G$ is the automorphism group of a homogeneous relational signature structure with arity at most $m$.
\end{enumerate}
\end{lemma}

\begin{proof}
	Follows essentially from~\cite{bodor2024classification}, Lemma 2.21.
\end{proof}

	Lemma~\ref{when_finhom} can also be formulated for structures.

\begin{lemma}\label{when_finhom2}
	Let $\fa$ be an $\omega$-categorical structure, and $m\in \omega$. Then the following are equivalent.
\begin{enumerate}
\item $\fa$ is $m$-homogenizable.
\item $\aut(\fa)$ is $m$-sensitive.
\end{enumerate}
\end{lemma}

	In the study of homogeneous structures we are particularly interested in the ones which are \emph{finitely bounded}, i.e., whose age can be described by finitely many forbidden substructures. The precise definition is as follows.

\begin{definition}\label{def:fbh}
	Let $\tau$ be a relational signature, and let $\mathscr{C}$ be a set of finite $\tau$-structures. Then we denote by $\forb(\mathscr{C})$ the class of those finite $\tau$-structures which do not embed any structure from $\mathscr{C}$.

	A relational structure $\fa$ is called \emph{finitely bounded} if its signature $\tau:=\sign(\fa)$ is finite, and there is a finite set $\mathscr{C}$ of $\tau$-structures such that $\age(\fa)=\forb(\mathscr{C})$.
	
	We write $\fbh$ for the class of structures which are interdefinable with a finitely bounded homogeneous relational structure.
\end{definition}

	We close this subsection by recalling some results about finite-index reducts of homogeneous structures.
	
\begin{lemma}[\cite{bodor2024classification}, Lemma 4.2]\label{fin_index}
	Let $H\leq G$ be permutation groups with $|G:H|=d<\infty$, and let us assume that $H$ is $m$-sensitive for some $m\geq 2$. Then $G$ is $dm$-sensitive.
\end{lemma}

	As usual, Lemma~\ref{fin_index} can also be formulated in terms of structures.

\begin{corollary}\label{finite_index_s}
	Let $\fa$ be an $m$-homogenizable structure, and let $\fb$ be a reduct of $\fa$ with $|\aut(\fb):\aut(\fa)|=d$. Then $\fb$ is a $dm$-homogenizable. In particular, the class of finitely homogenizable structures are closed under taking finite-index reducts.
\end{corollary}

	We know in fact that the second statement of Corollary~\ref{finite_index_s} also holds for the class $\fbh$.
	
\begin{lemma}[\cite{bodor2024classification}]\label{fin_index_fbh}
	The class $\fbh$ is closed under taking finite-index reducts, finite disjoint unions and wreath products with the pure set.
\end{lemma}

	Note that Lemma~\ref{fin_index_fbh} together with~\ref{thm:build_groups} immediately implies that $\gones\subset \fbh$. In this article we will generalize this result to the class $\subexs$ by following a similar argument.

\subsection{Highly set-transitive groups}\label{sect:high}

	We say that a permutation group $G$ is \emph{highly set-transitive} if $\os_n(G)=1$, that is, $G$ acts transitively on the set of $n$-subsets of $A$ for all $n\in \omega$. In~\cite{Cameron5} Cameron classified all closed highly set-transitive groups.
	
\begin{theorem}\label{cameron}
	Let $G$ be a closed highly set-transitive permutation group acting on a countable set. Then $G$ is isomorphic to the automorphism group of one of the following 5 structures.
\begin{itemize}
\item $(\mathbb{Q};<)$: the rationals equipped with the usual ordering;
\item $(\mathbb{Q};\betw)$ where $\betw$ denotes the \emph{betweenness} relation: $$\betw(x,y,z)\Leftrightarrow (x<y<z)\vee (z<y<x);$$
\item $(\mathbb{Q};\cyc)$ where $\cyc$ denotes the \emph{circular order}: $$\cyc(x,y,z)\Leftrightarrow (x<y<z)\vee (y<z<x)\vee (z<x<y);$$
\item $(\mathbb{Q};\sep)$ where $\sep$ denotes the \emph{separation relation}: 
\begin{align*}
\sep(x,y,z,t)\Leftrightarrow &\left(\cyc(x,y,z)\wedge \cyc(y,z,t)\wedge \cyc(z,t,x)\wedge \cyc(t,x,y)\right)\vee \\
&\left(\cyc(z,y,x)\wedge \cyc(t,z,y)\wedge \cyc(x,t,z)\wedge \cyc(y,x,t)\right);
\end{align*}
\item the pure set $(\omega;=)$.
\end{itemize}
\end{theorem}

	Theorem~\ref{cameron} immediately implies that the structure $(\mathbb{Q};<)$ has exactly 5 reducts up to interdefinability.

\begin{theorem}\label{reducts_of_q}
	Let $\fa$ be a reduct of $(\mathbb{Q};<)$. Then $\fa$ is interdefinable with $(\mathbb{Q};R)$ where $R\in \{<,\betw,\cyc,\sep,=\}$.
\end{theorem}

	Next we recall how the automorphism groups of the structures $(\mathbb{Q};R)\colon R\in \{\betw,\cyc,\sep\}$ can be generated over $\aut(\mathbb{Q};<)$.
		
\begin{notation}\label{not:ft}
	Let $r$ be a fixed irrational number. Then we denote by $\varphi_+$ $(\varphi_-)$ an order preserving (resp. order reversing) bijection from $\mathbb{Q}(<r)=\{q\in \mathbb{Q}: q<r\}$ to $\mathbb{Q}(>r)=\{q\in \mathbb{Q}: q>r\}$. Then we define
\begin{itemize}
\item $\flip:=\varphi_-\cup \varphi_-^{-1}$,
\item $\turn:=\varphi_+\cup \varphi_+^{-1}$.
\end{itemize}
\end{notation}

	The following statements are easy to check and well-known in the literature.
	
\begin{facts}
\ 
\begin{itemize}
\item $\aut(\mathbb{Q};\betw)=\aut(\mathbb{Q};<)\cup \aut(\mathbb{Q};<)\!\flip$
\item $\aut(\mathbb{Q};\cyc)=\aut(\mathbb{Q};<)\cup \aut(\mathbb{Q};<)\turn\aut(\mathbb{Q};<)$
\item $\aut(\mathbb{Q};\sep)=\aut(\mathbb{Q};\cyc)\cup \aut(\mathbb{Q};\cyc)\!\flip$.
\end{itemize}
\end{facts}
	
	We will need the following two observations in our later discussion.

\begin{lemma}\label{normal_highly}
	Let $N\leq \sym(X)$ be a closed highly set-transitive permutation group, and let us assume that $N\triangleleft G\leq \sym(X)$. Then we can rename the elements of $X$ such that one of the following holds.
\begin{itemize}
\item $G=N=\aut(\mathbb{Q};R)$ where $R\in \{<,\betw,\cyc,\sep,=\}$.
\item $N=(\mathbb{Q};<), G=(\mathbb{Q};\betw)$ and $|G:N|=2$.
\item $N=(\mathbb{Q};\cyc), G=(\mathbb{Q};\sep)$ and $|G:N|=2$.
\end{itemize}
\end{lemma}

\begin{proof}
	Clearly, $N$ is also a normal subgroup of $\bar{G}$. Therefore, by Theorem~\ref{cameron} we may assume that $N,\overline{G}\in \{\aut(\mathbb{Q};R)\colon R\in \{<,\betw,\cyc,\sep,=\}\}$. It is well-known that the full symmetric group does not have any nontrivial normal closed subgroup (see for instance Corollary 2.2 in~\cite{cameron1981normal}), and it is easy to check that none of the groups $\aut(\mathbb{Q};<)$ and $\aut(\mathbb{Q};\betw)$ is a normal subgroup of either $\aut(\mathbb{Q};\cyc)$ or $\aut(\mathbb{Q};\sep)$. In the remaining cases either $N=G=\bar{G}$; or
\begin{itemize}
\item $N=(\mathbb{Q};<), G=(\mathbb{Q};\betw)$; or
\item $N=(\mathbb{Q};\cyc), G=(\mathbb{Q};\sep)$.
\end{itemize}
	In the last 2 cases we have $|\bar{G}:N|=2$, and thus in fact $G=\bar{G}$.
\end{proof}

\begin{lemma}\label{normal_highly}\label{high_const}
	Let $G\leq \sym(X)$ be a closed highly set-transitive permutation group different from $\sym(X)$, and let $s,t\in X$ with $s\neq t$. Then $G_{s,t}$ is isomorphic to either $\aut(\mathbb{Q};<)^3\times \id(\{0,1\})$ or $\aut(\mathbb{Q};<)^2\times \id(\{0,1\})$.
\end{lemma}

\begin{proof}
	By Theorem~\ref{cameron} we can assume that $G=\aut(\mathbb{Q};R)$ where $R\in \{<,\betw,\cyc,\sep,=\}$. We can also assume without loss of generality that $s<t$. Note that $\flip$ cannot be contained in $G_{s,t}$. This implies that $G_{s,t}=\aut(\mathbb{Q};<)_{s,t}$ or $G_{s,t}=\aut(\mathbb{Q};\cyc)_{s,t}$.
	
	In the first case it is clear that 
\begin{align*}
G_{s,t}=&\aut(\mathbb{Q};<)_{s,t}\\=&((-\infty,s);<)\times ((s,t);<)\times ((t,\infty);<)\times \id(\{s,t\})\simeq \aut(\mathbb{Q};<)^3\times \id(\{0,1\}).
\end{align*}
	
	Now let us assume that $G_{s,t}=\aut(\mathbb{Q};\cyc)_{s,t}$. Note that $(s,t)=\{a: (s,a,t)\in \cyc\}$, thus this interval is preserved by $\aut(\mathbb{Q};\cyc)$. Moreover, for $a,b\in (s,t)$ we have $a<b$ if and only if $(s,a,b)\in \cyc$ holds. Therefore, $(\aut(\mathbb{Q};\cyc)_{s,t})|_{(s,t)}\subset \aut((s,t);<)\simeq \aut(\mathbb{Q};<)$. On the other hand, we have $$\aut((s,t);<)\subset (\aut(\mathbb{Q};<)_{\mathbb{Q}\setminus (s,t)})|_{(s,t)}\subset (\aut(\mathbb{Q};<)_{s,t})|_{(s,t)}\subset (\aut(\mathbb{Q};\cyc)_{s,t})|_{(s,t)},$$ therefore all of these groups are equal. This implies that $$\aut(\mathbb{Q};\cyc)_{s,t}=\aut((s,t))\times \id(\{s,t\})\times (\aut(\mathbb{Q};\cyc)_{s,t})|_{\mathbb{Q}\setminus (s,t)}.$$ Let $\sigma$ be an automorphism of $(\mathbb{Q};\cyc)$ that maps $s$ to $t$ and $t$ to $s$. Then $\sigma((s,t))=\mathbb{Q}\setminus (s,t)$. This implies that the restrictions of $\aut(\mathbb{Q};\cyc)_{s,t}$ to $(s,t)$ and $\mathbb{Q}\setminus (s,t)$ are isomorphic. Therefore, we have 
\begin{align*}
G_{s,t}=&\aut(\mathbb{Q};\cyc)_{s,t}\\=&\aut((s,t);<)\times \id(\{s,t\})\times (\aut(\mathbb{Q};\cyc)_{s,t})|_{\mathbb{Q}\setminus (s,t)}\simeq \aut(\mathbb{Q};<)^2\times \id(\{0,1\}).
\end{align*}
\end{proof}

\section{Covers with finite fiber factors}\label{sect:fff}

	In this section we discuss how covers of a given group $G^*$ with finite fiber factors can be recovered from linked finite covers of $G^*$ together with some additional information on the fiber groups and the pointwise binding groups. This reconstruction is reflected in the following construction.
	
\begin{definition}\label{def:recover}
	Let $G^*\leq \sym(X^*)$ be a permutation group, and let $\tilde{G}$ be a linked finite cover of $G^*$ with domain $\tilde{X}$. Let $\data$ be a map defined on $X^*$ such that
\begin{itemize}
\item for all $a\in X^*$ we have $\data(a)=(B_a,F_a,\phi_a)$ where
\begin{itemize}
\item $F_a$ and $B_a$ are permutation groups with $B_a\triangleleft F_a$.
\item $\phi_a$ is a surjective homomorphism from $F_a$ to the fiber group of $(\tilde{G},G^*,\pi_2)$ at $a$ with kernel $B_a$.
\end{itemize}
\item $\data$ is constant on each orbit of $G^*$.
\end{itemize}

	Then we define a set $\lift(\tilde{G},G^*,\data)$ to be the set of those permutations $g$ of $\bigcup_{a\in X^*}(\dom(F_a)\times \{a\})$ which satisfy the following conditions
\begin{enumerate}
\item $g$ preserves $\sim_{\pi_2}$,
\item for all $a\in X^*$ the permutation $\eta_a(g)\colon u\mapsto \pi_1(g((u,a)))$ is contained in $F_a$, and
\item there exists a $\tilde{g}\in \tilde{G}$ such that
\begin{enumerate}
\item $\mu_{\pi_2}(g)=\mu_{\pi_2}(\tilde{g})$,
\item for all $a\in X^*$ and $u\in \pi_1(\pi_2^{-1}(a)\cap \tilde{X})$ we have $\phi_a(\eta_a(g))(u)=\pi_1(\tilde{g}((u,a)))$.
\end{enumerate}
\end{enumerate}
\end{definition}


	Next we show that the construction above indeed gives rise to a cover of $G^*$ with finite fiber factors.

\begin{lemma}\label{lift}
	Under the assumptions of Definition~\ref{def:recover} the following hold.
\begin{enumerate}
\item\label{it:lift_g} $\lift(\tilde{G},G^*,\data)$ is a group.
\item\label{it:lift_fff} $\lift(\tilde{G},G^*,\data)$ is a cover of $G^*$ with finite fiber factors.
\item\label{it:lift_closed} $\lift(\tilde{G},G^*,\data)$ is closed if and only if $\tilde{G}$ is closed and all groups $B_a$ are closed.
\end{enumerate}
\end{lemma}

\begin{proof}
	Put $G\coloneqq \lift(\tilde{G},G,\data)$.
	
	(\ref{it:lift_g}). Clearly, $\id(\bigcup_{a\in X^*}(\dom(F_a)\times \{a\}))\in G$. Now let $g,h\in G$ where condition (3) is witnessed by the permutations $\tilde{g},\tilde{h}\in \tilde{G}$, respectively. We have to show that $h^{-1}g\in G$. Clearly, $h^{-1}g$ preserves $\sim_{\pi_2}$, and it follows from the definition that $\eta_a(h^{-1}g)=(\eta_{h^{-1}g(a)}(h))^{-1}\eta_a(g)$. Moreover, we have $\mu(h^{-1}g)=\mu(h)^{-1}\mu(g)$ since $\mu$ is a homomorphism. Using this one can easily check that the map $\tilde{h}^{-1}\tilde{g}$ witnesses the fact that $h^{-1}g\in G$.
	
	Note that the map $\tilde{g}$ as in item (3) of Definition~\ref{def:recover} is unique, and in fact it follows from our argument above that the map $g\mapsto \tilde{g}$ is a surjective homomorphism from $G$ to $\tilde{G}$. We denote this homomorphism by $\psi$.
	
	(\ref{it:lift_fff}) It is clear from the definition that $(G,G^*,\pi_2)$ is a cover. Now let $a\in X^*$, and let $\data(a)=(F_a,B_a,\phi_a)$. We claim that the fiber group and the pointwise binding group at $a$ are $F_a$ and $B_a$, respectively. The former is clear from the definition. It is also clear that $B_a$ is contained in the pointwise binding group. Now let us assume that $g\in G$ fixes all elements outside $\pi_2^{-1}(a)$. Then $\mu(g)=\mu(\psi(g))=\id(\dom(G^*))$, and by definition $\psi(g)$ also fixes all elements outside $\pi_2^{-1}(a)$. Since $(\tilde{G},G,\pi_2)$ is linked, this implies that in fact $\psi(g)$ is the identity map, and therefore by definition $\mu_{\pi_1}(g|_{\pi^{-1}(a)})\in B_a$. Finally, by definition the groups $F_a/B_a$ are finite which means exactly that $(G,G^*,\pi_2)$ has finite fiber factors.
	
	(\ref{it:lift_closed}) Let us first assume that the groups $\tilde{G}$ and $B_a$ are closed. Since we have $|F_a/B_a|<\infty$ this also implies that the groups $F_a$ are also closed. Let $g\in \overline{G}$, we show that in fact $g\in G$. Let $(g_i)_i$ be a sequence in $G$ converging to $g$. The fact that $g$ preserves $\sim_{\pi_2}$ is clear. Also, for all $a\in X^*$ the sequence $(\eta_a(g_i))_i$ converges to $\eta_a(g)$. On the other hand, $\lim_{i\rightarrow \infty}(\eta_a(g_i))\in \overline{F_a}=F_a$. Thus, $g$ satisfies conditions (1) and (2) in Definition~\ref{def:recover}. Finally, it is easy to check that the sequence $(\psi(g)_i)$ is convergent, and its limit is exactly $\psi(g)$.
	
	For the other direction let us assume that $G$ is closed. Then, as we have seen before, the groups $B_a$ are exactly the pointwise binding groups of $(G,G^*,\pi_2)$, so these groups are also closed. Now let $(\tilde{g}_i)_i$ be a convergent sequence in $\tilde{G}$ with limit $\tilde{g}$. We have to show that $\tilde{g}\in \psi(G)$. Let $g_i\colon i\in \omega$ be arbitrary permutations in $G$ with $\psi^{-1}(g_i)=\tilde{g}_i$. Let $a\in G^*$. Then the sequence $(\eta_a(g_i))_i$ stabilizes eventually modulo $B_i$. In this case we can modify the elements $g_i$ by first permuting the elements of $\pi^{-1}(a)$ so that $\eta_a(g_i))(u)$ stabilizes eventually. By doing this modification for each fiber we can obtain a sequence $(g_i')_i$ in $G$ which converges, and still has the property that $\psi^{-1}(g_i')=\tilde{g_i}$. Let $g':=\lim_{i\rightarrow \infty}g_i'$. Then $g_i'\in G$ since $G$ is closed, and thus $\psi(g')=\lim_{i\rightarrow \infty}\tilde{g}_i=\tilde{g}$.
\end{proof}

	Next we show that in fact every closed cover of a permutation group with finite fiber factors arises from the construction above.
	
\begin{lemma}\label{recover2}
	Let $G^*$ be a permutation group. Then every closed cover of $G^*$ with finite fiber factors is isomorphic to $\lift(\tilde{G},G^*,\data)$ where $\tilde{G}$ and $\data$ are as in Definition~\ref{def:recover}.
\end{lemma}

\begin{proof}
	Let $(G,G^*,\pi)$ be a cover with finite fiber factors. We assume as usual that $G,G^*,\pi$ are as in conditions (a)-(c). For $a\in X^*\coloneqq \dom(G^*)$ we write $F_a$ and $B_a$ for the fiber and pointwise binding groups of $(G,G^*,\pi)$ at $a$, respectively. As before, these groups are considered by their action of the first coordinates. Since the groups $F_a/B_a$ are finite, they are isomorphic to some permutation group with finite degree. Let us pick some permutation groups $P_a: a\in X^*$ with finite degree, and some surjective homomorphism $\phi_a\colon F_a\rightarrow P_a$ with $\ker(\phi_a)=B_a$. We can do this in such a way that if $a,b\sim_{G^*}$ then $(P_a,\phi_a)=(P_b,\phi_b)$. Let us denote the domain set of $P_a$ by $D_a$, and let $\tilde{X}=\bigcup_{a\in X^*}D_a\times \{a\}$.
	
	We define a map $\psi\colon G\rightarrow \sym(\tilde{X})$ as follows. Let $g\in G$. Then, as above, for all $a\in X^*$ we write $\eta_a\colon u\mapsto \pi_1(g((u,a)))$. By using condition (c) we can conclude that $\eta_a\in F_a$. Now for a pair $(d,a)\in \tilde{X}$ we define $\psi(g)((d,a))=(\phi(\eta_a(g))(d),\mu(g)(a))$. We claim $\psi$ is a homomorphism. Indeed, for $g,h\in G$ and $(d,a)\in \tilde{X}$ we have $\eta_a(gh)=\eta_{\mu(h)(a)}(g)\eta_a(h)$, and therefore
\begin{align*}
\psi(gh)((d,a))=&(\phi(\eta_a(gh))(d),\mu(gh)(a))=(\phi(\eta_{\mu(h)(a)}(g)\eta_a(h))(d),\mu(g)(\mu(h)(a)))=\\
=&(\phi(\eta_{\mu(h)}(g))(\phi(\eta_a(h))(d)),\mu(g)(\mu(h))(a))=\\
=&\psi(g)(\phi(\eta_a(h)),\mu(h)(a))=\psi(g)(\psi(h)(d,a)).
\end{align*}

	Let $\tilde{G}=\psi(G)$. Then $\tilde{G}$ is a finite cover of $G^*$. 
	
	We claim that $(\tilde{G},G^*,\pi_2)$ is linked. Let $a\in X^*$ and $A\coloneqq D_a\times \{a\}$. We have to show that $\tilde{G}_{((A))}=\id(A)$. Let $\tilde{g}\in \tilde{G}|_{A}$. By definition this means that there exists some $g\in G$ such that $\psi(g)=\tilde{g}$, that is, $\mu_{\pi_2}(g)=\id(X^*)$ and $\eta_b(g)=B_b$ for all $b\in X^*\setminus a^*$. Let $h\in G$ such that $h|_{X\setminus A}=g|_{X\setminus \pi^{-1}(a)}$ and $h|_{\pi^{-1}(a)}=\id(\pi^{-1}(a))$. Then since $G$ is closed it follows that $h\in G$. Let $g'\coloneqq gh^{-1}$. Then $g'$ fixes every element outside $A$, and $B_a\ni \eta_a(g')=\eta_a(g)$. Therefore, $\eta_a(\tilde{g})=\eta_a(\psi(g))=\phi(\eta_a(g))=\id(A)$, and this is what we needed to prove.
	
	Finally, we show that $G=\lift(\tilde{G},G^*,\data)$ where $\data\colon a\mapsto (F_a,B_a,\phi_a)$. For a permutation $g\in G$ items (1) and (2) of Definition~\ref{def:recover} are clear and item (3) is witnessed by $\psi(g)\in \tilde{G}$. This shows $G\subset \lift(\tilde{G},G^*,\data)$. For the other direction let $g\in \lift(\tilde{G},G^*,\data)$ be arbitrary. As we noted before, the map $\psi\colon g\mapsto \tilde{g}$ is also a surjective homomorphism from $\lift(\tilde{G},G^*,\data)$ to $\tilde{G}$. Thus, we can find some $g'\in G$ such that $\psi(g')=\psi(g)$, that is, $\psi(g(g')^{-1})=\id(\tilde{X})$. In this case however $g(g')^{-1}$ must be contained in the direct product of pointwise binding groups of $B_a: a\in X^*$. Since $G$ is closed, this implies that $gg'\in G$, and thus $g=g(g')^{-1}g'\in G$.
\end{proof}

	We close out this section by an important observation taken from~\cite{EvansIvanovMacpherson} (Lemma 1.4.2). For the sake of completeness we also provide a proof here.

\begin{lemma}\label{closed_factor}
	Let $(G,G^*,\pi)$ be a finite cover. Then if $G$ is closed, then so is $G^*$.
\end{lemma}

	Before showing Lemma~\ref{closed_factor} we first show an easy observation about sequences of permutations which we will also need later.

\begin{lemma}\label{fin_closed}
	Let $(g_i)_{i\in \omega}$ be a sequence of permutations of some countable set $X$ such that for all $x\in X$ the set $F_X:=\{g_i(x)\colon i\in \omega\}$ is finite. Then $(g_i)_{i\in \omega}$ has a convergent subsequence.
\end{lemma}

\begin{proof}
	Let $x_1,x_2,\dots$ be an enumeration of $X$. Then by induction on $j$ we can define subsequences $(g_i^j)_{i\in \omega}$ such that $g_i^j=g_i^{j-1}$ if $i\leq j$ and $|\{g_i^j(x_j)\colon i>j\}|=1$. Then the sequence $h_i:=g_i^i$ is a convergent subsequence of $g_i$. 
\end{proof}

\begin{proof}[Proof of Lemma~\ref{closed_factor}]
	Let $(g_i^*)_{i\in \omega}$ a convergent sequence in $G^*$ with limit $g^*$. We have to show that $g^*\in G^*=\mu(G)$. Let us choose $g_i\in G$ with $g_i/E=g_i^*$. Clearly, the set $\{g_i(x)\colon i\in \omega\}$ is a subset of $\pi^{-1}(\{g_i^*(\pi(x))\colon i\in\omega\})$, in particular, it is finite. Thus, by using Lemma~\ref{fin_closed} we can take a subsequence $(h_i)_{i\in \omega}$ of $(g_i)_{i\in \omega}$ which is convergent. Then we have $$\mu(\lim_{i\rightarrow \infty}h_i)=\lim_{i\rightarrow \infty}(\mu(h_i))=\lim_{i\rightarrow \infty}(g_i^*)=g^*,$$ which finishes the proof of the lemma.
\end{proof}

	Using the results of this section Lemma~\ref{closed_factor} can also be generalized to cover with finite fiber factors.
	
\begin{corollary}\label{closed_factor_fff}
	Let $(G,G^*,\pi)$ be a cover with finite fiber factors. If $G$ is closed, then so is $G^*$.
\end{corollary}

\begin{proof}
	By Lemma~\ref{recover2} we can assume without loss of generality that $G=\lift(\tilde{G},G^*,\data)$ where $\tilde{G}$ and $\data$ are as in Definition~\ref{def:recover}. Then by Lemma~\ref{lift}, item~\ref{it:lift_closed} we can conclude that $\tilde{G}$ is closed, and thus $G^*$ is also closed by Lemma~\ref{closed_factor}.
\end{proof}
\section{Groups with not too fast unlabelled growth}

	We denote by $\subex$ the class of closed groups $G$ such that for no polynomials $p$ do we have $\os_n(G)\geq \frac{2^n}{p(n)}$. As usual, the corresponding class of structures is denoted by $\subexs$. The class $\subexs$ was introduced and discussed in~\cite{simon2025omega}, and it is the main topic of this article. In this paper we refer to the elements of $\subex$ (resp. $\subexs$) as groups (resp. structures) with \emph{not too fast} unlabelled growth.\footnote{The word `slow' is a bit overloaded and usually indicates something slower than what we want to have, for instance polynomial or slower than exponential. Our upper bound is slightly faster than `slower than exponential', this is why we picked the term `not too fast'.} 
	
	Next we recall one of the main results of~\cite{simon2025omega} rephrased in a way that it is more suitable for our discussion.
	
\begin{lemma}\label{simon}
	Let $G\in \subex$. Then there exists a permutation group $G^*$ with $\overline{G^*}\in \gone$, and a cover $(G,G^*,\pi)$ such that for all $a\in \dom(G^*)$
\begin{itemize}
\item either $|\pi^{-1}(a)|=1$ or
\item for every finite subset $B$ of $\dom(G)\setminus \pi^{-1}(a)$ the group $(G_B)|_{\pi^{-1}(a)}$ is a finite cover of a highly set-transitive, but not highly transitive group.
\end{itemize}
	Moreover, for a given $a\in \dom(G^*)$ satisfying item (2) the fibers of the finite cover do not depend on our choice of $B$.
\end{lemma}

\begin{proof}
	  We follow the proof of Theorem 5.1. in~\cite{simon2025omega}. Let $\fm$ be a structure with $\aut(\fm)=G$. We define $F$ to be the equivalence relation of equialgebraicity on elements of order type (i.e. $(a,b)\in F$ iff $a\in \acl_G(\{b\})$ and $b\in \acl_G(\{a\})$) extended by the equality relation on the rest of the structure. 
	
	Then let us consider the equivalence relations $F$ on $M$, and $E$ on $M/F$ as defined in the proof of Theorem 5.1 in~\cite{simon2025omega}. The exact details of the definition of $F$ and $E$ are not important for our argument, we only need that they satisfy the following properties which are argued for in the aforementioned proof.
	
	\emph{$(*)$
\begin{itemize}
\item $F$ has finite classes,
\item the group $\overline{(G/F)/E}$ is the automorphism group of a stable structure, and
\item for all $C\in M/F$ one of the following two conditions holds.
\begin{enumerate}
\item $|C|=1$.
\item for any finite subset $B$ of $M/F\setminus C$ the group $(\aut(M/F)_{B,\{C\}})|_C$ is isomorphic to a highly set-transitive, but not highly transitive group.
\end{enumerate}
\end{itemize}
}

	The unlabelled growth of the group of $(G/F)/E$ is clearly slower than that of the group $G$ which implies that $\overline{(G/F)/E}\in \subex$. Let $\fn$ be a structure with $\aut(\fn)=\overline{(G/F)/E}$. Then $\fn$ is stable and has at most exponential unlabelled growth. By Theorem~\ref{sam} this implies that in fact $\fn\in \gones$, and thus $\overline{(G/F)/E}\in \gone$.

	Now let us define $E'$ to be the pullback of $E$ to $M$, that is, $E'=\{(a,b)\colon (a/F,b/F)\in E\}$. Then we can translate statement $(*)$ for classes of $E'$ as follows.
	
\emph{Let $C'\in M/E'$. Then of the following holds.
\begin{itemize}
\item $|C'|=|\bigcup_{C\in C'} C|=1$.
\item for any finite subset $B'$ of $M\setminus C$ the group $(\aut(M)_{B',\{C'\}})|_{C'}$ is a finite covering of a highly set-transitive, but not highly transitive group.
\end{itemize}
}

	Then by defining $G^*=G/E'$ and $\pi: a\mapsto [a]_E$ we obtain the statement of the lemma. The ``moreover'' part follows from the fact that the fibers as in item (2) are exactly the $F$ classes contained in the corresponding $E$ class.
\end{proof}

	Note that it does not follow directly from the argument presented in~\cite{simon2025omega} that any of the groups $(G_B)|_{\pi^{-1}(a)}$ are closed. In the following we show that not only all these groups, including the fiber groups of the cover $(G,G^*,\pi)$, are closed, but they are in fact finite-index supergroups of the corresponding pointwise binding groups (which are always closed). This implies in particular that the cover $(G,G^*,\pi)$ has finite fiber factors which makes it possible to use the results of Section~\ref{sect:fff}. In particular, we obtain by Corollary~\ref{closed_factor_fff} that the group $G^*$ also needs to be closed.
	
	In order to show the claims above we need the following two lemmas about finite covers of highly set-transitive groups both of which will be proven in Section~\ref{sect:fin_q}. For convenience, we introduce the notation $\fh$ for the class of closed finite covers of highly set-transitive groups. Note that by Lemma~\ref{closed_factor} we know that all groups in $\fh$ are in fact finite covers of some \emph{closed} highly set-transitive group.

\begin{lemma}[Corollary~\ref{no_inf_chain_original}]\label{no_inf_chain}
	There is no infinite strictly descending sequence of groups in $\fh$ with the same fibers.
\end{lemma}

\begin{lemma}\label{factors}
	Let $N,G\in \fh$ be permutation groups with the same fibers so that $N$ is a normal subgroup of $G$. Then $|G:N|$ is finite.
\end{lemma}

	Now we prove a stronger version of Lemma~\ref{simon} that was indicated in our discussion above.
	
\begin{theorem}\label{simon_final}
	Let $G\in \subex$. Then there exists a permutation group $G^*\in \gone$, and a cover $(G,G^*,\pi)$ such that for all $a\in \dom(G^*)$
\begin{enumerate}
\item either $|\pi^{-1}(a)|=1$ or
\item the fiber group $F_a$, and the pointwise binding group $B_a$ at $a$ are both contained in $\fh$, $F_a$ and $B_a$ have the same fibers, and $|F_a:B_a|$ is finite.
\end{enumerate}
	In particular, all fiber groups of $(G,G^*,\pi)$ are closed.
\end{theorem}

	Before proving Lemma~\ref{simon_final} we first show the following easy observation about restrictions of pointwise stabilizers of invariant subsets.

\begin{lemma}\label{tricky_restriction}
	Let $G$ be a closed permutation group, and let $(A,B)$ be a partition of $\dom(G)$ with $G$-invariant subsets. Let us assume that for all finite subsets $T$ of $B$ we have $\overline{(G_T)|_A}=\overline{G|_A}$. Then ${G_B}|_A=G|_A=\overline{G|_A}$.
\end{lemma}

\begin{proof}
	We first show the equality ${G_B}|_A=G|_A$. The ``$\subset$'' containment is obvious. For the other direction let us take an element $h\in G$. We have to show that $h|_A\cup \id(B)\in G$. Since the group $G$ is closed, for this it is enough to show that for all finite sets $S\subset A$ and $T\subset B$ there exists some $h'\in G$ such that $h'|_S=h|_S$ and $h'|_T=\id(T)$. This means in other words that there exists some $h''\in (G_T)|_A$ such that $h''$ and $h'$ agree on $S$. The existence of such a permutation follows from the fact that $h|_A\in G|_A\subset \overline{(G_T)|_A}$.
	
	Finally, since the group $G_B$ is closed, the equality $G|_A=\overline{G|_A}$ follows.
\end{proof}

\begin{proof}[Proof of Theorem~\ref{simon_final}]
	Let $\pi$ and $G$ be as in the conclusion of Lemma~\ref{simon}. Let $a\in X\coloneqq \dom(G^*)$ and $C=\pi^{-1}(a)$. If $|C|=1$ then we are done, so let us assume that this is not the case. Let $A_0\subset A_1\subset \dots$ be finite sets such that $\bigcup_{i=0}^{\infty}(A_i)=X\setminus C$. By our construction, we know that the groups $(G_{A_i,\{C\}})|_C$ are all finite covers of highly set-transitive groups. Thus, by Lemma~\ref{no_inf_chain} we can conclude that the sequence $\overline{(G_{A_i,\{C\}})|_C}$ stabilizes after finitely many steps. This implies that there exists a finite subset $B$ of $X\setminus C$ such that for all finite $S\subset X\setminus C$ we have $\overline{(G_{B,\{C\}})|_C}=\overline{(G_{B\cup S,\{C\}})|_C}$. Thus, we can apply Lemma~\ref{tricky_restriction} to the group $\overline{G_{B,\{C\}}}$ and the partition $(C,X\setminus C)$. We obtain that $G_{((C))}=G_{B,\{C\}}=\overline{G_{B,\{C\}}}$. This in turn implies that $G_{((C))}\in \fh$. Moreover, $G_{((C))}$ and $G_{(C)}$ must have the same fibers. Clearly, $G_{((C))}$ is a normal subgroup in $G_{(C)}$, and thus also in $\overline{G_{(C)}}$. Therefore, by Lemma~\ref{factors} we obtain that $|\overline{G_{(C)}}:G_{((C))}|<\infty$, and thus also $|G_{(C)}:G_{((C))}|<\infty$. In particular, $G_{(C)}$ is closed, and thus $G_{(C)}\in \fh$.
	
	We have shown that $(G,G^*,\pi)$ is a cover with finite fiber factors. Therefore, by Corollary~\ref{closed_factor_fff} the group $G^*$ is also closed.
\end{proof}

	Using Theorem~\ref{simon_final} we will obtain a classification of groups in $\subex$ in the following steps. 

\begin{itemize}
\item\label{it:proof1} We first identify all possible pairs $(B_a,F_a)$ as in item (2) of Theorem~\ref{simon_final}. This boils down to the classification of those pairs $(N,G)$ such that
\begin{enumerate}[(i)]
\item $N,G\in \fh\cap \subex$
\item $N$ and $G$ have the same fibers, and
\item $N\triangleleft G$.
\end{enumerate}

	This will be done in Section~\ref{sect:fin_q}. The reason why the groups $B_a$ and $F_a$ need to be contained in $\subex$ follows from an easy orbit counting argument (see Lemma~\ref{class1}). We will also show that the conditions above are sufficient in the sense that if $G^*\in \gone$, and $G$ is a closed cover of $G^*$ all of whose pointwise binding groups and fiber groups are as in items (i)-(iii) above then $G\in\subex$ (see Section~\ref{sect:finish}).

\item\label{it:proof2} Then we give a description of finite covers of hereditarily cellular groups. This will be done in Section~\ref{sect:fin_h}. In particular, we will see that closed finite covers of hereditarily cellular groups are also hereditarily cellular.

\item\label{it:proof3} Finally, by using the results of Section~\ref{sect:fff} we can recover all groups in $\subex$ from the information provided in the previous items.
\end{itemize}

	The full details of the classification outlined above will be presented in Section~\ref{sect:finish}.

%% file: fincover_q.tex

\section{Finite covers of closed highly set-transitive groups}\label{sect:fin_q}

\subsection{Classification}\label{sect:class_fq}

	Finite covers of reducts of $\mathbb{Q}$ are described in~\cite{ivanov1999finite}, Section 2. We start this section by recalling this description in terms of automorphism group, and as a result, we give a complete characterization of groups in $\fh$.

	We first need to introduce some further notation. Let $r$ be an arbitrary irrational number that will be fixed throughout this section, and let us consider some permutations $\flip$ and $\turn$ as in Notation~\ref{not:ft}.
	
\begin{notation}
	Let $F$ be a finite set, and let $\sigma,\tau\in \sym(F)$. Then we define the permutations $\turn_{\sigma},\flip_{\tau}\in \sym(F)\wr \sym(\mathbb{Q})$ so that
\begin{itemize}
\item $\mu_{\pi_2}(\turn_{\sigma})=\turn$,
\item $\delta_a(\turn_{\sigma})=\sigma$ for $a>r$,
\item $\delta_a(\turn_{\sigma})=\id(F)$ for $a<r$,
\item $\mu_{\pi_2}(\flip_{\tau})=\flip$,
\item $\delta_a(\flip_{\tau})=\tau$ for all $a\in \mathbb{Q}$.
\end{itemize}
\end{notation}

\begin{theorem}\label{fin_covers_q}
	Let $G$ be a closed finite cover of some highly set-transitive group $G^*$. Then we can rename the elements of $\dom(G)$ and $\dom(G^*)$ so that the following hold.
\begin{enumerate}[(i)]
\item $\dom(G^*)=\mathbb{Q}$
\item $\dom(G)=F\times \mathbb{Q}$,
\item $\pi=\pi_2$ is a finite cover,
\end{enumerate}
	and there exist $H\triangleleft L\leq \sym(F)$ such that one of the following holds.
\begin{enumerate}[(a)]
\item\label{it:c_order} $G^*=\aut(\mathbb{Q};<)$ and $G=\kernel{L}{H}(\id(F)\wr \aut(\mathbb{Q};<))$.
\item\label{it:c_betw} $G^*=\aut(\mathbb{Q};\betw)$ and $G=\langle \kernel{L}{H}(\id(F)\wr \aut(\mathbb{Q};<)),\flip_{\tau}\rangle$ for some $\tau\in \sym(F)$.
\item\label{it:c_cyc} $G^*=\aut(\mathbb{Q};\cyc)$ and $G=\langle \kernel{L}{H}(\id(F)\wr \aut(\mathbb{Q};<)),\turn_{\sigma}\rangle$ for some $\sigma\in \sym(F)$.
\item\label{it:c_sep} $G^*=\aut(\mathbb{Q};\sep)$ and $G=\langle \kernel{L}{H}(\id(F)\wr \aut(\mathbb{Q};<)),\flip_{\tau},\turn_{\sigma}\rangle$ for some $\sigma,\tau\in \sym(F)$.
\item\label{it:c_full} $G^*=\sym(\mathbb{Q})$ and $G=\kernel{L}{H}(\id(F)\wr \sym(\mathbb{Q}))$.
\end{enumerate}
\end{theorem}

\begin{proof}
	By Lemma~\ref{fin_closed} we know that $G^*$ is also closed. Thus, by Theorem~\ref{cameron} we can assume that $G^*$ is the automorphism group of one of the five reducts of $(\mathbb{Q};<)$ listed in the same theorem. The finite covers of the full symmetric group are classified in~\cite{ivanov1994some} and in~\cite{bodor2024classification}. From now on we assume that $G^*\neq \sym(\mathbb{Q})$.
	
	By~\cite{ivanov1999finite}, Theorem 2.1 we know that every finite cover of $\aut(\mathbb{Q};<)$ is strongly split which means that the domain of $G$ can be renamed such that $G=\Ker(\pi)(\id(F)\wr \aut(\mathbb{Q};<)$ for some finite set $F$. Furthermore, from the rest of the results of Section 2.1 in~\cite{ivanov1999finite} we know that
\begin{itemize}
\item if $G^*=\aut(\mathbb{Q};\betw)$ then $\flip_{\tau}\in G$ for some $\tau\in \sym(F)$
\item if $G^*=\aut(\mathbb{Q};\cyc)$ then $\turn_{\sigma}\in G$ for some $\sigma\in \sym(F)$, and
\item if $G^*=\aut(\mathbb{Q};\sep)$ then $\flip_{\tau},\turn_{\sigma}\in G$ for some $\sigma,\tau\in \sym(F)$.
\end{itemize}

	Let us define $\phi(\alpha)\coloneqq \hat{\alpha}$ for $\alpha\in \aut(\mathbb{Q};<)$, $\phi(\flip)\coloneqq \flip_{\tau}$ if $\flip\in G^*$ and $\phi(\turn)\coloneqq \turn_{\sigma}$ if $\turn\in G^*$. Then we can extend $\phi$ to $G^*$ in a straightforward way so that $\pi_2\circ \phi=\id$.
	
	Finally, by the results of Section 2.2 in~\cite{ivanov1999finite} we know that $\Ker(\pi_2)$ can be written as $\kernel{L}{H}$ for some $H\triangleleft L\leq \sym(F)$. Considering that the group $G$ can be written as $G=\Ker(\pi_2)\phi(G^*)$, the statement of the theorem follows in all cases.
\end{proof}

	The classification given in Theorem~\ref{fin_covers_q} immediately implies the following.
	
\begin{corollary}\label{no_inf_chain_original}
	There is no infinite strictly descending sequence of groups in $\fh$ with the same fibers.
\end{corollary}

\begin{proof}
	Let us suppose for contradiction that there exists an infinite sequence $G_0\supsetneq G_1\supsetneq$ of groups in $\fh$ with the same fibers. By Theorem~\ref{fin_covers_q} it follows easily that the groups $G_0,G_1,\dots$ can only realize finitely isomorphism types. Therefore, there exist $G_i\subsetneq G_j$ which are isomorphic. In this case, however, for all $n\in \omega$ we have $\orb_n(G_i)=\orb_n(G_j)$ and since these are oligomorphic we have in fact $\orbits_n(G_i)=\orbits_n(G_j)$. Since the groups $G_i$ and $G_j$ are closed this implies $G_i=G_j$, a contradiction.
\end{proof}

	For the rest of this subsection we assume that the groups $G,G^*,H,L$ are as in the conclusion of Theorem~\ref{fin_covers_q}.
	
\begin{lemma}\label{about_flips}
	Let us assume that $\flip_{\tau},\flip_{\tau'}\in G$. Then the following hold.
\begin{enumerate}
\item\label{it:square} $\tau^2\in L$.
\item\label{it:2index} either $\tau\in L$ or $\tau L$ is an index 2 supergroup of $L$.
\item\label{it:2coset} $\tau L=\tau'L$.
\item\label{it:flip_normal} $H\triangleleft \tau L$.
\end{enumerate}
\end{lemma}

\begin{proof}
	Let $\alpha:=\flip_{\tau'}^{-1}\flip_{\tau}$. Then $\alpha=\Ker(\pi)$, thus $\alpha\in \kernel{L}{H}$. On the other hand, we have $\delta_a(\alpha)=(\tau')^{-1}\tau$ for all $a\in \mathbb{Q}$. Therefore, $(\tau')^{-1}\tau\in L$ which implies item~\ref{it:2coset} immediately.
	
	Item~\ref{it:square} follows immediately from item~\ref{it:2coset} applied in the case when $\tau'=\tau^{-1}$. Item~\ref{it:2index} is clear from item~\ref{it:square}.
	
	Finally, let $a\in \mathbb{Q}$ and $\sigma\in H$ be arbitrary. Then $$\kernel{L}{H}\ni\, \flip_{\tau}^{-1}\!\kappa(\sigma,a)\flip_{\tau}=\kappa(\tau^{-1}\sigma\tau, \flip\!(a)).$$ Therefore, $\tau^{-1}\sigma\tau\in H$, that is, $\tau$ normalizes $H$. We already know that $H\triangleleft L$, thus item~\ref{it:flip_normal} follows.
\end{proof}

	Lemma~\ref{about_flips} implies that $\tilde{L}\coloneqq \tau L$ is an index 1 or 2 supergroup of $L$ group which, together with $H$ and $L$, already determine $G$ in the case when $G^*=\aut(\mathbb{Q};\betw)$.
	
\begin{lemma}\label{about_turns}
	Let us assume that $\turn_{\sigma},\turn_{\sigma'}\in G$. Then the following hold.
\begin{enumerate}
\item\label{it:turn} $\sigma\in L$
\item\label{it:turn_cosets} $\sigma H=\sigma'H$.
\item\label{it:turn_comm} $\sigma^{-1}\ell^{-1}\sigma \ell\in H$ for all $\ell\in L$.
\end{enumerate}
\end{lemma}

\begin{proof}
	(\ref{it:turn}) It is easy to see that $\turn_{\sigma}^2=\kappa(\sigma,\mathbb{Q})$. This implies immediately that $\sigma\in L$.
	
	(\ref{it:turn_cosets}) Let $\alpha:=\turn_{\sigma'}^{-1}\turn_{\sigma}$. Then $\alpha=\kappa(\sigma^{-1}\sigma,\mathbb{Q}(>r))$, and therefore $\sigma^{-1}\sigma\in H$.
	
	Now let us assume that $\ell\in L$. Let $\beta\coloneqq \kappa(\ell,\mathbb{Q})\in \kernel{L}{H}\leq G$. Then we have $G\ni \beta^{-1}\!\turn_{\sigma}\!\beta=\turn_{\ell^{-1}\sigma \ell}$. Thus, by item~\ref{it:turn_cosets} we obtain that $\ell^{-1}\sigma \ell=\sigma H$, and thus item~\ref{it:turn_comm} holds. 
\end{proof}

	Lemma~\ref{about_turns} implies that in the case when $G^*=\aut(\mathbb{Q};\cyc)$ then the groups $H$, $L$, and the coset $\sigma H$ already determine the group $G$.
	
	The following observation will be useful in our later discussions.
	
\begin{lemma}\label{fix_some}
	Let $S\subset \mathbb{Q}$ with $3\leq |S|<\infty$. We write $G_S$ for the group $\{g\in G\colon \mu_{\pi_2}(g)\in G^*_S\}$. Then $(G_S)|_{\pi^{-1}(S)}=\kernel{L}{H}|_{\pi^{-1}(S)}$.
\end{lemma}

\begin{proof}
	Let $g\in G_S$. Since $g^*:=\mu_{\pi_2}(g)$ fixes 3 elements of $\mathbb{Q}$ it follows that $g^*$ cannot be contained in $\tau\aut(\mathbb{Q};\betw)$. Thus, we can assume without generality that $G^*=\aut(\mathbb{Q};R)$ for some $R\in \{<,\cyc,=\}$. In these cases it follows from the definition of the groups listed in Theorem~\ref{fin_covers_q} that $\delta_s(g)\in L$ for all $s\in \mathbb{Q}$. It remains to show that $\delta_s(g)H=\delta_t(g)H$ for all $s,t\in S$. By Theorem~\ref{fin_covers_q} this is again clear if $G^*=\aut(\mathbb{Q},<)$ or $G^*=\sym(\mathbb{Q})$. We are left with the case when $G^*=\aut(\mathbb{Q},\cyc)$ in which case $$G=\langle \kernel{L}{H}(\id(F)\wr \aut(\mathbb{Q};<)),\turn_{\sigma}\rangle$$ for some $\sigma\in L$ by Theorem~\ref{fin_covers_q}.
	
	Let $u=(u_1,\dots,u_k)$ be a tuple listing the elements of $S$, and let
\begin{align*}
O_1:=&\{((\ell(u_1),s),\dots,(\ell(u_k),s),(\ell h(u_1),t),\dots,(\ell h(u_k),t))\colon \ell\in L,h\in H,s,t\in \mathbb{Q}, s<t\},\\
O_2:=&\{((\ell(u_1),s),\dots,(\ell(u_k),s),(\sigma \ell h(u_1),t),\dots,(\sigma \ell h(u_k),t))\colon \ell\in L,h\in H,s,t\in \mathbb{Q}, s>t\},\\
O:=&O_1\cup O_2.
\end{align*}

	Then it follows from the construction of $G$ that $O$ is a $2k$-orbit of $G$. Now let $s,t\in S$. Then clearly $G_S$ cannot map a tuple $O_1\cap \pi^{-1}(\{s,t\})$ to $O_2$. Hence, $G_S$ preserves $O_1$ which means exactly that $\delta_s(g)^{-1}\delta_t(g)\in H$.
\end{proof}

\subsection{Normal subgroups}\label{sect:normal}

	In this subsection we examine all those pairs $(N,G)$ such that
\begin{enumerate}
\item $N,G\in \fh$
\item $N\triangleleft G$, and
\item $N$ and $G$ have the same fibers.
\end{enumerate}

	In particular we show that for all such pairs of groups the index $|G:N|$ is finite.
	
	So let us consider some pair $(N,G)$ as above which we will fix throughout this subsection. By renaming the elements of $\dom(N)$ we can assume that there exists some $N^*\in \aut(R)$ for some $R\in \{<,\betw,\cyc,\sep,=\}$ such that the pair $(N^*,N)$ is as in the conclusion of Theorem~\ref{fin_covers_q}. Let $G^*:=\mu(G)=\mu_{\pi_2}(G)$. Then $(G^*,G,\pi_2)$ is also a finite cover, and it is also as in the conclusion of Theorem~\ref{fin_covers_q}. For the rest of the subsection we assume that the $G,G^*,N,N^*,\pi$ are as above.
	
	We start with a few observations.
	
\begin{lemma}\label{image_normal}
	$N^*\triangleleft G^*$.
\end{lemma}

\begin{proof}
	Let $n^*\in N^*$ and $g^*\in G^*$. Then there exist $n\in N$ and $g\in G$ with $\mu(n)=n^*$ and $\mu(g)=g^*$. Then $g^{-1}ng\in N$, so $(g^*)^{-1}n^*g^*=\mu(g^{-1}ng)\in N^*$.
\end{proof}

	By Lemma~\ref{normal_highly} we know that the groups $G^*$ and $N^*$ are as in the conclusion of Lemma~\ref{normal_highly}. By Theorem~\ref{fin_covers_q} we know that $\Ker(\mu|_G)=\kernel{L}{H}$ and $\Ker(\mu|_N)=\kernel{M}{H'}$ for some $H'\leq H,H'\triangleleft M, H\triangleleft L\leq \sym(F)$.

\begin{proposition}\label{kernel}
	$H=H'$.
\end{proposition}

\begin{proof}
	Clearly, $\ker(\mu|_N)=\ker(\mu|_G)\cap N$ is a normal subgroup of $\ker(\mu|_G)$. Let $s\in \mathbb{Q}$ and $\alpha\in \aut(\mathbb{Q};<)$ such that $\alpha(s)\neq s$, and let $h\in H$ be arbitrary. Then $\hat{\alpha}\in N$ and $\kappa(h,s)\in \kernel{L}{H}\leq G$. Hence, $\beta:=\hat{\alpha}^{-1}\kappa(h,s)\hat{\alpha}\in N$. Now $\delta_s(\beta)=h$ and $\delta_t(\beta)=\id(F)$ for all $t\neq s,\alpha(s)$. This implies that $h$ must be contained in $H'$. Therefore, $H=H'$.
\end{proof}

	By Lemma~\ref{normal_highly} we know that $\turn\in N^*$ if and only if $\turn\in G^*$. In this case we know by Theorem~\ref{fin_covers_q} that $\turn_{\sigma}\in N^*$ for some $\sigma\in M$. Since $\turn_{\sigma}$ is also contained in $G$ we can assume that it is also the permutation which appears in either item~\ref{it:c_cyc} or item~\ref{it:c_sep} in the definition of the group $G$. 
	
	Let $$\tilde{L}\coloneqq \{\ell \in \sym(F)\colon \exists \alpha\in G\,\forall a\in \mathbb{Q}(\delta_a(\alpha)=\ell)\}.$$ It follows from Theorem~\ref{fin_covers_q} that in cases~\ref{it:c_order}, \ref{it:c_cyc} and \ref{it:c_full} we have $\tilde{L}=L$, and in the remaining cases if $\flip_{\tau}\in G$ then $\tilde{L}=\tau L$. The group $\tilde{M}$ is defined similarly but with the group $N$.
	
\begin{lemma}\label{factors}
	$\tilde{M}$ is a normal subgroup of $\tilde{L}$, and
\begin{enumerate}
\item $G/N\simeq \tilde{L}/\tilde{M}\times Z_2=L/M\times Z_2$ if $|G^*:N^*|=2$ and $\tilde{L}=L$.
\item $G/N\simeq \tilde{L}/\tilde{M}$ otherwise. 
\end{enumerate}
	In particular $|G:N|$ is finite.
\end{lemma}

\begin{proof}
	The fact that $\tilde{M}$ is a normal subgroup of $\tilde{L}$ is clear from the definition. Let $\alpha\in G$ be arbitrary. Then if $a,b\in \mathbb{Q}$ then one can easily check that $\delta_{a}(\alpha) H=\delta_{b}(\alpha) H$. Therefore, the homomorphism $\delta_{\bullet}\colon G\rightarrow \tilde{L}/H$ is well-defined. Let $N':=(\delta_{\bullet})^{-1}(\tilde{M}/H)$. Then $G/N'\simeq (\tilde{L}/H)/(\tilde{M}/H)\simeq \tilde{L}/\tilde{M}$ by the third isomorphism theorem for groups. In the case when $G^*=N^*$ it is clear that $N'=N$, so in this case we are done.
	
	Now let us assume that $|G^*:N^*|=2$. In this case $\flip_{\tau}\in G\setminus H$ for some $\tau\in \tilde{L}\setminus L$. If $\tau\not\in L$, that is, $|\tilde{L}:L|=2$, then we have $N=N'$ again. If $\tau\in L$, that is, $\tilde{L}=L$ then $\mu_{\pi}(N')=G^*$, and $N$ can be written as $N=N'\cap \pi^{-1}(\aut(\mathbb{Q};\cyc))$. In other words $N$ is the kernel of the surjective homomorphism $$\delta_{\bullet}'\coloneqq G\rightarrow \tilde{L}/\tilde{M}\times G^*/N^*, \alpha\mapsto (\delta_{\bullet}(\alpha)\tilde{M},\mu_{\pi}(\alpha)N^*).$$ Therefore, by the homomorphism theorem for groups it follows that $G/N\simeq \tilde{L}/\tilde{M}\times G^*/N^*\simeq \tilde{L}/\tilde{M}\times Z_2$.
\end{proof}

	The following theorem summarizes the results of this subsection.
	
\begin{lemma}\label{factors_coverse}
	Let $N,G$ be permutation groups with the same domain. Then the following are equivalent.
\begin{enumerate}
\item\label{it:normal1}
\begin{itemize}
\item $N,G\in \fh$
\item $N\triangleleft G$, and
\item $N$ and $G$ have the same fibers.
\end{itemize}

\item\label{it:normal2} There exist some permutation groups $H,L,M$ with some finite domain $F$ such that $H\triangleleft M,L$ and $M\triangleleft L$, and we can rename the elements of $\dom(G)$ such that one of the following holds.
\begin{enumerate}
\item $N=\kernel{M}{H}(\id(F)\wr \aut(\mathbb{Q};R))$ and $G=\kernel{L}{H}(\id(F)\wr \aut(\mathbb{Q};R))$ with $R\in \{<,=\}$.
\item $N=\kernel{M}{H}(\id(F)\wr \aut(\mathbb{Q};<))$ and\\ 
$G=\langle \kernel{L}{H}(\id(F)\wr \aut(\mathbb{Q};<)),\flip_{\tau}\rangle$.
\item $N=\langle \kernel{M}{H}(\id(F)\wr \aut(\mathbb{Q};<)),\flip_{\tau}\rangle$ and\\
$G=\langle \kernel{L}{H}(\id(F)\wr \aut(\mathbb{Q};<)),\flip_{\tau}\rangle$.
\item $N=\langle \kernel{M}{H}(\id(F)\wr \aut(\mathbb{Q};<)),\turn_{\sigma}\rangle$ and\\ 
$G=\langle \kernel{L}{H}(\id(F)\wr \aut(\mathbb{Q};<)),\turn_{\sigma}\rangle$.
\item $N=\langle \kernel{M}{H}(\id(F)\wr \aut(\mathbb{Q};<)),\turn_{\sigma}\rangle$ and\\
$G=\langle \kernel{L}{H}(\id(F)\wr \aut(\mathbb{Q};<)),\turn_{\sigma},\flip_{\tau}\rangle$ 
\item $N=\langle \kernel{M}{H}(\id(F)\wr \aut(\mathbb{Q};<)),\turn_{\sigma},\flip_{\tau}\rangle$ and\\
$G=\langle \kernel{L}{H}(\id(F)\wr \aut(\mathbb{Q};<)),\turn_{\sigma},\flip_{\tau}\rangle$ 
\end{enumerate}
	where
\begin{itemize}
\item $\sigma\in M$ with $\{\sigma^{-1}m^{-1}\sigma m\colon m\in M\}\subset H$, and
\item either $\tau\in L$; or $H,L\triangleleft \tau L$ and $|\tau L:L|=2$.
\end{itemize}
\end{enumerate}
\end{lemma}

\begin{proof}
	The implication (\ref{it:normal1})$\rightarrow$ (\ref{it:normal2}) follows directly from the argument presented in this subsection. The converse is a straightforward calculation.
\end{proof}

\subsection{Finite covers with not too fast unlabelled growth}\label{sect:qslow}

	In this subsection we classify those closed finite coverings of highly set-transitive groups which are also contained in the class $\subex$, and we describe when we have a normal subgroups relation between two such groups with the same fibers. We will see that a group $G\in \fh$ is contained in $\subex$ if and only if its pointwise binding group is highly transitive.

	We start by introducing a sequence $\subex_d: d\in \mathbb{N}=\{1,2,\dots\}$ of subclasses of $\subex$ by imposing some stronger, but still exponential, upper bounds on the unlabelled profile. Within the class $\fh$ these classes correspond exactly to those groups whose fiber size is at most $d$.
		
\begin{definition}
	Let $d\in \mathbb{N}$. Then we denote by $\gamma_d$ the largest real root of the polynomial $f_d(x):=x^d-x^{d-1}-\dots-x-1$, and we denote by $\subex_d$ the class of those closed groups $G$ such that for all $\varepsilon>0$ we have $\os_n(G)<(\gamma_d+\varepsilon)^n$ if $n$ is large enough.
\end{definition}

\begin{example}
	$\gamma_1=1$ and $\gamma_2=\frac{1+\sqrt{5}}{2}\approx 1.618$, the golden ratio. The next values of this sequence are approximately $\gamma_3\approx 1.839, \gamma_4\approx 1.928, \gamma_5\approx 1.966.$
\end{example}

\begin{lemma}\label{two}
	The sequence $(\gamma_d)_d$ is strictly increasing and $\lim_{d\rightarrow \infty}\gamma_d=2$.
\end{lemma}

\begin{proof}
	Clearly, $\gamma_d>0$ since $f_d(0)=-1<0$. Since $f_d(\gamma_d)=0$ we have $f_{d+1}(\gamma_d)=\gamma_nf_n(\gamma_n)-1=-1<0$. Therefore, $\gamma_{n+1}$ must be greater than $\gamma_n$. 

	Put $\gamma=\lim_{d\rightarrow \infty}\gamma_d$. If $t\geq 2$ then we can easily show by induction that $f_n(t)\geq 1$. Indeed, $f_1=1$ and $f_{n+1}(t)=tf_n(t)-1>2-1=1$. Therefore, $\gamma_d<2$, and thus $\gamma\leq 2$. On the other hand, if $p<2$ then $$\sum_{k=1}^\infty\Bigl({\frac{1}{p}}\Bigr)^k=\frac{1}{p(1-\frac{1}{p})}=\frac{1}{p-1}>1,$$ and thus $\sum_{k=1}^d({\frac{1}{p}})^k>1$ for some $d\in \omega$. In this case we have $$f_d(p)=p^n\Bigl(1-\sum_{k=1}^d\Bigl({\frac{1}{p}}\Bigr)^k\Bigr)<0,$$ and hence $\gamma_d>p$. This proves that $\gamma=2$.
\end{proof}

	It is clear from the definition that $\subex_1\subset \subex_2\subset \dots$, and by Lemma~\ref{two} we have $\bigcup_{d\in \omega}\subex_d\subset \subex$. Later we will see that in fact the equality holds.
	
	For convenience, we sometimes write $\subex_{\infty}$ for $\subex$.

\begin{lemma}\label{adding_constants}
	The classes $\subex_d\colon d\in \mathbb{N}\cup \{\infty\}$ are closed under taking stabilizers of finite sets.
\end{lemma}

\begin{proof}
	For $d=\infty$ this is exactly Lemma 2.2 in~\cite{simon2025omega}, and for $d<\infty$ essentially the same proof works.
\end{proof}

	Our next goal is to give a description of all groups in $\fh\cap \subex_d$ for $d\in \mathbb{N}\cup \{\infty\}$. 
	
	Let us fix a group $G\in \fh$. Then we can rename the elements of $\dom(G)$ such that there exists a finite cover $(G^*,G,\pi)$ and groups $H\triangleleft L\leq \sym(F)$ as in the conclusion of Theorem~\ref{fin_covers_q}. We also use the notation from Subsection~\ref{sect:class_fq}. Next we show that $G\in \subex_d$ if and only if $H$ is highly set-transitive and $|F|\leq d$. We know that all finite coverings of $\sym(\omega)$ have a polynomial orbit growth (see~\cite{falque2020classification}). In particular all these groups are contained in $\subex_1$. This handles the case when $G^*=\sym(\mathbb{Q})$. Next we examine what happens in the other cases.
	
	We first need the following general observation about the unlabelled growth of finite-index subgroups.
	
\begin{lemma}\label{fin_index_not_too_fast}
	Let $H\leq G\leq \sym(X)$. Then $\os_n(H)\leq |G:H|\os_n(G)$.
\end{lemma}

\begin{proof}
	Let $\gamma_i:i\in I$ be such that $G=\bigsqcup_{i\in I}\gamma_iH$. If the $A_j\in {X \choose n}:j\in J$ represent all orbits of the action $G\curvearrowright {X \choose n}$ then the sets $\gamma_i(A_j)$ represent all orbits of the action $H\curvearrowright {X \choose n}$. Therefore, the inequality $\os_n(H)\leq |G:H|\os_n(G)$ follows.
\end{proof}

	The following is a straightforward consequence of Lemma~\ref{fin_index_not_too_fast}.

\begin{corollary}\label{subex_finex}
	Let $G$ be a permutation group, and let $H$ be a finite-index subgroup of $G$. Then $G\in \subex_d$ if and only if $H\in \subex_d$.
\end{corollary}
	
\begin{lemma}\label{h_enough}
	Let us assume that $G^*\neq \sym(\mathbb{Q})$, and let $s_1,s_2,s_3\in \mathbb{Q}$ with $s_1<s_2<s_3$. Then $G_{F\times \{s_1,s_2,s_3\}}$ fixes the set $F\times (s_1,s_2)$, and the restriction $$K:=(G_{F\times \{s_1,s_2,s_3\}})|_{F\times (s_1,s_2)}$$ is isomorphic to $\kernel{H}{H}(\id(F)\wr \aut(\mathbb{Q};<))=H\wr \aut(\mathbb{Q};<).$
\end{lemma}

\begin{proof}
	By Theorem~\ref{reducts_of_q} we know that $G^*$ preserves the relation $\sep$. It is easy to check that $(t,s_2,s_3,s_1)\in \sep$ if and only if $s_1<t<s_2$. Moreover, $$K^*\coloneqq (G^*_{s_1,s_2,s_3})|_{s_1,s_2,s_3}=\aut((s_1,s_2);<)\simeq \aut(\mathbb{Q};<).$$ This implies that $G_{F\times \{s_1,s_2,s_3\}}\leq \mu_{\pi}^{-1}(G^*_{s_1,s_2,s_3})$ preserves the set $F\times (s_1,s_2)$ and $\mu(K)|_{(s_1,s_2)}\leq K^*$. We know that $G$ contains the group $\kernel{H}{H}(\id(F)\wr \aut(\mathbb{Q};<))$ from which we can easily conclude that $K$ contains both $\kernel{H}{H}|_{F\times (s_1,s_2)}$ and $(\id(F)\wr \aut((s_1,s_2);<))$. Therefore, $K^*$ is a finite cover of $K$. The only thing we have to check that the kernel of this cover is $\kernel{H}{H}|_{F\times (s_1,s_2)}$, that is, $\delta_a(\alpha)\in H$ for all $\alpha\in K$ and $s_1<a<s_2$. By our definition $\alpha$ extends to and element of $G$ fixing $F\times \{s_1,s_2,s_3\}$ pointwise. But then $\delta_{s_1}(\alpha)=\id$ which implies $\delta_a(\alpha)\in H$.
\end{proof}

\begin{corollary}\label{gnull}
	Let us assume that $G^*\neq \sym(\mathbb{Q})$, and let
\[
G_0\coloneqq H\wr \aut(\mathbb{Q};<)).
\]
	Then for all $d\in \mathbb{N}\cup \{\infty\}$ we have $G\in \subex_d$ if and only if $G_0\in \subex_d$.
\end{corollary}

\begin{proof}
	The ``if'' direction is obvious since $G_0\leq G$. The other direction follows from Lemmas~\ref{adding_constants} and \ref{h_enough}.
\end{proof}

	Now let us fix the group $G_0$ as in Corollary~\ref{gnull}. Then we can write the following recursion formula for the unlabelled growth of $G_0$.
	
\begin{lemma}\label{recursion}
	$\os_n(G_0)=\sum_{i=1}^{|F|}{\os_i(H)\os_{n-i}(G_0)}$. 
\end{lemma}

\begin{proof}
	Let $A$ be an $n$-element subset of $F\times \mathbb{Q}$. Let $m$ the maximum of those numbers which appear as the second coordinate of an element of $A$, and let $A_{\max}\coloneqq \{a\in A\colon \pi(a)=m\}$. Then the $G_0$-orbits of $A_{\max}$ and $A\setminus A_{\max}$ can be recovered from the orbit of $A$. Conversely, by the construction of $G_0$ it follows that the orbit of $A$ is uniquely determined by the orbits of $A_{\max}$ and $A\setminus A_{\max}$. Whenever $|A_{\max}|=i$ then we have $\os_i(H)$ many choices for the orbit of $A_{\max}$ and $\os_{n-i}(G_0)$ many choices for the orbit of $A\setminus A_{\max}$. Hence, the recursion formula follows.
\end{proof}

	The following statement follows from the resolution of linear recursion. For the sake of completeness, we also give a direct proof.
	
\begin{lemma}\label{solve_recursion}
	Let $(a_n)_n$ be a sequence of positive numbers satisfying the recursion $a_{n+k}:=\sum_{i=0}^{k-1}c_ia_{n+i}$ for some $c_i>0$. Let $q$ be the largest real root of $f(x):=x^k-\sum_{i=0}^{k-1}c_ix^i$. Then there exist $C_1,C_2>0$ such that $C_1q^n<a_n<C_2q^n$ for all $n\in \mathbb{N}$.
\end{lemma}

\begin{proof}
	We can clearly find $C_1,C_2>0$ such that $C_1q^n<a_n<C_2q^n$ holds for $n\leq k$. We show by induction that this choice suffices.
	
	Let us suppose that $C_1q^m<a_m<C_2q^m$ holds for $m<n$. Then we have
\[
a_n=\sum_{i=0}^{k-1}c_ia_{n-k+i}\geq C_1q^{n-k}\sum_{i=0}^{k-1}c_iq^i=C_1q^{n-k}(q^k-f(q))=C_1q^n.
\]
	The other inequality follows similarly.
\end{proof}

\begin{lemma}\label{when_subex}
	Let us assume that $G^*\neq \sym(\mathbb{Q})$, and let $d\in \mathbb{N}\cup \{\infty\}$. Then $G\in \subex_d$ if and only if $H$ is highly set-transitive and $|F|\leq d$. Moreover, if $d\geq 2$ and $G\in \subex_d\setminus \subex_{d-1}$ then $\os_n(G)\geq C\gamma_d^n$ for some $C>0$.
\end{lemma}

\begin{proof}
	$G_0$ be as in Corollary~\ref{gnull}, and let $G_0^*\coloneqq \mu_{\pi}(G_0)$. By Corollary~\ref{gnull} it is enough to prove the statement for the group $G_0$.
	
	Let us define the polynomial $f(x):=x^{|F|}-\sum_{i=0}^{|F|-1}\os_i(H)x^i$, and let $q$ be the largest real root of $f$. Then by Lemmas~\ref{recursion} and~\ref{solve_recursion} we can conclude that $\os_n(G_0)=\Theta(q^n)$. Now if $H$ is not highly set-transitive, then $\os_j(H)\geq 2$ for some $0\leq j\leq k-1$. In this case $$f(2)=2^n-\sum_{i=0}^{|F|-1}\os_i(H)2^i\leq 2^n-\sum_{i=0}^{|F|-1}2^i-2^j=1-2^j\leq 0.$$ Therefore, $q\geq 2$ and thus $G_0\not\in \subex$. On the other hand, if $H$ is highly set-transitive then $f(x)=x^{|F|}-\sum_{i=0}^{|F|-1}x^i$, and thus $q=\gamma_{|F|}$. Therefore, $G_0\in \subex_d$ if and only if $\gamma_{|F|}\leq \gamma_d$, that is, $|F|\leq d$.
	
	The ``moreover'' part of the lemma is also clear from the calculation above.
\end{proof}

\begin{remark}
	By solving the recursion in Lemma~\ref{solve_recursion} one can easily conclude that in fact $u_n(G_0)\sim C\gamma_{|F|}^n$ for some $C$. We do not need this fact in our discussion. 
\end{remark}

	Highly set-transitive groups of finite degree are classified in~\cite{livingstone1965transitivity}.
	
\begin{theorem}\label{hh_finite}
	Let $K$ be a highly set-transitive group with degree $n\in \omega$. Then one of the following holds.
\begin{enumerate}
\item $K$ is the full symmetric group $S_n$.
\item $n\geq 3$ and $K$ is the alternating group $A_n$.
\item $n=5$ and $K$ is the affine group $\agl(1,5)$.
\item $n=9$ and $K$ is the projective linear group $\pgl(2,8)$.
\item $n=9$ and $K$ is the projective semilinear group $\pgammal(2,8)$.
\end{enumerate}
\end{theorem}

	Now let us assume that $H$ is one of the groups above. Since $\tilde{L}$ contains $H$ it follows that $\tilde{L}$ is also one of the groups on the list. By item~\ref{it:flip_normal} of Lemma~\ref{about_flips}  we also know that $H$ is a normal subgroup of $\tilde{L}$. This is only possible in the following cases. $(\dag)$.

\begin{enumerate}
\item $H=\tilde{L}$;
\item $H=A_n, \tilde{L}=S_n$ in which case $|\tilde{L}:H|=2$. $(*)_1$
\item $H=\pgl(2,8), L=\tilde{L}=\pgammal(2,8)$ in which case $|\tilde{L}:H|=|L:H|=3$. $(*)_2$
\end{enumerate}

	Considering this and Lemma~\ref{fin_covers_q} we can give a list of every group in $\fh\cap \subex$.
	
\begin{theorem}\label{summary_fcover}
	Let $G\in \fh\cap \subex$. Then $G$ is isomorphic to one of the following groups.
\begin{enumerate}[(i)]
\item $\kernel{L}{H}(\id(F)\wr \sym(\mathbb{Q})$ for some $H\triangleleft L\leq \sym(F)$.
\item $\kernel{H}{H}(\id(F)\wr \aut(\mathbb{Q};R))$ where $H$ is one of the groups listed in Theorem~\ref{hh_finite} and $R\in \{<,\betw,\cyc,\sep\}$.
\item $\kernel{L}{H}(\id(F)\wr \aut(\mathbb{Q};R))$ where the groups $H$ and $L$ are as in $(*)_1$ or $(*)_2$ and $R\in \{<,\betw,\cyc,\sep\}$.
\item $\kernel{A_n}{A_n}(\id(F)\wr \aut(\mathbb{Q};R))\!\!\flip_{\tau}$ where $\tau$ is an odd permutation and $R\in \{<,\cyc\}$.
\item $\kernel{L}{H}(\id(F)\wr \aut(\mathbb{Q};R))\!\turn_{\sigma}$ where the groups $H$ and $L$ are as in $(*)_1$ or $(*)_2$, $\sigma\in L\setminus H$ and $R\in \{<,\betw\}$.
\end{enumerate}
	Conversely, all the groups listed above are in $\fh\cap \subex$.
\end{theorem}

	Using the results of Subsection~\ref{sect:normal} we can easily determine in which cases two groups from the list as above are in a normal subgroup relationship.
	
\begin{lemma}\label{normal_subex}
	Let $N$ and $G$ be groups as in items (i)-(v) in Theorem~\ref{summary_fcover}. Then $N\triangleleft G$ if and only if one of the following holds.
\begin{enumerate}[(i)]
\item\label{it:n_equal} $N=G$.
\item\label{it:n_sym} $N=\kernel{M}{H}(\id(F)\wr \sym(\mathbb{Q})), G=\kernel{L}{H}(\id(F)\wr \sym(\mathbb{Q}))$ and $M\triangleleft L$.
\item\label{it:n_order} $N=\kernel{H}{H}(\id(F)\wr \aut(\mathbb{Q};R)), G=\kernel{L}{H}(\id(F)\wr \aut(\mathbb{Q};R))$ where the groups $H$ and $L$ are as in $(*)_1$ or $(*)_2$ and $R\in \{<,\betw,\cyc,\sep\}$.
\item\label{it:n_flip} $N=\kernel{M}{H}(\id(F)\wr \aut(\mathbb{Q};R)), G=\kernel{L}{H}(\id(F)\wr \aut(\mathbb{Q};R))\flip_{\id}$ where 
\begin{itemize}
\item either $H=L$ or $H$ and $L$ are as in $(*)_1$ or $(*)_2$
\item $M\in \{H,L\}$, and
\item $R\in \{<,\cyc\}$.
\end{itemize}
\item\label{it:n_an} $N=\kernel{A_n}{A_n}(\id(F)\wr \aut(\mathbb{Q};R))$ and $G=\kernel{A_n}{A_n}(\id(F)\wr \aut(\mathbb{Q};R))\!\!\flip_{\tau}$ where $\tau$ is an odd permutation and $R\in \{<,\cyc\}$.
\end{enumerate}
\end{lemma}

\begin{proof}
	It follows from Lemma~\ref{factors_coverse} that for all pairs of groups $(N,G)$ satisfying any of the items above we have $N\triangleleft G$.
	
	For the other direction, let us assume that $N\triangleleft G$.

	We already handled the case when $G^*=\sym(\mathbb{Q})$, so let us assume that this is not the case.
	
	We are using the notation of Subsection~\ref{sect:normal} for the groups $N$ and $G$. Then we have $H\leq M\leq L\leq \tilde{L}$ where $H\triangleleft L$ and $|L:\tilde{L}|=2$. Considering that $H$ and $\tilde{L}$ are as $(\dag)$ it follows that $|\tilde{L}:H|\leq 3$, so either $L=H$ or $\tilde{L}=L$.
	
	\emph{Case 1. $\tilde{L}=H$.} In this case we have $N=G$.
	
	\emph{Case 2. $\tilde{L}=L\neq H$.} In this case we have either $M=H$ or $M=L$. If $M=H$ then depending on whether $|G^*:N^*|$ is 1 or 2 either item~\ref{it:n_order} or item~\ref{it:n_flip} holds. If $M=L$ then $\flip_{\id}$ must be contained in $G$, otherwise we would have $N=G$. In this case item~\ref{it:n_flip} holds again.
	
	\emph{Case 3. $\tilde{L}\neq L$.} In this case $|\tilde{L}:L|=|G^*:N^*|=2$ and $L=M=H$. This is only possible if $\tilde{L}=S_n$ and $L=M=H=A_n$, and then item\ref{it:n_an} holds.
\end{proof}

\begin{corollary}\label{factors_options}
	Let $N$ and $G$ be groups as in items (i)-(v) in Theorem~\ref{summary_fcover}. Let us assume that $N\triangleleft G$ and $G^*\neq \sym(\mathbb{Q})$. Then $|G:N|$ is isomorphic to one of the following groups: $Z_1,Z_2,Z_3,Z_2\times Z_2,Z_3\times Z_2\simeq Z_6$.
\end{corollary}

\begin{proof}
	We have seen already that $|L:H|\in \{1,2,3\}$. Thus, the statement of corollary follows from Lemma~\ref{factors}.
\end{proof}

	The results of this subsection are summarized in the following theorem.
	
\begin{theorem}\label{fs_normal}
	Let $N$ and $G$ be permutation groups with the same domain, and let $d\in \mathbb{N}\cup \{\infty\}$. Then the following are equivalent.
\begin{enumerate}
\item
\begin{itemize}
\item $N,G\in \fh\cap \subex_d$
\item $N\triangleleft G$, and
\item $N$ and $G$ have the same fibers.
\end{itemize}

\item\label{it:normal2} The size of fibers of $N$ and $G$ is at most $d$, and we can rename the elements of $\dom(G)=\dom(N)$ such that both $N$ and $G$ are as in one of the items (i)-(v) of Theorem~\ref{summary_fcover} and one of the items (i)-(v) of Lemma~\ref{normal_subex} holds.
\end{enumerate}
\end{theorem}

\subsection{Finite homogenizability of finite covers}

	We finish this section by showing that all finite covers of reducts of $\mathbb{Q}$ are finitely homogenizable. By Lemma~\ref{when_finhom} we can translate this fact to finite sensitivity of groups in $\fh$.
	
	We first show that for groups in $\fh$ finite sensitivity already follows from a seemingly weaker condition.
	
\begin{lemma}\label{hom_group1}
	Let $G,G^*,L,H,F$ be as in Theorem~\ref{fin_covers_q}, and let $m\geq 4|F|$. Then $G$ is $m$-sensitive if and only if the following condition holds.
	
\emph{$(*)_m$ For all $u,v\in (F\times \mathbb{Q})^n$ with $\pi_2(u)=\pi_2(v)$ if for all $p\colon [m]\rightarrow [n]$ we have $u_p\same_Gv_p$ then $u\same_Gv$.}
\end{lemma}

\begin{proof}
	The forward direction is clear. Now let us assume that $(*)_m$ holds, and let $u,v\in (F\times \mathbb{Q})^n$ be arbitrary such that $u_p\same_Gv_p$ holds for all $p\colon [m]\rightarrow [n]$. We have to show that $u\same_Gv$. We can assume without loss of generality that the tuples $u$ and $v$ are injective, i.e., their entries are pairwise distinct. Since $m\geq 2$ it also follows that $\pi_2(u_i)=\pi_2(u_j)$ if and only if $\pi_2(v_i)=\pi_2(v_j)$ for all $i,j\in [n]$. Let $U=\{u_1^*,\dots,u_{n^*}^*\}=\pi_2(\{u_1,\dots,u_n\})$ and $V=\{v_1^*,\dots,v_{n^*}^*\}=\pi_2(\{v_1,\dots,v_n\})$. Then by rearranging the tuples $u$ and $v$ we can assume without loss of generality that
\begin{align*}
u=&(u_{11},\dots,u_{1\ell_1},u_{21},\dots,u_{2\ell_2},\dots,u_{n^*\ell_{n^*}}) \text{, and}\\
v=&(v_{11},\dots,v_{1\ell_1},v_{21},\dots,v_{2\ell_2},\dots,v_{n^*\ell_{n^*}})
\end{align*}
	where $\pi_2(u_{ij})=u_i^*$ and $\pi_2(v_{ij})=v_i^*$ for all $1\leq i\leq n^*, 1\leq j\leq \ell_i$. We show that the tuples $u^*=(u_1^*,\dots,u_{n^*}^*)$ and $v^*=(v_1^*,\dots,v_{n^*}^*)$ are in the same orbit of $G^*$. Recall that $G^*$ is the automorphism of some structure $(\mathbb{Q};R)$ with $R\in \{<,\betw,\cyc,\sep,=\}$. Note that all these structures are homogeneous. Therefore, by Lemma~\ref{when_finhom} if $u^*\notsame_{G^*}v^*$ then $(u_{i_1}^*,u_{i_2}^*,u_{i_3}^*,u_{i_4}^*)\notsame_{G^*}(v_{i_1}^*,v_{i_2}^*,v_{i_3}^*,v_{i_4}^*)$. However, this is impossible since by our choice the tuples 
$$(u_{i_11},\dots,u_{i_1\ell_{i_1}},u_{i_21},\dots,u_{i_2\ell_{i_2}},u_{i_31},\dots,u_{i_3\ell_{i_3}},u_{i_41},\dots,u_{i_4\ell_{i_4}}),$$ and
$$(v_{i_11},\dots,v_{i_1\ell_{i_1}},u_{i_21},\dots,v_{i_2\ell_{i_2}},v_{i_31},\dots,v_{i_3\ell_{i_3}},v_{i_41},\dots,v_{i_4\ell_{i_4}})$$ are contained in the same $G$-orbit, and the arity of these tuples are at most $4|F|\leq m$. Now let $\alpha^*\in G^*$ such that $\alpha^*(v^*)=u^*$, and let $\alpha\in G$ with $\mu_{\pi_2}(\alpha)=\alpha^*$. Then the tuples $u$ and $\alpha(v)$ are as in $(*)_m$. Hence, by our assumption, we have $u\same_G\alpha(v)$, and thus $u\same_Gv$. 
\end{proof}

\begin{lemma}\label{hom_group2}
	Let $G,G^*,L,H,F$ be as in Theorem~\ref{fin_covers_q}. Then $G$ is $4|F||L:H|$-sensitive.
\end{lemma}

\begin{proof}
	Let $d\coloneqq|F|$ and $m\coloneqq 4d|L:H|$. We prove that $G$ is $m$-sensitive by showing that it satisfies condition $(*)_m$ from Lemma~\ref{hom_group1}.
	
	\emph{Case 1. $G=\kernel{H}{H}(\id(F)\wr \sym(\mathbb{Q}))=H\wr \sym(\mathbb{Q})$.} In this case $m=4d$. By Lemma~\ref{hom_group1} it is enough to check condition $(*)_m$. Let $u,v\in (F\times \mathbb{Q})^n$ be as in $(*)_m$. We use the notation $u^*,v^*$ and $u_{ij}$ as in the proof of Lemma~\ref{hom_group1}. By our assumption, it follows that the tuples $(\pi_1(u_{i1}),\dots,\pi_1(u_{i\ell_i}))$ and $(\pi_1(v_{i1}),\dots,\pi_1(v_{i\ell_i}))$ are in the same $H$-orbit. This implies that $u$ and $v$ are in the same $\kernel{H}{H}$-orbit, so in this case we are done.
	
	\emph{Case 2. $G=\kernel{L}{H}(\id(F)\wr \sym(\mathbb{Q}))$.} Let $N:=H\wr \sym(\mathbb{Q})$. Then $N\triangleleft G$ and $|G:N|=|L:H|$ by Lemma~\ref{factors}, and we already know by Case 1 that $N$ is $4d$-sensitive. Therefore, by Lemma~\ref{fin_index} we obtain that $G$ is $4d|L:H|$-sensitive.

	Finally, for the general case we check the condition $(*)_m$ again. Let $G_0:=\kernel{L}{H}(\id(F)\wr \sym(\mathbb{Q}))$. Let $u,v$ be as in the assumption of $(*)_m$. We claim that for all $p\colon [m]\rightarrow [n]$ we have $u_p\same_{G_0}v_p$. Indeed, by our assumption, we already know that $u_p\same_Gv_p$, and thus in fact $u_p\same_{G_{\pi_2(u_p)}}v_p$. Since $m>2d$, it follows that $|\pi_2(u_p)|\geq 3$. Thus, by Lemma~\ref{fix_some} it follows that $G_{\pi_2(u_p)}=G_{\pi_2(v_p)}=(\kernel{L}{H})_{\pi_2(u_p)}$ which in particular implies that $u_p\same_{G_0}v_p$. By Claim 2 we already know that $G_0$ is $m$-sensitive. Therefore, $u\same_{G_0}v$, in particular $u\same_Gv$, and this is what we needed to show.
\end{proof}

\begin{corollary}\label{ff_fincom}
	Every finite cover of a reduct of $\mathbb{Q}$ is finitely homogenizable.
\end{corollary}

\begin{proof}
	Theorem~\ref{fin_covers_q} and Lemma~\ref{hom_group2} together imply that every finite cover of $\aut(\mathbb{Q})$ is finitely sensitive. Then the statement of the corollary follows from Lemma~\ref{when_finhom2}.
\end{proof}

	Let $\fa$ be a homogeneous structure in $\fhs$. By Theorem~\ref{fin_covers_q} we know that in this case $\fa$ is \emph{multichainable} in the sense of~\cite{pouzet2007profile}. Then by Theorem 39 in~\cite{pouzet2007profile} we know that $\fa$ is in fact finitely bounded. Combining this fact with Corollary~\ref{ff_fincom} we obtain the following.
	
\begin{lemma}\label{ff_finbound}
	$\fhs \subset \fbh$.
\end{lemma}

%% file: fincover_h.tex

\section{Finite covers of hereditarily cellular groups}\label{sect:fin_h}

	In this section we describe all pairs $(G,G^*)$ where $G^*$ is hereditarily cellular and $G$ is a finite cover of $G^*$. We first show that in case $G$ is also hereditarily cellular.

\begin{lemma}\label{ms_covers}
	The classes $\gone_n$ are closed under taking closed finite covers.
\end{lemma}

\begin{proof}
	We show the statement by induction on $n$. For $n=0$ the statement is trivial. For the induction step let $G^*\in \gone_n$, and let $(G,G^*,\pi)$ be a finite cover. In this proof we use the recursive description of the classes $\gone_n$ given in Theorem~\ref{omega_part}.

	So let $\partition^*:=(K^*,\bigeq^*,\smalleq^*)$ be an $\omega$-partition of $G^*$ with components contained in $\gone_{n-1}$. Let $$K=\pi^{-1}(K^*), \bigeq:=\pi^{-1}(\bigeq^*), \smalleq:=\pi^{-1}(\smalleq^*),\partition:=(K,\bigeq,\smalleq).$$

	Then it is clear from the definition that $\partition$ is an $\omega$-partition of $G$. It remains to prove that the components of $\partition$ are contained in $\gone_{n-1}$. We put $X\coloneqq \dom(G)$ and $X^*\coloneqq \dom(G^*)$. Let $Y$ be a $\smalleq$-class, and let $Y^*:=\pi(Y)$. Let $$N^*:=(G^*)_{((Y^*))},N:=\mu_{\pi}^{-1}(N^*).$$ Then $N$ is a finite cover of $N^*$ witnessed by the map $\pi|_Y$. Since $N^*$ is a component of $\partition^*$ we know that $N^*\in \gone_{n-1}$. Thus, by the induction hypothesis $N\in \gone_{n-1}$.

	By an abuse of notation we write $N_F$ for the group $$\{g|_Y\colon g\in G_F\wedge \mu_{\phi}(g)\in N^*\}$$ where $F$ is any subset of $X$. We claim that for every finite set $F\subset X\setminus Y$ the group $N_F$ is a finite-index closed subgroup of $N$. The fact that $N_F$ is closed follows from Lemma~\ref{fin_closed} (see the proof of Lemma~\ref{closed_factor}). Now let $a_1,\dots,a_k$ be an enumeration of $F$. Then the orbit of $a_i$ with respect to the group $\mu_{\pi}^{-1}((G^*)_{X^*\setminus Y^*})$ is contained in the $\sim_{\pi}$ class of $a_i$, in particular it is finite. Therefore, $$|(\mu_{\pi}^{-1}((G^*)_{X^*\setminus Y^*}))_{a_1,\dots,a_i}\colon (\mu_{\pi}^{-1}((G^*)_{X^*\setminus Y^*}))_{a_1,\dots,a_{i+1}}|<\infty$$ for all $i=1,\dots,k$. This in turn implies that $$|(\mu_{\pi}^{-1}((G^*)_{X^*\setminus Y^*})):(\mu_{\pi}^{-1}((G^*)_{X^*\setminus Y^*}))_F|<\infty,$$ and thus $|N:N_F|$ is finite. 
	
	We obtained that for every finite subset $F$ of $X\setminus Y$ we have $N_F\supset N_{\fin}$. On the other hand, by definition we have
$$G_{((Y))}=\bigcap\{N_F\colon F\subset X\setminus Y, |F|<\infty\}.$$ Therefore, $G_{((Y))}\supset N_{\fin}$. Then by Lemma~\ref{fin_index_gone} we can conclude that $G_{((Y))}\in \gone_{n-1}$, and this is what we needed to show.
\end{proof}

	Next we give a recursive description of all closed finite covers of hereditarily cellular groups up to isomorphism.
	
\begin{definition}\label{def:fcoverh}
	We denote by $\mathbf{A}$ the minimal class of finite covers which contains all finite covers of finite degree groups and which satisfies the following condition.

	Let us assume that $(H,H^*,\rho)\in \mathbf{A}$, let $N_0,\dots,N_k$ groups such that $|\dom(N_0)|< \infty, N_0=\id(\dom(N_0)), N_i\in \gone$, and $N_i\coloneqq \prod_{i=0}^kN_i\triangleleft H$. Let $N_i^*\coloneqq \mu_{\pi}(N_i)$. Then $(G,G^*,\pi)\in \mathbf{A}$ where
\begin{itemize}
\item $G=\biggroup(H;N_0,\dots,N_k)$ (see Definition~\ref{def:bigg}),
\item $G^*=\biggroup(H^*;N_0^*,\dots,N_k^*)$,
\item $\pi((a,n))=(\rho(a),n)$.
\end{itemize}
\end{definition}

\begin{lemma}\label{ms_covers_class}
	Let $(G,G^*,\pi)$ a finite cover. Then the following are equivalent.
\begin{enumerate}
\item\label{it:fg1} $G\in \gone$.
\item\label{it:fg2} $G^*\in \gone$ and $G$ is closed.
\item\label{it:fg3} $(G,G^*,\pi)$ is isomorphic to some element in $\mathbf{A}$.
\end{enumerate}
\end{lemma}

\begin{proof}
	The implication~\ref{it:fg3}$\rightarrow$\ref{it:fg2} is clear from definition, and the implication~\ref{it:fg2}$\rightarrow$\ref{it:fg1} follows from Lemma~\ref{ms_covers}.
	
	It remains to show that item~\ref{it:fg1} implies item~\ref{it:fg3}. Let us assume that $(G,G^*,\pi)$ is a finite cover, and $G\in \gone_n$. We show by induction on $\rk(G)$ that $(G,G^*,\pi)$ is as in item~\ref{it:fg3}. For $\rk(G)=0$ there is nothing to prove. Otherwise, we can rename the elements of $\dom(G)$ such that $G=\biggroup(H;N_0,\dots,N_k)$ for some $H,N_0,\dots,N_k\in \gone_{n-1}$. Then since $G\supset \id(\dom(N_0))\times \prod_{i=1}^k(\id(\dom(N_i))\wr \sym(\omega))$ it follows that for all $(a,m),(b,n)\in \dom(H)\setminus \dom(N_0)$ with $n\neq m$ we have $(a,m)\not\sim_{\pi}(b,n)$. Indeed, if $(a,m)\sim_{\pi}(b,n)$ would hold for some $n\neq m$ then we would have $(a,m')\sim_{\pi}(b,n)$ for any $m'\neq m,n$, and thus the $\sim_{\pi}$-class of $(a,m)$ would contain infinitely many elements. Using the fact that $\sim_{\pi}$ is a congruence of $G$ it also follows that $(a,0)\sim_{\pi}(b,0)$ if and only if $(a,n)\sim_{\pi}(b,n)$ for all $n\in \omega$. Let us define $\sim'$ on $\dom(H)$ such that $a\sim'b$ if and only if $(a,0)\sim_{\pi}(b,0)$. Then $\sim'$ has finite classes, that is $H$ is a finite cover of $H^*:=H/\unsim$. Let $\rho\colon \dom(H)\rightarrow \dom(H^*), a\mapsto [a]_{\sim'}$. Then by the induction hypothesis we can rename the elements of $\dom(H)$ and $\dom(H^*)$ such that $(H,H^*,\rho)\in \mathbf{A}$. We obtain in particular that $H^*\in \gone$. Let $N_i^*:=N_i/\unsim'|_{\dom(N_i)}$. Then similarly we have $N_i^*\in \gone$. It is also clear from the definition that $|\dom(N_0)|< \infty, N_0=\id(\dom(N_0))$, and $N^*\coloneqq\prod_{i=0}^kN_i^*$ is a normal subgroup of $H^*$. In this case the map $\overline{\rho}\coloneqq (a,n)\mapsto (\rho(a),n)$ defines a finite covering map from $G$ to $\biggroup(H^*;N_0^*,\dots,N_k^*)$. It is also clear from the definition that $\pi$ and $\overline{\rho}$ have the same fibers. This implies that $G^*$ is isomorphic to $\biggroup(H^*;N_0^*,\dots,N_k^*)$, and thus the cover $(G,G^*,\pi)$ is isomorphic to some cover in $\mathbf{A}$.
\end{proof}

    Lemma~\ref{ms_covers_class} has the following consequence which we will need later.

\begin{lemma}\label{ms_covers_trivial}
	Let $(G,G^*,\pi)$ a finite cover with $G\in \gone$. Then there exist $G_0\leq G$ and $G_0^*\leq G^*$ with $G_0\in \gone$ such that $(G_0,G_0^*,\pi)$ is strongly trivial.
\end{lemma}

\begin{proof}
    By Lemma~\ref{ms_covers_class} we can assume that $(G,G^*,\pi)\in \mathbf{A}$. We show the statement by induction of $\rk(G)$. If $\rk(G)=0$, that is $\dom(G)$ is finite then we can take both $G_0$ and $G_0^*$ to be the identity group. Otherwise, let us assume that $H,H^*,\rho,N_i,N_i^*: i\in \{1,\dots,k\}$ are as in Definition~\ref{def:fcoverh}. We can assume without loss of generality that $H=N=\prod_{i=0}^kN_i$ and $H^*=N^*=\prod_{i=0}^kN_i$. By the induction hypothesis, we know that there exist some $H_0\leq H, H_0\in \gone$ and $H_0^*\leq H^*$ such that $(H_0,H_0^*,\rho)$ is a strongly trivial cover. Let 
\begin{align*}
G_0\coloneqq& \biggroup(H_0;(H_0)|_{\dom(N_0)},\dots,(H_0)|_{\dom(N_k)}),\\
G_0^*\coloneqq& \biggroup(H_0^*;(H_0^*)|_{\dom(N_0^*)},\dots,(H_0^*)|_{\dom(N_k^*)}).
\end{align*}
    Then it is easy to check that the groups $G_0$ and $G_0^*$ are as stated in the lemma.
\end{proof}

%% file: finish.tex

\section{Finishing the classification}\label{sect:finish}

	In this section we put together the pieces from the previous sections and present a classification of groups in all the classes $\subex_d\colon d\in \mathbb{N}\cup \{\infty\}$. We start with a necessary condition which follows directly from the results of the previous sections.
	
\begin{lemma}\label{class1}
	Let $d\in \mathbb{N}\cup \{\infty\}$. Then every permutation group in $\subex_d$ is isomorphic to some group $G=\lift(\tilde{G},G^*,\data)$ where
\begin{enumerate}
\item the triple $(\tilde{G},G^*,\data)$ is as in Definition~\ref{def:recover},
\item $(\tilde{G},G^*,\pi_2)$ is isomorphic to some finite cover in $\mathbf{A}$ (see Definition~\ref{def:fcoverh}),
\item for all $a\in X^*$ we have $\data(a)=(B_a,F_a,\phi_a)$ where either $|\dom(F_a)|=1$ or
\begin{enumerate}[(i)]
\item $B_a,F_a\in \fh\cap \subex_d$
\item $B_a$ and $F_a$ have the same fibers, and
\item $B_a\triangleleft F_a$.
\end{enumerate}
\end{enumerate}
\end{lemma}

\begin{proof}
	The statement of the lemma essentially follows from Theorem~\ref{simon_final} and the discussion in Section~\ref{sect:fff}. The only thing we have to check is that if $G=\lift(\tilde{G},G^*,\data)$ and $G\in \subex_d$ then for the groups $B_a$ and $F_a$ as above we have $B_a,F_a\in \subex$. Since $|F_a:B_a|<\infty$ it is enough to check that $B_a\in \subex$. The finiteness of $|F_a:B_a|$ also implies that there exists some finite subset $S$ of $\dom(G)\setminus \pi_2^{-1}(a)$ such that $G_S|_{\pi_2^{-1}(a)}$ is isomorphic to $B_a$ (see for instance the proof of Lemma~\ref{tricky_restriction}). On the other hand, clearly $\os_n(G_S|_{\pi_2^{-1}(a)})$ is at most $\os_n(G)$. Therefore, $B_a\in \subex$.
\end{proof}

\begin{remark}\label{at_most5}
	Note that in the construction of the group $\tilde{G}$ (see the proof of Lemma~\ref{recover2}) the fibers $F$ of $\tilde{G}$ can be chosen arbitrarily as long as $F_a/B_a$ embeds into $\sym(F)$. On the other hand, we know by Corollary~\ref{factors_options} that all these quotient groups are isomorphic to $Z_k$ for $k\in \{1,2,3,4,6\}$ or $Z_2\times Z_2$. Since all these groups embed into $S_5$, this means that we can always choose the group $\tilde{G}$ above so that all of its fibers are of size at most 5.
\end{remark}

	Next we show that the converse of Lemma~\ref{class1} is also true: every group $G=\lift(\tilde{G},G^*,\data)$ as described above is contained in the class $\subex_d$. We will show this via a third equivalent description of the classes $\subex_d$ analogous to Theorem~\ref{thm:build_groups}. Namely, we show that $\subex_d$ is exactly the smallest class of permutation groups which contains $(\subex_d\cap \fh)\cup \id(\{\emptyset\})$ and is closed under taking finite direct products, wreath products with $\sym(\omega)$, finite-index supergroups, and isomorphisms. We present a slightly stronger formulation of this statement in Theorem~\ref{final} below.
	
	We start by showing that the classes $\subex_d$ are in fact closed under the aforementioned constructions. By definition, this is clear for taking supergroups, and isomorphisms.
	
\begin{lemma}\label{upper_bound}
	The classes $\subex_d\colon d\in \mathbb{N}$ are closed under taking finite direct products and wreath product with $\sym(\omega)$.
\end{lemma}

\begin{proof}
	We first show that $\subex_d$ is closed under taking finite direct products. For this, it is enough to show that if $G_1,G_2\in \subex_n$ then $G_1\times G_2\in \subex_d$. Let $X_i\coloneqq \dom(G_i), X\coloneqq X_1\cup X_2$. Let $A\in {X\choose n}$. Then the orbit of $A$ is uniquely determined by $|A\cap X_1|$ and the $G_i$-orbit of $A\cap X_i$ for $i=1,2$. For $|A\cap X_1|$ we have $n+1$ choices, and if $k\coloneqq |A\cap X_1|$ is given then we have $\os_k(G_1)$ choices for the $G_1$-orbit of $A\cap X_1$ and $\os_{n-k}(G_2)$ choices for $A\cap X_2$. Therefore, $\os_n(G)=\sum_{k=0}^n\os_k(G_1)\os_{n-k}(G_2)$. Let $\varepsilon>0$. Then by our assumption we have $\os_n(G_i)<c(\gamma_d+\varepsilon)^n$ for some $c\in \mathbb{R}$ (which depends on $\varepsilon$). By the formula above we obtain 
$$\os_n(G)=\sum_{k=0}^n\os_k(G_1)\os_{n-k}(G_2)\leq (n+1)c^2(\gamma_d+\varepsilon)^n$$
	which is clearly smaller than $(\gamma_d+2\varepsilon)^n$ if $n$ is large enough.

	For the second part of the lemma let us assume that $G\in \subex_d$. We show that $G\wr \sym(\omega)\in \subex_d$. Let $X=\dom(G)$, and let $A\in {{X\times \omega}\choose n}$. We can assume without changing the orbit of $A$ that $$|A\cap (X\times \{0\})|\geq |A\cap (X\times \{1\})|\geq \dots \geq|A\cap (X\times \{n\})|=0.$$ We put $k_i\coloneqq |A\cap (X\times \{i\})|$. Then $\sum_{i=0}^{n-1}k_i=n$. Observe that the $(G\wr \sym(\omega))$-orbit of $A$ is uniquely determined by the sequence $(k_0,\dots,k_{n-1})$ and the $G$-orbits of the sets $A_i:=\{a\colon (a,i)\in A\}$. The number of choices for the sequence $(k_0,\dots,k_{n-1})$ is $p(n)$, that is the number of partitions of $n$. The asymptotic behavior of the partition function $p(n)$ is well-known: $p(n)\sim \frac{1}{4n\sqrt{3}}\exp\bigl(\pi\sqrt{\frac{2n}{3}}\bigr)$ (see for instance~\cite{natanson2003elementary}, Section 16.4). In particular, for all $c>0$ we have $p(n)<c^n$ if $n$ is large enough. By our assumption, we know that $\os_n(G)\leq c(\gamma_d+\varepsilon)^n$ for some $c\in \mathbb{R}$. Then based on our argument above we can give the following upper bound on $\os_n(G\wr \sym(\omega))$:
\begin{align*}
\os_n(G\wr \sym(\omega))=&\sum\Bigl\{\prod_{i=0}^{n-1}\os_{k_i}(G): k_0\geq \dots\geq k_{n-1}\wedge \sum_{i=0}^{n-1}k_i=n\Bigr\}\\
\leq& \sum\Bigl\{\prod_{i=0}^{n-1}c(\gamma_d+\varepsilon)^{k_i}: k_0\geq \dots\geq k_{n-1}\wedge \sum_{i=0}^{n-1}k_i=n\Bigr\}\\
=&c^np(n)(\gamma_d+\varepsilon)^n
\end{align*}
	which is smaller than $(\gamma_d+2\varepsilon)^n$ if $n$ is large enough.
\end{proof}

\begin{remark}
	Note that for $d=\infty$ our proof above does not work. However, as we will see, $\subex=\bigcup_{d\in \mathbb{N}}\subex_d$, so the statement of the lemma is still true in this case.
\end{remark}

	Now let us assume that $(\tilde{G},G^*,\pi)$ where $\tilde{G}\in \gone$. Then by Lemma~\ref{ms_covers_class} we can find a finite cover $(\tilde{H},H^*,\rho)\in \mathbf{A}$ and groups $\tilde{N_0},\dots,\tilde{N_k},N_0^*,\dots,N_k^*$ such that after renaming the elements of the domains of $\tilde{G}$ and $G^*$ we have
	
\begin{itemize}
\item $\tilde{H},H^*\in \gone$,
\item $|\dom(\tilde{N_0})|< \infty$,
\item $\tilde{N}_i\in \gone$,
\item $N_i^*\coloneqq \mu_{\pi}(\tilde{N_i})$,
\item $\tilde{G}=\biggroup(\tilde{H};\tilde{N_0},\dots,\tilde{N_k})$,
\item $G^*=\biggroup(H^*;N_0^*,\dots,N_k^*)$,
\item $\pi((a,n))=(\rho(a),n)$.
\end{itemize}

	Now let $G=\lift(\tilde{G},G^*,\data)$ be as in the conclusion of Lemma~\ref{class1}. Note that in this case $\data$ is constant on elements of the form $(a,n)\colon n\in \omega$. Moreover, the map $a\mapsto (a,n)$ gives rise to a bijection between the fiber groups of $H^*$ and $G^*$. Thus, for $a\in Y^*\coloneqq \dom(H^*)$ we can define $B_a\coloneqq B_{(a,\cdot)},F_a\coloneqq F_{(a,\cdot)}$ and $\phi_a\colon \gamma\mapsto \pi_1(\phi_{(a,0)}(\bar{\gamma}))$ where $\bar{\gamma}\colon (a,n)\mapsto (\gamma(a),n)$. Finally, we put $\data(a)\coloneqq (B_a,F_a,\phi_a)$. Then it is easy to check that the triples $(\tilde{H},H^*,\data)$ and $(\tilde{N_i},N_i^*,\data)$ also satisfy the conditions in Definition~\ref{def:recover}. Let $H\coloneqq \lift(\tilde{H},H^*,\data)$, $N_i\coloneqq \lift(\tilde{N_i},N_i^*,\data)$. As in Subsection~\ref{sect:her} we write 
\begin{itemize}
\item $N\coloneqq \prod_{i=0}^kN_i$
\item $X\coloneqq \dom(G)$,
\item $Y_i\coloneqq \dom(N_i)$,
\item $Y\coloneqq \dom(N)=\bigcup_{i=0}^kY_i,$
\end{itemize}
	and we also use the analogous notation with tildes and asterisks.
	
\begin{lemma}\label{generate_all}
	With the above notation $N\triangleleft H$ and the group $G$ is isomorphic to $\biggroup(H;N_0,\dots,N_k)$.
\end{lemma}

\begin{proof}
	Let $n\in N$ and $h\in H$, and let $\tilde{n}\in \tilde{N}$ and $\tilde{h}\in \tilde{H}$ witnessing this fact as in Definition~\ref{def:recover}. Then since $\tilde{N}\triangleleft \tilde{H}$ we have $\tilde{h}^{-1}\tilde{n}\tilde{h}\in \tilde{N}$ which witnesses the fact that $h^{-1}nh\in N$. Therefore, $N\triangleleft H$.

	The isomorphism from right to left is induced by the bijection $\iota\colon ((u,a),n)\mapsto (u,(a,n))$. Let $G'\coloneqq \biggroup(H;N_0,\dots,N_k)$.
	
	We introduce the following notation:

\begin{itemize}
\item $T\coloneqq \id(Y_0\times \{0\})\times \prod_{i=1}^k(\id(Y_i)\wr \sym(\omega))$
\item $N_{n,i}\coloneqq\{\kappa(\sigma,n): \sigma\in N_i\}$
\item $K\coloneqq\{\kappa(\sigma,\omega): \sigma\in H\}$,
\end{itemize}

	and we also use the analogous notation with tildes and asterisks. Then we know that $G'$ is the closure of the group generated by $T,N_{n,i}$ and $K$, and $\tilde{G}$ is the closure of the group generated by $\tilde{T},\tilde{N_{n,i}}$ and $\tilde{K}$.
	
	Let us fix for all $a\in Y$ a map $\psi_a$ which is a right inverse of $\phi_a$ such that
\begin{itemize}
\item $\phi_a(\id)=\id$, and
\item $\phi_a$ only depends on the orbit of $a$.
\end{itemize}

	Then we define a map $\gamma$ by the formula $\gamma(g)\coloneqq (u,a)\mapsto (\phi_a(\eta_a)(u),\mu_{\pi_2}(a))$. This definition works on all groups $\tilde{H},\tilde{N_i},\tilde{G}$, and for all $g$ in any of these groups we have $\gamma(\tilde{g})=g$.
	
	Let $B$ be the kernel of the homomorphism $G\rightarrow \tilde{G}, g\mapsto \tilde{g}$. Then $B$ is equal to the direct product of the groups $\{(u,(a,n))\mapsto (\beta(u),(a,n)): \beta\in B_{(a,n)}\sigma\}$. Since $B_a=B_{(a,n)}$ this implies immediately by the definition of $\iota$ that $B\subset \mu_{\iota}(G')$.
	
	Using the description of the group $B$ above one can check easily that

\begin{itemize}
\item $\mu_{\iota}(T)=\gamma(\tilde{T})$
\item $\gamma(N_{n,i})\subset \mu_{\iota}(\widetilde{N_{n,i}})\subset \gamma(N_{n,i})B$,
\item $\gamma(K)\subset \mu_{\iota}(\tilde{K})\subset \gamma(K)B$.
\end{itemize}

	Therefore, we have
\begin{align*}
\mu_{\iota}(G')=&\overline{\langle \mu_{\iota}(T),\mu_{\iota}(N_{n,i}): (n,i)\in \omega\times \{0,\dots,k\},\mu_{\iota}(K)\rangle}\\
=&\overline{\langle \gamma(\tilde{T}),\gamma(\widetilde{N_{n,i}}): (n,i)\in \omega\times \{0,\dots,k\},\gamma\tilde{K}\rangle}B=\gamma(\tilde{G})B=G.\qedhere
\end{align*}
\end{proof}

\begin{lemma}\label{general_finite_index}
	With the above notation $|\biggroup(H;N_0,\dots,N_k):\biggroup(N;N_0,\dots,N_k)|$ is finite.
\end{lemma}

\begin{proof}
	Let $G=\biggroup(H;N_0,\dots,N_k)$ and $G_0=\biggroup(N;N_0,\dots,N_k)$. Then the map homomorphism $g\mapsto \tilde{g}$ defined on $G$ and $G_0$ have the same kernel and they are equal where both of them are defined. Therefore, $|G:G_0|=|\tilde{G}:\tilde{G_0}|$ which is finite by our argument in Lemma~\ref{fin_index_gone}.
\end{proof}

\begin{lemma}\label{subex_build}
	Let $d\in \mathbb{N}\cup \{\infty\}$, and let $G\in \subex_d$. Then $G$ can be constructed from the following groups by taking isomorphism, finite-index supergroups, finite direct products and wreath products with $\sym(\omega)$.
\begin{enumerate}
\item $\id(\{\emptyset\})$.
\item\label{it:subex_b2} $G=H\wr \aut(\mathbb{Q};<)$ where $H$ is a finite-degree highly set-transitive group with degree at most $d$.
\item\label{it:subex_b3} $G=\langle H\wr \aut(\mathbb{Q};<), \turn_{\sigma}\rangle$ where $H$ is a highly set-transitive group with degree at most $d$ and $H\triangleleft H\sigma$.
\end{enumerate}
\end{lemma}

\begin{remark}
	Clearly, the ``finite-degree'' condition is only needed if $d=\infty$.
\end{remark}

\begin{proof}
	We can assume that $G=\lift(\tilde{G},G^*,\data)$ as in Lemma~\ref{class1}. We prove the statement by induction on $\rk(G^*)$. 
	
	We first show the induction step. By defining the groups $H,H^*,N_0,\dots,N_k,$ as above we know by Lemma~\ref{generate_all} that $G$ is isomorphic to $\biggroup(H;N_0,\dots,N_k)$ where $H,N_0,\dots,N_k\in \subex$, and $\rk(H^*)<\rk(G^*)$. By Lemma~\ref{general_finite_index} and Corollary~\ref{subex_finex} we can conclude that $\biggroup(N;N_0,\dots,N_k)\in \subex_d$. Clearly, $\os_n(N_i),\os_n(H)\leq \os_n(\biggroup(N;N_0,\dots,N_k))$. Therefore, all the groups $N_0,\dots,N_k,N,H$ are contained in $\subex_d$. Then the induction step follows from Lemma~\ref{general_finite_index} and the definition of the group $\biggroup(H;N_0,\dots,N_k)$.
	
	By our argument above it is enough to show the lemma in the case when $\tilde{G}$ and $G^*$ have finite domains. Let $B$ be the kernel of the natural homomorphism $G\rightarrow \tilde{G}$. Then $|G:B|$ is finite and $B$ is isomorphic to the direct product of the groups $B_a: a\in \dom(G^*)$. By our construction all groups in $B_a$ are either isomorphic to $\id(\{\emptyset\})$ or it is contained in $\fh\cap \subex$. In the latter case it follows from Theorem~\ref{summary_fcover} and Lemma~\ref{normal_subex} that $B_a$ has a finite-index normal subgroup as in items~\ref{it:subex_b2} or~\ref{it:subex_b3}. It remains to show that the fiber size of $B_a$ is at most $d$. By using Corollary~\ref{subex_finex} we know that $B_a\in \subex_d$, and clearly $\os_n(B_a)\leq \os_n(B)$ for all $a\in \dom(G^*)$. Thus, we are done by Lemma~\ref{when_subex}.
\end{proof}

\begin{lemma}\label{subex_reduct}
	Let $d\in \mathbb{N}\cup \{\infty\}$, and let $G\in \subex_d$. Then $G$ is isomorphic to a closed supergroup of a finite direct product of groups of the form $\Gamma\wr N$ where $N\in \gone$ and $\Gamma=\id(\{\emptyset\})$ or $\Gamma=H\wr \aut(\mathbb{Q};<)$ where $H$ is a finite-degree highly set-transitive group with degree at most $d$.
\end{lemma}

\begin{proof}
    We can assume that $G=\lift(\tilde{G},G^*,\data)$ as in Lemma~\ref{class1}. Let $\tilde{G_0}\leq \tilde{G}$ and $G_0^*\leq G^*$ as in the conclusion of Lemma~\ref{ms_covers_trivial}, and let $G_0$ be the set of those permutations in $G$ which are witnessed by some permutation $\tilde{g}\in \tilde{G}_0$ as Definition~\ref{def:recover}. Then $(G_0,G_0^*,\pi_2)$ is still a cover and since $(\tilde{G}_0,G_0^*,\pi_2)$ is strongly trivial it follows that for all $a\in \dom(G^*)$ both the pointwise binding group and the fiber group of $(G_0,G_0^*,\pi_2)$ at $a$ is $B_a$. This implies that $G_0$ can be written as a direct product of the wreath products $B_a\wr (G_0^*)|_O$ where $O$ goes through the orbits of $G_0^*$ and $a\in O$. Since $G_0^*$ is hereditarily cellular, so are the groups $(G_0^*)|_O$. Finally, by the same argument as in the proof of Lemma~\ref{subex_build} we can conclude that each $B_a$ is a closed supergroup of some $\Gamma_a$ as described in the lemma.
\end{proof}

	We summarize the results of this section in the following theorem.
	
\begin{theorem}\label{final}
	Let $d\in \mathbb{N}\cup \{\infty\}$, and let $G$ be a permutation group. Then the following are equivalent.
\begin{enumerate}
\item\label{it:final1} $G\in \subex_d$.
\item\label{it:final2} $G$ is isomorphic to some group $\lift(\tilde{G},G^*,\data)$ where
\begin{itemize}
\item the triple $(\tilde{G},G^*,\data)$ is as in Definition~\ref{def:recover},
\item $(\tilde{G},G^*,\pi_2)$ is isomorphic to some finite cover in $\mathbf{A}$,
\item for all $a\in X^*$ we have $\data(a)=(F_a,B_a,\phi_a)$ where either $|\dom(F_a)|=1$ or
\begin{enumerate}[(i)]
\item $B_a,F_a\in \fh\cap \subex_d$
\item $B_a$ and $F_a$ have the same fibers, and
\item $B_a\triangleleft F_a$.
\end{enumerate}
\end{itemize}
\item\label{it:final3} $G$ can be constructed from the groups
\begin{enumerate}\label{it:build_up}
\item $\id(\{\emptyset\})$,
\item $G=H\wr \aut(\mathbb{Q};<)$ where $H$ is a finite-degree highly set-transitive group with degree at most $d$,
\item $G=\langle H\wr \aut(\mathbb{Q};<), \turn_{\sigma}\rangle$ where $H$ is a highly set-transitive group with degree at most $d$ and $H\triangleleft H\sigma$
\end{enumerate}
by taking isomorphism, finite-index supergroups, finite direct products and wreath products with $\sym(\omega)$.
\item\label{it:final4} $G$ can be constructed from the groups
\begin{enumerate}\label{it:build_up2}
\item $\id(\{\emptyset\})$,
\item $G=H\wr \aut(\mathbb{Q};<)$ where $H$ is a finite-degree highly set-transitive group with degree at most $d$
\end{enumerate}
by taking isomorphism, closed supergroups, finite direct products and wreath products with $\sym(\omega)$.
\item\label{it:final5} $G$ is isomorphic to a closed supergroup of a finite direct product of groups of the form $\Gamma\wr N$ where $N\in \gone$ and $\Gamma=\id(\{\emptyset\})$ or $\Gamma=H\wr \aut(\mathbb{Q};<)$ where $H$ is a finite-degree highly set-transitive group with degree at most $d$.
\end{enumerate}
\end{theorem}

	Note that all possible pairs $(B_a,F_a)$ as in items (i)-(iii) are classified in Theorem~\ref{summary_fcover} and Lemma~\ref{normal_subex} up to renaming the elements of $\dom(F_a)$. This finishes the classification of all classes $\subex_d: d\in \mathbb{N}\cup \{\infty\}$.

\begin{proof}
	We prove the theorem by showing the implications~\ref{it:final1}$\rightarrow$\ref{it:final2}$\rightarrow$\ref{it:final3}$\rightarrow$\ref{it:final4}$\rightarrow$\ref{it:final1} and~\ref{it:final2}$\rightarrow$\ref{it:final5}$\rightarrow$\ref{it:final4}

	The implication~\ref{it:final1}$\rightarrow$\ref{it:final2} is Lemma~\ref{class1}. The implications~\ref{it:final2}$\rightarrow$\ref{it:final3} and~\ref{it:final2}$\rightarrow$\ref{it:final5} are shown in Lemmas~\ref{subex_build} and~\ref{subex_reduct}, respectively. The implications~\ref{it:final3}$\rightarrow$\ref{it:final4} and~\ref{it:final5}$\rightarrow$\ref{it:final4} are obvious. 

	Finally, we show the implication~\ref{it:final4}$\rightarrow$\ref{it:final1}. Let us first assume that $d\in \mathbb{N}$. By Lemma~\ref{upper_bound} we know that the class $\subex_d$ class is closed under taking finite direct products and wreath products with $\sym(\omega)$. Also, clearly $\subex_d$ is closed under isomorphisms and taking closed supergroups. Finally, by Lemma~\ref{when_subex} we know that the groups listed in item~\ref{it:final4} are all contained in $\subex_d$. Now for $d=\infty$ if $G$ as in~\ref{it:final4} then since the construction of $G$ is finite, it follows that it also satisfies item~\ref{it:final4} for some $d'\in \mathbb{N}$. Then by the argument above we can conclude that $G\in \subex_{d'}\subset \subex$.
\end{proof}

\begin{corollary}\label{union}
	$\subex=\bigcup_{d\in \omega}\subex_d$, that is for all $G\in \subex$ we have $\os_n(G)\leq c^n$ for some $c<2$ if $n$ is large enough.
\end{corollary}

\begin{proof}
	Clear from the equivalence of items~\ref{it:final1} and~\ref{it:final3} (or~\ref{it:final4}) of Theorem~\ref{final}.
\end{proof}

	The following is a straightforward consequence of Theorem~\ref{final}.
	
\begin{theorem}\label{final_struct}
	Let $d\in \mathbb{N}\cup \{\infty\}$, and let $\fa\in \subexs_d$. Then the following hold.
\begin{enumerate}
\item\label{it:struct1} $\fa$ can be constructed from the one-element pure set and finite covers of $(\mathbb{Q};<)$ and $(\mathbb{Q};\cyc)$ of fiber size at most $d$ by taking isomorphism, finite index reducts, finite disjoint unions and wreath products with the pure set.
\item\label{it:struct2} $\fa$ is bidefinable with a reduct of a finite disjoint union of wreath products of the from $\ff\wr \mathfrak{H}$ where each $\mathfrak{H}$ is hereditarily cellular and each $\ff$ is either the one-element pure set or a finite cover of $(\mathbb{Q};<)$ of fiber size at most $d$.
\end{enumerate}
    Moreover, if $d=1$ then either of the items above implies $\fa\in \subexs_1$.
\end{theorem}

\begin{remark}
    Note that in general items~\ref{it:struct1} and~\ref{it:struct2} above do not imply $\fa\in \subexs_d$ since not all finite covers of $(\mathbb{Q};<)$ are in $\subexs$.
\end{remark}

    Note that the equivalence of $\fa\in \subexs_1$ and item~\ref{it:struct2} of Theorem~\ref{final_struct} in the $d=1$ case confirms Conjecture 1 in~\cite{braunfeld2022monadic}. Indeed, in this case we can choose $<$ to be the union of the orders $\{((s,a),(t,a)): s<t\}: a\in H$.

\begin{corollary}\label{fin_hom}
	$\subexs\subset \fbh$, i.e., every structure in $\subexs$ is interdefinable with a finitely bounded homogeneous structure.
\end{corollary}

\begin{proof}
	Follows from Lemma~\ref{fin_index_fbh}, Lemma~\ref{ff_finbound}, and item~\ref{it:struct1} of Theorem~\ref{final_struct}.
\end{proof}

\section{Orbit growth analysis}\label{sect:orbit}

	By taking a closer look at our calculations on the unlabelled profile in this section we obtain the following.
	
\begin{theorem}\label{orbit_growth}
	Let $G$ be a permutation group acting on a countable set. Then one of the following holds.
\begin{enumerate}
\item\label{it:slow} $\os_n(G)<c^n$ for all $c>1$ if $n$ is large enough, and $\os_n(G)=\os_n(G_*)$ for some $G_*\in \gone$.
\item\label{it:medium} There exists some $d\geq 2$ and $c\in \mathbb{R}^+$ such that $c\gamma_d^n<\os_n(G)$ and for all $\varepsilon>0$ we have $\os_n(G)<(\gamma_d+\varepsilon)^n$ if $n$ is large enough.
\item\label{it:fast} $\os_n(G)\geq \frac{2^n}{p(n)}$ for some polynomial $p$.
\end{enumerate}
\end{theorem}

\begin{proof}
	We can assume without loss of generality that $G$ is closed. If $G\not\in \subex$ then by definition, item~\ref{it:fast} holds, so let us assume that $G\in \subex$. Let $\tilde{G},G^*,\data$ be as in item~\ref{it:final2} in Theorem~\ref{final}. We can assume without loss of generality that $G=\lift(\tilde{G},G^*,\data)$. Let $d$ be the largest fiber size of the groups $F_a=\pi_1(\data(a)):a\in \dom(G^*)$. Then by Theorem~\ref{final} we know that $G\in \subex_d$ which gives us the required upper bound for $\os_n(G)$.
	
	If $d=1$ then it follows from~\cite{simon2025omega}, Corollary 5.5 that $\os_n(G)=\os_n(\fm_*)$ for some stable structure $\fm_*$. By Theorem~\ref{sam} we know that in this case $\fm_*$ is monadically stable, and thus it is hereditarily cellular. Then $G_*\coloneqq \aut(\fm_*)\in \gone$, so in this case we are done.
	
	Now let us assume that $d\geq 2$. Then we show a lower bound as stated in item~\ref{it:medium}. Let us pick some $a\in \dom(G^*)$ such that for the fiber group $F_a=\pi_1(\data(a))$ has fiber size $d$, and let $C\coloneqq \pi^{-1}(a)$. We claim that $u\same_G v$ some tuples $u,v$ in $C$ then $\pi_1(u)\same_C \pi_1(v)$. Indeed, if $g(u)=v$ for some $G$ then $g$ fixes the whole fiber $C$, and thus $\mu_{\pi_1}(g|_C)\in F_a$. This implies immediately that $\os_n(G)\geq \os_n(F_a)$. Together with the ``moreover'' part of Lemma~\ref{when_subex} this gives us the required lower bound.
\end{proof}

	Note that Theorem~\ref{orbit_growth} immediately implies Theorem~\ref{gaps_intro} which confirms Conjecture 3 in~\cite{braunfeld2022monadic}.

\thmgaps*
	
	By following the proof of Theorem 5.1 in~\cite{simon2025omega} we can describe the group $G_*$ in the $d=1$ case in the proof of Theorem~\ref{orbit_growth} as follows. Let us assume that $G=\lift(\tilde{G},G^*,\data)$ as in item~\ref{it:final2} in Theorem~\ref{final}, let $Y^*\coloneqq \{a\in X^*: |\dom(F_a)|\neq 1\}$. Then clearly $Y^*$ is an invariant subset of $G^*$. Let $X_*\coloneqq \omega\times Y^*\cup \{0\}\times (X^*\setminus Y^*)$. Then $G_*$ can be described as the set of permutations $\sigma\in \sym(X^*)$ for which there exists a $\sigma^*\in G^*$ such that $\pi_1(\sigma(n,a))=\sigma^*(a)$. Note that in this case the second projection defines a covering map from $G_*$ to $G^*$. It is easy to see that in this case we have
\begin{align*}
\rk(G^*)=&\max(\rk(G^*|_{Y^*}),\rk(G^*|_{X^*\setminus Y^*})),\\
\rk(G_*)=&\max(\rk(G^*|_{Y^*})+1,\rk(G^*|_{X^*\setminus Y^*})).
\end{align*}

	In particular $\rk(G_*)\geq 1$ if and only if $\rk(G^*)=1$ and $Y^*$ is a subset of the union of finite orbits of $G^*$. Combining this with Theorem~\ref{cell} we can give a description of closed groups with at most polynomial unlabelled growth using our notation. This classification is also given in~\cite{falque2019classification} but our presentation is different.
	
\begin{theorem}\label{poly_full}
	Let $G$ be a permutation group acting on a countable set. Then the following are equivalent.
\begin{enumerate}
\item\label{it:poly1} $G$ is $P$-oligomorphic.
\item\label{it:poly2} $G$ is isomorphic to some group $\lift(\tilde{G},G^*,\data)$ where
\begin{itemize}
\item the triple $(\tilde{G},G^*,\data)$ is as in Definition~\ref{def:recover},
\item $(\tilde{G},G^*,\pi_2)$ is isomorphic to some finite cover in $\mathbf{A}$ with $\rk(G^*)\leq 1$,
\item for all $a\in X^*$ we have $\data(a)=(F_a,B_a,\phi_a)$ where either $|\dom(F_a)|=1$; or $a\in \acl_{G^*}(\emptyset)$ and one the following holds.
\begin{enumerate}[(i)]
\item $B_a=F_a=\aut(\mathbb{Q};R)$ where $R\in \{<,\betw,\cyc,\sep\}$,
\item $B_a=(\mathbb{Q};<)$ and $F_a=(\mathbb{Q};\betw)$,
\item $B_a=(\mathbb{Q};\cyc)$ and $F_a=(\mathbb{Q};\sep)$.
\end{enumerate}
\end{itemize}
\end{enumerate}
\end{theorem}

\begin{proof}
	Note that condition~\ref{it:poly1} implies that $G\in \subex_1$. Therefore, it is enough to show that a group $G=\lift(\tilde{G},G^*,\data)$ as described in item~\ref{it:final2} of Theorem~\ref{final} for $d=1$ satisfies the conditions specified in item~\ref{it:poly2} of this theorem if and only if $\os_n(G)$ has a polynomial upper bound. First, let us notice that by Lemma~\ref{normal_highly} the conditions described for a pair $(F_a,B_a): |\dom(F_a)|\neq 1$ in the $d=1$ case holds if and only if one of the conditions (i)-(iii) above holds.
	
	Since the largest fiber size of the groups $F_a$ is 1, the groups $G$ and $G_*$ have the same unlabelled growth profile where $G_*$ is the group as described above. By Theorem~\ref{poly} this holds if and only if $\rk(G_*)\leq 1$. As we have seen, this hold if and only if $\rk(G^*)\leq 1$ and for all $a\in \dom(G^*)$ which has an infinite orbit we have $|\dom(F_a)|=1$.
\end{proof}

\begin{remark}
	By taking a closer look we can give a more detailed description of the finite covers $(\tilde{G},G^*,\pi_2)$ showing up in item~\ref{it:poly2} of Theorem~\ref{poly_full}. Let $(\tilde{H},H^*,\pi_2)\in \mathbf{A}$ and $\tilde{N}_0,\dots,\tilde{N}_k,N_0^*,\dots,N_k^*$ be as in Definition~\ref{def:fcoverh}. (The case when $\dom(G^*)$ is finite can be covered by the case when $k=0$.) Then since $\rk(G^*)\leq 1$ it follows that all groups $N_i^*$ and $H^*$ and therefore also $\tilde{N}_i$ and $\tilde{H}_i$ must be finite. Moreover, if $u\in \dom(H^*)\setminus \dom(N_0^*)$ then the groups in $\data((u,0))$ must be trivial, therefore in this case we can assume that $|\pi_2^{-1}((u,0))|=|\pi_2^{-1}(u)|=1$. On the other hand, if $u\in \dom(N_0^*)$ then $|F_{(u,0)}:B_{(u,0)}|\leq 2$, so we can also assume that all fiber sizes $|\pi_2^{-1}((u,0))|=|\pi_2^{-1}(u)|$ are at most 2.
\end{remark}

\begin{example}
	Let us consider the groups $Z_2\subset \sym(\{0,1\})$ and $Z_4\subset \sym((0,0),(0,1),(1,0),(1,1))$ generated by the cycle $((0,0)\,(0,1)\,(1,0)\,(1,1))$. Then $(Z_4,Z_2,\pi_2)$ is a finite cover which corresponds to the homomorphism $\pi:Z_4\rightarrow Z_2, 1\mapsto 1$. Clearly, this cover is non-split. Now let $F\coloneqq (\mathbb{Q};\betw), B\coloneqq (\mathbb{Q};<)$, let $\phi$ be the only isomorphism from $B/F$ to $\sym(\{0,1\})\simeq Z_2$, and let $D$ be the constant map on $\{0,1\}$ with value $(B,F,\phi)$. Then by Theorem~\ref{poly_full} the group $G\coloneqq \lift(Z_4,Z_2,D)$ is $P$-oligomorphic. This group is isomorphic to the group generated by $\aut(\mathbb{Q};<)^2\leq \sym(\mathbb{Q})^2$ and some operation $\sigma$ (which corresponds to the generator of $Z_4$) that flips the two copies of $\mathbb{Q}$ and on one side it is an isomorphism and on the other side it is an anti-isomorphism (i.e., it reverses the order). Note that in this case $\sigma^2$ preserves the two copies of $\mathbb{Q}$ but it reverses the order on both.
\end{example}
	
	Next we show that if we are allowed to add finitely many constants then closed $P$-oligomorphic groups have a much simpler description.
	
\begin{theorem}\label{poly_const}
	Let $G$ be a permutation group acting on a countable set. Then the following are equivalent.
\begin{enumerate}
\item~\label{it:poly_c1} $G$ is $P$-oligomorphic and closed.
\item~\label{it:poly_c2} There exists a finite set $F\subset \dom(G)$ such that $G_F$ is isomorphic to a direct product of finitely many copies of $\aut(\mathbb{Q};<)$ and some cellular group.
\end{enumerate}
\end{theorem}

\begin{proof}
	The implication~\ref{it:poly_c2}$\rightarrow$~\ref{it:poly_c1} follows directly from Theorem~\ref{poly_full} considering that $P$-oligomorphic groups are closed under taking supergroups.
	
	Now let us assume that~\ref{it:poly_c1} holds. Then we can assume that $G=\lift(\tilde{G},G^*,\data)$ as in item~\ref{it:poly2} of Theorem~\ref{poly_full}. Let $X^*\coloneqq \dom(G^*)$, and let $Y^*$ be the union of finite orbits of $G^*$. Since $G^*$ is oligomorphic it follows that $Y^*$ is itself finite. Let $$Z^*\coloneqq \{a\in X^*: |\pi^{-1}(a)|\neq 1\}.$$ Then by the conditions listed in item~\ref{it:poly2} of Theorem~\ref{poly_full} we can conclude that $Z^*\subset Y^*$, and thus $Z^*$ is also finite, and for all $a\in Z^*$ we have
$$\aut(\mathbb{Q};R)\times \id(\{a\})\leq G_{((\pi^{-1}(a)))}\leq G_{(\pi^{-1}(a))}\leq (\aut(\mathbb{Q};R)\!\flip)\times \id(\{a\})$$ for $R\in \{<,\cyc\}$.

	Now let $H^*\coloneqq (G^*)_{Z^*}$ and let $H$ be the intersection of the setwise stabilizers of $\pi^{-1}(a):a\in Z^*$. Then $H$ is a cover of $H^*$. Let $s,t\in \mathbb{Q}$ with $s<t$, and let $F\coloneqq\{s,t\}\times Z^*$. Then $F\subset X\coloneqq \dom(G)$ is finite. For all $a\in Z^*$ we have the following inclusions:
$$(B_a)_{s,t}\times \id(\{a\})\leq (H_F)|_{\pi^{-1}(a)}=(H_{F\cup (X\setminus \pi^{-1}(a))})|_{\pi^{-1}(a)}\leq (F_a)_{s,t}\times \id(\{a\}).$$

	We have seen that $|F_a:B_a|\leq 2$ and $|F_a:B_a|=2$ is only possible if $\flip\,\in F_a$. Note, however, that $\flip$ cannot be contained in the stabilizer of two different elements of $(\mathbb{Q};<)$ therefore in either case $(B_a)_{s,t}=(F_a)_{s,t}$, and therefore the equality holds everywhere in the chain of inclusions above. This implies that the group $H_F$ can be written as the direct product of the groups $(B_a)_{s,t}\times \id(\{a\})$ and $(H_F)|_{X\setminus \pi^{-1}(Z^*)}$. By definition $\pi$ is one-to-one outside $\pi^{-1}(Z^*)$ so the latter group is isomorphic to $G^*|_{X^*\setminus Z^*}$ which is an invariant subset of a cellular group, so it is cellular. On the other hand, by Lemma~\ref{high_const} it follows that the groups $(B_a)|_{s,t}$ must be isomorphic to either $\aut(\mathbb{Q};<)^3\times \id(\{0,1\})$ or $\aut(\mathbb{Q};<)^2\times \id(\{0,1\})$. This implies that $H|_F$ is isomorphic to a direct product of finitely many copies of $\aut(\mathbb{Q};<)$ and some cellular group. Finally, let us observe that if a permutation in $G$ fixes an element of $\pi^{-1}(a)$ then it also fixes the whole class $\pi^{-1}(a)$. This implies that in fact $G_F=H_F$ which finishes the proof of the theorem.
\end{proof}

\begin{remark}
	It follows from an easy calculation that $G$ is $P$-oligomorphic then so is $G_F$ for every finite $F\subset \dom(G)$. 
\end{remark}

	As usual Theorem~\ref{poly_const} can also be formulated in terms of structures which gives us exactly Theorem~\ref{poly_const_struct}.

\polyconst*

	The relevance of Theorem~\ref{poly_const} (or Theorem~\ref{poly_const_struct}) in the context of CSPs is as follows. We say that a countable $\omega$-categorical $\fa$ is a \emph{model-complete} core if the closure of $\aut(\fa)$ (in the set of all functions $A\rightarrow A$) is equal to its endomorphism monoid. We say that two structures $\fa$ and $\fb$ are \emph{homomorphically equivalent} if there exist homomorphisms between them in both directions. Note that in this case $\csp(\fa)=\csp(\fb)$. By the results of~\cite{Cores-journal} we know that every $\omega$-categorical structure is homomorphically equivalent to a model-complete core which is unique up to isomorphism and is also $\omega$-categorical. We also know that this model-complete core always has at most as fast as the unlabelled growth profile of the original structure (see for instance~\cite{bodor2024classification}, Lemma 3.16). For this reason in order to solve the infinite-domain CSP dichotomy conjecture for structures with $P$-oligomorphic automorphism groups it is enough to solve them for the ones which are model-complete cores. We also know that for $\omega$-categorical model-complete cores adding constants does not change the complexity of its CSP up to $\mathbf{LOGSPACE}$ reductions~\cite{wonderland}. Combining the argument above with Theorem~\ref{poly_const} we obtain that every CSP over a structure with $P$-oligomorphic automorphism group is $\mathbf{LOGSPACE}$ equivalent to a CSP over a model-complete core whose automorphism group is a direct product of finitely many copies of $\aut(\mathbb{Q};<)$ and some cellular group.

%% file: thomas.tex

\section{Counting structures with not too fast unlabelled growth}\label{sect:thomas}

	In this section we show that Thomas' conjecture holds for structures with not too fast unlabelled growth, that is, every structure in $\subexs$ has finitely many reducts up to interdefinability; or equivalently every group in $\subex$ has finitely closed supergroups. The idea of our proof is similar to the one given given in Section 5 of~\cite{bodor2024classification} where the statement is shown for the special case of hereditarily cellular structures: we define some parameters such that they never increase when taking a reduct, and such that up to bidefinability there are only finitely many structures in $\subexs$ with any given finite upper bound on these parameters. As a byproduct of our argument we also obtain that $\subexs$ contains countably many structures up to bidefinability; or equivalently $\subex$ contains countably many groups up to isomorphism.

\begin{definition}
	Let $G$ be a permutation group. Then we write
\begin{itemize}
\item $\epsilon(G)$ for the number of congruences of $G$, and
\item $w(G)$ for the supremum of the set $$\{|\acl_{G_x/E}(\emptyset)|: x\in \dom(G), \text{$E$ is a congruence of $G_x$}\}.$$
\end{itemize}
\end{definition}

\begin{observation}\label{we_reduct}
	The following are clear from the definition.
\begin{itemize}
\item If $G$ is oligomorphic then $\epsilon(G)$ and $w(G)$ are finite.
\item If $E$ is a congruence of $G$ then $\epsilon(G/E)\leq \epsilon(G)$ and $w(G/E)\leq w(G)$.
\item If $G\leq H$ then $\epsilon(H)\leq \epsilon(G)$.
\end{itemize}
\end{observation}

	We know the following from~\cite{bodor2024classification}.
	
\begin{theorem}\label{thomas_stab}
	For all $n,m\in \omega$ there exist finitely many groups $G\in \gone_n$ with $w(G)\leq m$ up to isomorphism.
\end{theorem}

	Our plan in this section is to show an analogous statement for the whole class $\subex$, and thereby show that every structure in $\subexs$ has finitely many closed supergroups. By using the following observation it is enough to check the finiteness ``up to isomorphism''.
	
\begin{lemma}\label{iso_enough}
	Let $G$ be a closed oligomorphic group. Then the following are equivalent.
\begin{enumerate}
\item\label{it:fin_inter} $G$ has finitely many closed supergroups.
\item\label{it:fin_bi} $G$ has finitely many closed supergroups up to isomorphism.
\end{enumerate}
\end{lemma}

\begin{proof}
	We know that $G$ is the automorphism group of some $\omega$-categorical structure $\fa$. Then item~\ref{it:fin_inter} is equivalent to $\fa$ having finitely many reducts up to \emph{interdefinability} and item~\ref{it:fin_bi} is equivalent to $\fa$ having finitely many reducts up to \emph{bidefinability}. The equivalence of these conditions then follows from Proposition 6.37 in~\cite{bodirsky2021permutation}.  
\end{proof}
	
	We start by a couple of easy observations about the functions $w$ and $\epsilon$.
	
\begin{lemma}\label{we_go_down}
	Let $(G,G^*,\pi)$ be a cover. Then $\epsilon(G^*)\leq \epsilon(G)$. Furthermore, if all binding groups of $(G,G^*,\pi)$ are transitive then $w(G^*)\leq w(G)$.
\end{lemma}

\begin{proof}
	The first claim is clear since for every congruence $E$ of $G^*$ the relation $\pi^{-1}(E)$ is a congruence of $G$. For the other claim, let us pick an $x\in X^*\coloneqq \dom(G^*)$, and a congruence of $E^*$ such that $G^*_x/E^*$ has $w(G^*)$ many finite orbits. Let $A\coloneqq \pi^{-1}(x)$ and pick an element $a\in A$. We claim that $\mu_{\pi}(G_a)=\mu_{\pi}(G_{\{A\}})$. Clearly, we have $G_a\subset G_{\{A\}}$, and thus $\mu_{\pi}(G_a)\subset \mu_{\pi}(G_{\{A\}})$. For the other containment, let us assume that $g\in G_{\{A\}}$. Then by our assumption we know that there exists a $h\in \ker(\mu_{\pi})$ such that $h(g(a))=a$. Then $hg\in G_a$ and $\mu_{\pi}(hg)=\mu_{\pi}(g)$. Now let $E$ be the inverse image if $E^*$ on $X\coloneqq \dom(G)$. Then clearly $E$ is a congruence of $G_{(A)}$ and thus it is also a congruence of $G_a\subset G_{(A)}$. Since $\mu_{\pi}(G_a)\subset \mu_{\pi}(G_{\{A\}})$ and $\sim_{\pi}\subset E$ we obtain that $G_a/E\simeq G_{\{A\}}/E\simeq G_x/E^*$. In particular, $|\acl_{G_a/E}(\emptyset)|=|\acl_{G_x/E^*}(\emptyset)|=w(G^*)$, and thus $w(G)\geq w(G^*)$.
\end{proof}

\begin{proposition}\label{ep_rk}
	Let $G\in \gone$. Then $\epsilon(G)\geq \rk(G)$.
\end{proposition}

\begin{proof}
	We show this by induction on $\rk(G)$. For $\rk(G)=0$ there is nothing to show. Otherwise, $\rk(G)$ can be written as $\biggroup(H;N_0,\dots,N_k)$ with $\rk(H)=\rk(G)-1$. Thus, by the induction hypothesis we have $\epsilon(H)\geq \rk(G)-1$. Let us observe that for all congruences $E$ of $H$ the relation $E'=\{((a,n),(b,n): (a,b)\in E\}$ is a congruence of $G$. Since none of these congruences include the full congruence $\dom(G)\times \dom(G)$ we have $\epsilon(G)\geq \epsilon(H)+1\geq \rk(G)$.
\end{proof}

	Next we show that for a fixed group $G^*\in \gone$ and a $d\in \mathbb{N}$ there exist only finitely many non-isomorphic groups $G=\lift(\tilde{G},G^*,\data)$ where the triple $(\tilde{G},G^*,\data)$ is as in item~\ref{it:final2} of Theorem~\ref{final}.
	
	Note that by Remark~\ref{at_most5} we can always choose the group $\tilde{G}$ such that all of its fibers have size at most 5. We first show that under these constraints there are finitely many options for $\tilde{G}$ up to isomorphism.
		
\begin{lemma}\label{w_fincover}
	Let $n\in \omega$. Then every hereditarily cellular group has only finitely many covers with maximal fiber size at most $n$ up to isomorphism.
\end{lemma}

\begin{proof}
	We prove the lemma by induction on the rank. For finite-degree groups there is nothing to prove. Now let us fix some group $G^*\in \gone$ of rank at least 1, and let $(G,G^*,\pi)$ be a finite cover of $G^*$. By Lemma \ref{ms_covers_class} we can assume without loss of generality that $(G,G^*,\pi)\in \mathbf{A}$. This means that we can find some groups $H,\dots,N_0,N_1,\dots,N_k,H^*,N_0^*,\dots,N_k^*$ and a finite cover $(H,H^*,\rho)\in \mathbf{A}$ as in Definition~\ref{def:fcoverh}. By the definition of the rank we know that ranks of all groups $H^*,N_0^*,\dots,N_k^*$ are at most $\rk(G^*)$. Clearly, all the fibers of all the covers $(H,H^*,\rho)$ and $(N_i,N_i^*,\rho|_{\dom(N_i)})$ have size at most $n$. Thus, by the induction hypothesis there are only finitely many options for the cover $(H,H^*,\rho)$ and the groups $N_0,\dots,N_k$ up to isomorphism. Since these groups determine the cover $(G,G^*,\pi)$ we obtain that $G^*$ has only finitely many finite covers with maximal fiber size at most $n$.
\end{proof}

\begin{lemma}\label{fixed_star}
	Let $G^*\in \gone$, and $d\in \mathbb{N}$. Then up to isomorphism there exist finitely many groups $G\in \subex_d$ which are isomorphic to some group $\lift(\tilde{G},G^*,\data)$ where $(\tilde{G},G^*,\data)$ is as in item~\ref{it:final2} of Theorem~\ref{final}.
\end{lemma}

\begin{proof}
	As we argued above, we can assume that the maximal fiber size of the cover $(\tilde{G},G^*,\pi_2)$ is at most 5. Since $G^*$ is oligomorphic, and $\data$ must be constant on each orbit it is enough to show that $\data(a)=(F_a,B_a,\phi_a)$ can only take finitely many values for a given $a$ (up to renaming the elements of $\dom(F_a)$). By Theorem~\ref{fs_normal} we know that there are only finitely many possibilities for the pair $(F_a,B_a)$, and since $|F_a:B_a|$ is finite there are also only finitely many options for the homomorphism $\phi_a$.
\end{proof}

\begin{lemma}\label{only_finite}
	Let $d,n,m\in \omega, d\geq 1$. Then there exist finitely many groups $G\in \subex_d$ with $\epsilon(G)\leq n$ and $w(G)\leq m$ up to isomorphism.
\end{lemma}

\begin{proof}
	Suppose that $G$ satisfies the conditions of the lemma. Let $(\tilde{G},G^*,\data)$ be as in item~\ref{it:final2} of Theorem~\ref{final}. Since $G$ is a cover of $G^*$ with transitive binding groups we have $\epsilon(G^*)\leq \epsilon(G)\leq n$ and $w(G^*)\leq w(G)\leq m$. Combining this with Proposition~\ref{ep_rk} and Theorem~\ref{thomas_stab} we obtain that there are finitely many possibilities for the group $G^*$ up to isomorphism. Then the statement of the lemma follows from Corollary~\ref{fixed_star}.
\end{proof}

	The following is a direct consequence of Lemma~\ref{only_finite}.

\begin{corollary}\label{countable_many}
	$\subex$ contain countable many groups up to isomorphism. Equivalently $\subexs$ contains countably many structures up bidefinability.
\end{corollary}

	Now we are ready to prove Thomas' conjecture for the class $\subexs$.

\begin{theorem}\label{thomas1}
	Every group in $\subex$ has finitely many closed supergroups. Equivalently, every structure in $\subexs$ has finitely many reducts up to interdefinability.
\end{theorem}

\begin{proof}
	Let $G\in \subex$. Then by Corollary~\ref{union} we know that $G\in \subex_d$ for some $d\in \mathbb{N}$. Now let $G^+$ be a closed supergroup of $G$. Then $G^+$ is also contained in $\subex_d$. Moreover, we have $\epsilon(G^+)\leq \epsilon(G)$ and $w(G^+)\leq w(G)$. Therefore, by Lemma~\ref{only_finite} there are finitely many possibilities for the group $G^+$ up to isomorphism. The statement of the lemma then follows from Lemma~\ref{iso_enough}.
\end{proof}

	Note that since every structure in $\subexs$ is finitely homogenizable (see Corollary~\ref{fin_hom}), it follows that all these structures indeed fall into the scope of Thomas' conjecture.
	
\section{Interpretation of structures with not too fast unlabelled growth}

	In this section we show that all structures in $\subexs$ are interpretable in $(\mathbb{Q};<)$. Our strategy is to use Theorem~\ref{final_struct} by showing that all constructions listed in the theorem preserve interpretability in $(\mathbb{Q};<)$.
	
\begin{definition}
	Let $\fa$ and $\fb$ be structures. We say that $I$ is an \emph{interpretation} of $\fb$ in $\fa$ if there exists some $d\in \omega$ such that $I$ is a partial surjective map from $A^d$ to $B$ such that for every atomic formula $\varphi(x_1,\dots,x_k)$ in the signature of $\fb$ the relation
\[
I^{-1}(R)\coloneqq \{(a_1^1,\dots,a_1^d,\dots,a_k^1,\dots,a_k^d): \fb\models \varphi(I(a_1^1,\dots,a_1^d),\dots,I(a_k^1,\dots,a_k^d))\}
\]
	is definable in $\fa$.
	
		
	We say that $\fa$ \emph{interprets} $\fb$, or $\fb$ is \emph{interpretable} in $\fa$, if there exists an interpretation of $\fb$ in $\fa$. For a fixed structure $\fa$ we write $\inter(\fa)$ for the class of structures which are interpretable in $\fa$.
\end{definition}

	Let us consider some interpretation $I$ as in the definition above. By applying the definition for the formula $(x=x)$ we see that the domain set of $I$ needs to be definable. Moreover, by applying the definition for the formula $(x=y)$ we obtain that the kernel of the map $I: \dom(I)\rightarrow B$ also needs to be definable in $\fa$. This means essentially that the domain set of $\fb$ can be thought of as a definable quotient of some definable substructure of some power of $\fa$ and every relations of $\fb$ needs to be represented by some set of tuples which are also definable in $\fa$.

	We start with some easy observations about interpretability in a given structure.

\begin{lemma}\label{inter:lab}
	Let $\fa,\fb,\fc$ be relational structures such that $\fa$ has at least two elements, $\fb,\fc\in \inter(\fa)$ and $\fb$ and $\fc$ have disjoint signatures. Then $\fb\uplus \fc$ and $\fb\wr\fc$ are interpretable in $\fa$.
\end{lemma}

\begin{proof}
	Let $I:A^d\rightarrow B$ and $J:A^e\rightarrow C$ be interpretations of $\fb$ and $\fc$, respectively. We can assume without loss of generality that $d=e$ and that $B$ and $C$ are disjoint.
	
	Let us define the partial map
$$
I+J: A^{d+2}\rightarrow B\cup C, (a_1,\dots,a_{d+2})\mapsto
\begin{cases*}
I(a_1,\dots,a_d) \text{ if $a_{d+1}=a_{d+2}$}\\
I(a_1,\dots,a_d) \text{ if $a_{d+1}\neq a_{d+2}$}
\end{cases*} 
$$

	wherever it is defined. Then it is easy to see that $I+J$ is an interpretation of $\fb\uplus \fc$ in $\fa$.
	
	For the other claim we simply consider the partial map $$I\times J: A^{2d}\rightarrow B\times C, (a_1,\dots,a_{2d})\mapsto (I(a_1,\dots,a_d),J(a_{d+1},\dots,a_{2d}))$$ wherever it is defined. Then it can be checked easily that $I\times J$ is an interpretation of $\fb\wr \fc$ in $\fa$.
\end{proof}

\begin{corollary}\label{fcover_int}
	Let us assume that $\fa$ has at least two elements, and $\fb\in \inter(\fa)$. Then every strongly split finite cover of $\fb$ is interpretable in $\fa$.
\end{corollary}

\begin{proof}
	Clearly, every structure interprets the 1-element pure set as witnessed by the constant map. It is also clear that $\inter(\fa)$ is closed under taking reducts, therefore it is enough to show that every strongly trivial finite cover of $\fb$ is interpretable in $\fa$, i.e., for all finite subsets $(F;=)\wr \fb\in \inter(\fa)$. This follows directly from Lemma~\ref{inter:lab}.
\end{proof}

\begin{theorem}\label{thm:inter}
	$\subexs\subset \inter(\mathbb{Q};<)$ and $\gones\subset \inter(\omega;=)$.
\end{theorem}

\begin{proof}
	By Theorem~\ref{fin_covers_q} we know that all finite covers of reducts of $(\mathbb{Q};<)$ are reducts of some strongly split finite cover of $(\mathbb{Q};<)$. Therefore, by Corollary~\ref{fcover_int} all these structures are interpretable in $(\mathbb{Q};<)$. Then the statement of the theorem follows from Theorems~\ref{final_struct}, \ref{thm:build_structures} and Lemma~\ref{inter:lab}.
\end{proof}